\documentclass[leqno]{amsart}
\usepackage{egastyle}
\usepackage{mymacros}
\usepackage[all,tips,dvips]{xy}
\usepackage{url}
\usepackage{bm}

\newcommand\A\bA
\newcommand\B\bB
\newcommand\D\bD
\newcommand\G\bG

\DeclareMathOperator\conv{conv}
\DeclareMathOperator\face{face}
\DeclareMathOperator\relint{relint}
\DeclareMathOperator\vertices{vert}
\DeclareMathOperator\minset{minset}
\DeclareMathOperator\maxset{maxset}
\DeclareMathOperator\Trop{Trop}
\DeclareMathOperator\New{New}
\DeclareMathOperator\trop{trop}
\DeclareMathOperator\val{val}
\DeclareMathOperator\inn{in}
\DeclareMathOperator\vol{vol}
\DeclareMathOperator\MV{MV}
\DeclareMathOperator\NP{NP}
\renewcommand\ord\val
\newcommand\sep{{\operatorname{sep}}}
\newcommand\rig{{\operatorname{rig}}}
\newcommand\an{{\operatorname{an}}}
\newcommand\berk{{\operatorname{berk}}}

\title{Tropical analytic geometry, Newton polygons, and tropical intersections}
\author{Joseph Rabinoff}

\begin{document}

\begin{abstract}
  In this paper we use the connections between tropical algebraic
  geometry and   rigid analytic geometry in order to prove two main
  results.  We use tropical methods to prove a theorem about the Newton
  polygon for convergent power series in several variables: if
  $f_1,\ldots,f_n$ are $n$ convergent power series in $n$ variables with
  coefficients in a non-Archimedean field $K$, we give a formula
  for  the valuations and multiplicities of the common zeros of
  $f_1,\ldots,f_n$.  We use rigid-analytic methods to show that stable
  complete intersections of tropical hypersurfaces compute algebraic
  multiplicities even when the intersection is not tropically proper.
  These results are naturally formulated and proved using the theory of
  tropicalizations of rigid-analytic spaces, as introduced by
  Einsiedler-Kapranov-Lind~\cite{ekl:non_arch_amoebas} and  
  Gubler~\cite{gubler:tropical}.
  We have written this paper  to be as 
  readable as possible both to tropical and arithmetic geometers.
\end{abstract}

\maketitle

\section{Introduction}

\paragraph Strong \label{par:firstpar}
connections between tropical algebraic geometry and the theory of
rigid-analytic spaces allows one to prove theorems in one field using
ideas from the other.  This paper establishes
two main results, the first result rigid-analytic in nature and the
second tropical:  
\begin{enum}
\item \textbf{A higher-dimensional theorem of the Newton polygon.} Let $K$ be a
  field that is complete with respect to a 
  nontrivial non-Archimedean valuation $\val: K^\times\to\R$ and let
  $f_1,\ldots,f_n$ be $n$ convergent power series (in a sense to be made
  precise later) in $n$ variables with coefficients in $K$.  Given $v\in\R^n$, we will
  give a formula~\parref{thm:multiplicity.ps} for the 
  number of common zeros (counted with multiplicity) 
  $\xi = (\xi_1,\ldots,\xi_n)\in(\bar K{}^\times)^n$ of $f_1,\ldots,f_n$
  such that $v = \trop(\xi)\coloneq(\val(\xi_1),\ldots,\val(\xi_n))$, in
  terms of the valuations of the coefficients of the $f_i$.
  (The set of all $v$ such that $v=\trop(\xi)$ for some common zero 
  $\xi\in V(f_1)\cap\cdots\cap V(f_n)(\bar K)$ is the
  tropicalization of $V(f_1)\cap\cdots\cap V(f_n)$, and can also be
  effectively calculated.)  This theorem generalizes the 
  classical theorem of the Newton polygon, which gives the valuations
  and multiplicities of the zeros of a convergent power series in
  one variable; see~\parref{eg:multiplicity.ps.dim1}.

\item \textbf{Tropically non-proper stable intersection multiplicities calculate
    algebraic multiplicities.}
  Let $f_1,\ldots,f_n\in K[x_1^{\pm 1},\ldots, x_n^{\pm 1}]$ be nonzero
  Laurent polynomials and let $C$ be a connected component of 
  $\bigcap_{i=1}^n \Trop(f_i)$.  We will
  show~\parref{thm:non.proper.multiplicity} that the sum of the stable 
  intersection multiplicities of the points of $C$ is equal to the
  sum of the algebraic multiplicities of the common zeros $\xi$ of
  $f_1,\ldots,f_n$ such that $\trop(\xi)\in C$ (assuming that the latter
  is finite), after passing to a suitable toric variety if $C$ is
  unbounded. 
\end{enum}
The above results are naturally formulated and proved in the realm of
\emph{tropical analytic geometry}, the theory of tropicalizations of
rigid-analytic spaces, as introduced by
Einsiedler-Kapranov-Lind~\cite{ekl:non_arch_amoebas} and 
Gubler~\cite{gubler:tropical}.  Much of this paper is dedicated to
extending their results and enriching this theory in several ways.

\paragraph 
Let us discuss~(\ref{par:firstpar},i) in more detail.
Let $K$ be as above, and assume for simplicity that $K=\bar K$ and
that $\val(K^\times)=\R$.  We (provisionally) define the rigid-analytic
unit ball $\B^n$ to be the set of all points $(\xi_1,\ldots,\xi_n)\in K^n$
such that $\val(\xi_i)\geq 0$ for all $i$.  Define a map
$\trop:\B^n\to(\R_{\geq 0}\cup\{\infty\})^n$ by
$\trop(\xi_1,\ldots,\xi_n) = (\val(\xi_1),\ldots,\val(\xi_n))$.  Let
$f = \sum_\nu a_\nu x^\nu\in K\ps{x_1,\ldots,x_n}$ be a function
converging on $\B^n$ (i.e.\ such that $\val(a_\nu)\to\infty$) and
define $\Trop(f)\subset(\R_{\geq 0}\cup\{\infty\})^n$ to be the set
$\{\trop(\xi)~:~f(\xi)=0\}$.  It is a fundamental fact
(see~\S\ref{sec:hypersurfaces}) that $\Trop(f)\cap\R_{\geq 0}^n$ is union of
finitely many polyhedra (generically of dimension $n-1$), and that the
ultrametric triangle inequality 
as applied to the equation $\sum_\nu a_\nu\xi^\nu = 0$ completely determines the set
$\Trop(f)$ (see~\ref{thm:main.tropical.thm.affinoid}).  In particular
the set $\Trop(f)$ is generally not hard to calculate.

\paragraph  Now \label{par:describe.mult.formula}
let $f_1,\ldots,f_n\in K\ps{x_1,\ldots,x_n}$ be power series
converging on $\B^n$.  If $\xi\in\B^n$ is a common zero of
$f_1,\ldots,f_n$ then $\trop(\xi)$ is contained in
$\bigcap_{i=1}\Trop(f_i)$, which is generically a finite set of points.
In other words, one gets very strong
restrictions on the valuations of the coordinates of the common zeros of
$f_1,\ldots,f_n$ via a simple combinatorial calculation, which when
$n=1$ reduces to finding the slopes of the Newton polygon of a power
series (see~\ref{eg:newtonpoly1}).  

The tropical hypersurfaces $\Trop(f_i)$ come equipped with
multiplicity information (the Newton complex), also determined by
the valuations of the 
coefficients of the $f_i$, which induces a notion of multiplicity on the
points of $\bigcap_{i=1}^n\Trop(f_i)$.  (When $n=1$ these multiplicities
amount to the horizontal lengths of the line segments in the Newton polygon.)
Osserman and Payne~\cite{osserman_payne:lifting} have proved a very general result
relating the multiplicities in the intersection theory of subvarieties of
a torus with the multiplicities in the intersection theory of tropical varieties,
which when applied to this case gives a formula for the number of common
zeros $\xi$ of an $n$-tuple of Laurent polynomials $f_1,\ldots,f_n$ (counted with
multiplicity) such that $\trop(\xi)$ is a specified point in
$\bigcap_{i=1}^n\Trop(f_i)$.  With enough of the framework of tropical
analytic geometry in place
(see~\S\S\ref{sec:polyhedral.subdomains}, \ref{sec:tropicalizations},
and~\ref{sec:hypersurfaces}), a continuity of
roots argument~\parref{thm:controots2} allows us to formulate and deduce
the corresponding result for power series~\parref{thm:multiplicity.ps}.  

\paragraph
From the perspective of a tropical geometer, the theory of rigid
spaces is useful because the analytic topology on $\R^n$ is much better
approximated by the rigid-analytic topology on the torus $\G_m^n$. For
example, the unit box 
$[0,1]^n$ is an analytic neighborhood in the Euclidean space $\R^n$, yet
$\trop\inv([0,1]^n)\subset|\G_m^n|$ is the $n$-fold product of the annulus
$\{\xi\in K^\times~:~\val(\xi)\in[0,1]\}$, which is a very
nicely behaved rigid-analytic object (it is a smooth affinoid space), but
is not the set of points underlying a subscheme.  
Similarly, $(\R\cup\{\infty\})^n$ can naturally be regarded as the
tropicalization of the affine space $\A^n$
(see~\S\ref{sec:extended.trop}), under which identification the unit ball
$\B^n$ is the inverse 
image under $\trop$ of $(\R_{\geq 0}\cup\{\infty\})^n$ (a
neighborhood of the point $(\infty,\ldots,\infty)$).

\paragraph
The following example is an application of
rigid-analytic methods to a tropically-local problem.
Let $U_{\{0\}} = \trop\inv(\{0\})\subset|\G_m^n|$.
This is an affinoid space, which implies
(see~\S\ref{sec:affinoids}) that it is the maximal spectrum of
the algebra $K\angles{U_{\{0\}}}$ of all
infinite Laurent series $\{\sum_{\nu\in\Z^n} a_\nu x^\nu~:~a_\nu\in K\}$
such that $\val(a_\nu)\to\infty$ as $|\nu_1|+\cdots+|\nu_n|\to\infty$.  If
$\fa\subset K[x_1^{\pm 1},\ldots,x_n^{\pm 1}]$ is an ideal and
$Y = V(\fa)\subset\G_m^n$ is the associated subscheme then 
$|Y|\cap U_{\{0\}}$ is identified with the set of maximal ideals of
$K\angles{U_{\{0\}}}$ containing the ideal $\fa K\angles{U_{\{0\}}}$, so
to show that $0\in\Trop(Y)$ is equivalent to showing that $\fa$ does not
generate the unit ideal in $K\angles{U_{\{0\}}}$.  This ends up being
equivalent to the well-known criterion that the initial ideal of $\fa$ at
$0$ contain no monomials.  The characterization of the tropicalization
(or rather the Bieri-Groves set) of a scheme by initial ideals was proved
by Einsiedler-Kapranov-Lind~\cite{ekl:non_arch_amoebas} using these
methods; we give a treatment below~\parref{thm:main.tropical.thm.affinoid}
which also applies to tropicalizations of analytic spaces.  (The first
complete proof of this theorem was given by
Draisma~\cite[Theorem~4.2]{draisma:tropical} and also uses affinoid
algebras, albeit in a different way; see~\ref{rem:trop.surj}.)

\paragraph 
A family of translations of a tropical variety parameterized by an
interval corresponds to a family of subvarieties of a torus
parametrized by a rigid-analytic annulus.  We study such families 
in order to prove the theorem indicated in~(\ref{par:firstpar},ii); to
illustrate the main idea we will sketch a special case.
Let $f_1,f_2\in K[x_1^{\pm 1},x_2^{\pm 1}]$, and suppose
that $\Trop(f_1)\cap\Trop(f_2)$ has a connected component $C$ of
positive dimension.  In this case there is a notion of the 
\emph{stable intersection multiplicity} of $\Trop(f_1)$ and $\Trop(f_2)$
along the component $C$, which is defined by translating $\Trop(f_2)$
in a generic direction by a small amount $\epsilon$ so that 
$\Trop(f_1)\cap\Trop(f_2)$ is a finite set of points, then
taking the limit as $\epsilon$ approaches zero.  This corresponds to
replacing $f_2(x_1,x_2)$ by $f_2(t^{a_1} x_1, t^{a_2} x_2)$
for generic $a_1,a_2\in\Z$ and some $t\in K^\times$ and then taking the
limit as $t$ approaches $1$. 
The results relating tropical and algebraic intersection multiplicities
mentioned above allow us to count the number of common zeros of
$f_1(x_1,x_2), f_2(t^{a_1} x_1, t^{a_2} x_2)$ with fixed tropicalization
when $\val(t) > 0$.

In order to relate these quantities with the number of common zeros $\xi$
of $f_1,f_2$ such that $\trop(\xi)\in C$ one is led to consider the family
of schemes $Y_t = V(f_1(x_1,x_2))\cap V(f_2(t^{a_1} x_1, t^{a_2} x_2))$
parametrized by the rigid-analytic annulus 
$\{t\in K^\times~:~\val(t)\in[0,\epsilon]\}$.  Under appropriate 
hypotheses the family $Y_t$ is automatically finite and
flat (at least after passing to an appropriate
toric compactification; see~\ref{thm:controots1}).
In particular, the length of the
fiber $Y_t$ is independent of $t$, which shows that algebraic intersection
multiplicities can be calculated after an analytically small perturbation.
We will make this kind of argument precise in~\S\ref{sec:non.transverse}.

\paragraph We now describe in more detail the contents of this paper.  As
the material in this paper bridges two different fields, we have made an effort
to ensure that it be as readable as possible
both to tropical geometers (who may not be familiar with
affinoid algebras or rigid spaces) as well as to arithmetic geometers
(who may not be familiar with convex or tropical geometry).
Hence we have included \S\S\ref{sec:convex.basic}--\ref{sec:extended.trop}
which are mainly expository, containing many examples 
and pictures.  In \S\ref{sec:convex.basic} we give definitions,
basic properties, and pictures of the convex-geometric objects that we
will encounter.  In \S\ref{sec:compactification} we describe
the compactification $N_\R(\Delta)$ of Euclidean space $N_\R$ associated
to a fan $\Delta$, as introduced by Kajiwara and Payne, which serves as the
tropicalization of the toric variety $X(\Delta)$.  We also introduce the
notion of a compactified polyhedron inside a space $N_\R(\Delta)$, which
will serve as the tropicalization of a so-called polyhedral affinoid subdomain of
$X(\Delta)$.  (The reader who is not familiar with toric varieties will
lose little on first reading by assuming throughout that $X(\Delta)$ is a torus and
$N_\R(\Delta)$ is Euclidean space.)  In \S\ref{sec:affinoids} we define,
give examples of, and state the basic properties of affinoid algebras and
rigid-analytic spaces.  We will emphasize the analogy with the theory of
finite-type schemes over a field.  In \S\ref{sec:extended.trop} we review
Kajiwara and Payne's notion of extended tropicalizations, in the process
defining the tropicalization map and setting our notation for toric varieties.

In \S\ref{sec:polyhedral.subdomains} we introduce the fundamental notion
of a polyhedral subdomain of a toric variety.  We will
show~\parref{prop:polyhedral.affinoids} that if 
$X(\Delta)$ is a toric variety adapted to a polyhedron $P\subset N_\R$ in
an appropriate sense, then the inverse image $U_P$ of the closure of $P$ in the
compactification $N_\R(\Delta)$ is the affinoid space associated to an
explicitly identified affinoid algebra.  This
extends the notion of a polytopal subdomain as defined
in~\cite{ekl:non_arch_amoebas} and~\cite{gubler:tropical} in a nontrivial
way.  In \S\ref{sec:tropicalizations} we define the tropicalization
$\Trop(Y)$ of a closed analytic subspace $Y$ of $U_P$ and characterize it
in terms of initial
ideals (\ref{thm:main.tropical.thm.affinoid},~\ref{thm:main.tropical.thm}).
The definition of $\Trop(Y)$ agrees with  
Gubler's notion when $P$ is a polytope.  In \S\ref{sec:hypersurfaces} we
review the canonical polyhedral complex structure on the tropical
hypersurface $\Trop(f)$ associated to a nonzero Laurent polynomial $f$, as
well as introducing the Newton complex $\New(f)$.  We will 
prove an important finiteness result~\parref{lem:inu.finiteness} which
implies that the tropicalization of an analytic hypersurface in a
polyhedral subdomain $U_P$ carries a \emph{finite} polyhedral structure.  

In \S\S\ref{sec:controots1}--\ref{sec:controots2} we prove two
``continuity of roots'' results which will be useful in 
\S\S\ref{sec:newton.polygons}--\ref{sec:non.transverse}.
Theorem~\parref{thm:controots1} is a tropical criterion for a 
rigid-analytic family of subvarieties (or analytic subspaces) of a toric
variety to be finite and flat.  Theorem~\parref{thm:controots2} is a local
continuity of roots criterion: it says that if $f_{1,t},\ldots,f_{n,t}$ is
a one-parameter family of power series in $n$ variables such that the
specializations $f_{1,0},\ldots,f_{n,0}$ have finitely many common zeros, then 
$f_{1,t},\ldots,f_{n,t}$ has the same number of common zeros when $|t|$ is 
small.  This result rests on Raynaud's approach to rigid geometry via
formal schemes.

In \S\ref{sec:newton.polygons} we prove a rigid-analytic intersection
multiplicity formula extending the corresponding result for subschemes of
a torus, as described in~\parref{par:describe.mult.formula}.  This result
is a strict generalization of the theorem of the 
Newton polygon that applies to convergent power series in several
variables.  More specifically, 
$f_1,\ldots,f_n$ are analytic functions on a polyhedral subdomain $U_P$ and
$v\in\bigcap_{i=1}^n\Trop(f_i)$ is an isolated point contained in the
interior of $P$ then we will give an explicit
formula~\parref{thm:multiplicity.ps} for the number of
common zeros $\xi$ of $f_1,\ldots,f_n$ such that $\trop(\xi) = v$.  

In \S\ref{sec:non.transverse} we prove a result relating algebraic
multiplicities and stable intersection multiplicities along a tropically
non-proper complete intersection of hypersurfaces.  That is,
$f_1,\ldots,f_n$ are nonzero Laurent polynomials and 
$C\subset\bigcap_{i=1}^n\Trop(f_i)$ is a connected component of positive
dimension then we will use the intersection multiplicity formula of 
\S\ref{sec:newton.polygons} to calculate the number of common zeros 
of $f_1,\ldots,f_n$ in an appropriate toric variety $X(\Delta)$ which lie
over the closure of $C$ in $N_\R(\Delta)$, in terms of stable tropical
intersection multiplicities.  The proof will involve families of
translations of tropical varieties parametrized by a rigid-analytic base, as
indicated above.

\paragraph[Tropical analytic geometry in the literature]  Several papers have
already appeared which take advantage of the connections between tropical
and rigid-analytic geometry. As mention above, Einsiedler, Kapranov,
and Lind~\cite{ekl:non_arch_amoebas} characterize the Bieri-Groves set of
a subvariety of a torus in terms of initial ideals; they also prove its 
connectedness using rigid-analytic results of Conrad~\cite{conrad:irredcomps}.
Payne~\cite{payne:analytification} has proved that the analytic space (in
the sense of Berkovich) associated to a subvariety of a toric variety is
naturally homeomorphic to the inverse limit of all of its
tropicalizations.  Gubler~\cite{gubler:tropical} has used the
combinatorial tructure on the tropicalization of a closed subspace of a
polytopal subdomain in order to prove special cases of
the Bogomolov conjecture over function fields~\cite{gubler:bogomolov}.  The author
has studied the tropicalization of the logarithm of a $p$-divisible
formal group in order to show that it has a canonical subgroup if its Hasse
invariant is small enough~\cite{jdr:thesis}.  

\paragraph[Acknowledgments]
It is a pleasure to thank Sam Payne for his input throughout the process
of writing this paper. 
The author thanks his advisor Brian Conrad for his comments
and suggestions, especially concerning the material in
\S\S\ref{sec:controots1}--\ref{sec:controots2}, as well as  his  advisor
Ravi Vakil.  The author is grateful to Walter Gubler and Matt Baker for
very interesting and relevant conversations and for their comments.

\tableofcontents

\paragraph[Notation] We will use the following general notation throughout
this paper.

\begin{bullets}
\item Let $S$ be a set and $T\subset S$ a subset, and let $f:S\to\R$ be any
function.  Define
\[ \minset(f,T) \coloneq \big\{ t\in T~:~ f(t) = \inf_{t'\in T} f(t') \big\}
\sptxt{ and }
\maxset(f,T) \coloneq \big\{ t\in T~:~ f(t) = \sup_{t'\in T} f(t') \big\}. \]
(These sets may be empty.)

\item If $X$ is a scheme (resp.\ a rigid space;
cf.~\S\ref{sec:affinoids}) we use $|X|$ to denote the set of closed points
of $X$ (resp. the underlying set of $X$).  For $\xi\in|X|$ we let
$\kappa(\xi) = \sO_{X,\xi}/\fm_{X,\xi}$ denote the residue field at $\xi$.

\item If $Y$ is a topological space and $P\subset Y$ is a subspace we let
  $P^\circ = P\setminus\del P$ denote the interior of $P$ in $Y$.

\item If $\Gamma$ is a subset of $\R$ and $r\in\Gamma$ we set
  \begin{align*} \Gamma_{\geq r} &= \{ v\in\Gamma~:~v\geq r \}, & 
  \Gamma_{>r} &= \{ v\in\Gamma~:~v > r \},  \\
  \Gamma_{\leq r} &= \{ v\in\Gamma~:~v\leq r \}, &
  \Gamma_{<r} &= \{ v\in\Gamma~:~v < r \}. 
  \end{align*}

\end{bullets}

\section{Basic notions from convex geometry}
\label{sec:convex.basic}

\paragraph
The tropicalization of an algebro-geometric or analytic-geometric object
is a convex-geometric object, which is combinatorial in nature.  In this
section we give definitions of, state some properties of, and draw some
pictures of the convex-geometric objects that will appear, for the benefit
of the reader who is not familiar with them.
Most of this material can be found in~\cite[\S\S1.2,1.5]{fulton:toric_varieties} 
and~\cite[Chapter~VI]{barvinok:convexity}, although almost all of it is
quite easy and instructive to prove on one's own.

Convex bodies live inside Euclidean space $\R^n$.  We prefer not to
choose a basis, so we fix the following notation for the rest of this paper:

\begin{notn} \label{notn:convex}

\begin{tabular}{ll}
  $N_\R\cong\R^n$ & is a real vector space of dimension $n$ \\
  $M_\R = N_\R^*$ & is its linear dual \\
  $\angles{\cdot,\cdot}: M_\R\times N_\R\to\R$ & is the canonical pairing \\
  $N\cong\Z^n$ & is a full-rank lattice in $N_\R$ \\
  $M = \Hom_\Z(N,\Z)$ & is the dual lattice in $M_\R$ \\
  $\Gamma \subset\R$ & is a nonzero additive subgroup \\
  $N_\Gamma = N\tensor_\Z \Gamma$ & the subgroup of $\Gamma$-rational
  points of $N_\R$ \\
  $M_\Gamma = M\tensor_\Z\Gamma$ & likewise for $M_\R$.
\end{tabular}\\
\end{notn}

The lattice $N\subset N_\R$ is called an \emph{integral structure}.  In
the sequel we will take the subgroup $\Gamma$ to be the value group of a
field equipped with a nontrivial non-Archimedean valuation.  

\begin{defn}
  \begin{enum}
  \item An (affine) \emph{half-space} in $N_\R$ is a subset of the form 
    \[ H = \{ v\in N_\R~:~ \angles{u,v} \leq a \}\sptxt{for some}
    u\in M_\R\setminus\{0\},~ a\in\R. \]
    The half-space $H$ is called \emph{integral} if we can take $u\in M$,
    and is \emph{integral $\Gamma$-affine} if we can take $u\in M$ and $a\in\Gamma$.
    The half-space $H$ is called \emph{linear} if we can take $a = 0$.

  \item With $H\subset N_\R$ as above, the \emph{complementary half-space}
    of $H$ is  
    \[ H^- = \{ v\in N_\R~:~\angles{u,v} \geq a \} \]
    and the \emph{boundary of $H$} is its topological boundary
    \[ \del H = \{ v\in N_\R~:~ \angles{u,v} = a \} = H\cap H^-. \]

  \item An \emph{affine space} in $N_\R$ is a translate of a linear
    subspace of $N_\R$.  Any affine space is of the form 
    $\bigcap_{i=1}^r \del H_i$, where the $H_i$ are half-spaces.

  \item A \emph{polyhedron} in $N_\R$ is a nonempty intersection 
    $P = \bigcap_{i=1}^r H_i$ of finitely many half-spaces 
    $H_i\subset N_\R$.  We say 
    that $P$ is \emph{integral} (resp. \emph{integral $\Gamma$-affine}) if we can
    take the $H_i$ to be integral (resp. integral $\Gamma$-affine).  If
    $P$ is integral $\Gamma$-affine we set $P_\Gamma = P\cap N_\Gamma$.

  \item An integral, resp.\  integral $\Gamma$-affine, \emph{polytope} is a bounded
    integral, resp.\  bounded integral $\Gamma$-affine, polyhedron.

  \item Let $P\subset N_\R$ be a polyhedron.  The \emph{affine span} of
    $P$, denoted $\spn(P)$, is
    the smallest affine subspace of $N_\R$ containing $P$.  The
    \emph{dimension} $\dim(P)$ of 
    $P$ is the dimension of $\spn(P)$.  The \emph{relative interior} of $P$,
    denoted $\relint(P)$, is the interior of $P$ as a subspace of $\spn(P)$.

  \item Let $S\subset N_\R$ be a subset.  The \emph{convex hull} of $S$ is the
  intersection $\conv(S)$ of all half-spaces in $N_\R$ containing $S$.  It
  is the smallest convex subset of $N_\R$ containing $S$.
  \end{enum}
\end{defn}

See~\parref{eg:polytopes} for examples.

\begin{defn}
  Let $P\subset N_\R$ be a polyhedron. For  $u\in M_\R$ we define
  \[ \face_u(P) = \maxset(u, P) = 
  \{ v\in P~:~ \angles{u,v}\geq\angles{u,v'}\text{ for all } v'\in P \}. \]
  A \emph{face} of $P$ is a nonempty subset of the form $F = \face_u(P)$ for some
  $u\in M_\R$.  We write $F\prec P$ to signify that $F$ is a face of $P$.
  A \emph{vertex} of $P$ is a face consisting of a single point; we let
  $\vertices(P)$ denote the set of vertices of $P$.

\end{defn}

In other words, a face of $P$ is a subset on which a linear form attains
its maximum.  Note that using these conventions we have $P\prec P$ but
$\emptyset\not\prec P$.  

\begin{eg} Let \label{eg:polytopes}
  $N = M = \Z^2\subset\R^2$, and let $\angles{\cdot,\cdot}$ be the dot
  product.
  The unit square $S = [0,1]^2$ is an integral $\Z$-affine polytope in
  $\R^2$, and the first 
  quadrant $Q =\R_{\geq0}^2$ is an integral $\Z$-affine polyhedron.  The
  four edges and four vertices of $S$ are faces; if 
  $u_1 = (-1,0), u_2 = (-1,-1),$ and $u_3 = (0,-1)$ then
  the left edge is $\face_{u_1}(S)$, the bottom edge is
  $\face_{u_2}(S)$, and $\{(0,0)\} = \face_{u_3}(S)$.  
  The polyhedron $Q$ has four faces: $Q$ itself, two edges $\face_{u_1}(Q)$ and
  $\face_{u_2}(Q)$, and $\{(0,0)\} = \face_{u_3}(Q)$.  Note that all faces
  are again integral $\Z$-affine, and that 
  $S$ is the convex hull of its vertices.  See
  Figure~\ref{fig:polytopes}. 

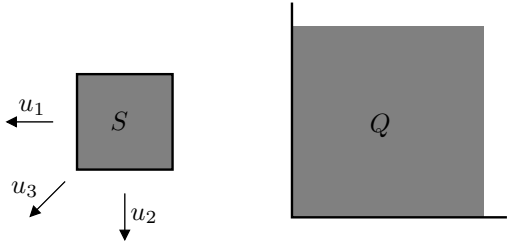
  \genericfig[ht]{polytopes}{The unit square is a polytope and the first
    quadrant is a polyhedron in $\R^2$.}

\end{eg}

Many statements about polyhedra can be deduced from the analogous results for
cones~\parref{prop:cone.properties} by considering the cone over a
polyhedron: see~\cite[p.24]{fulton:toric_varieties}.
Here we collect some of the basic properties of polyhedra:

\begin{prop} 
  Let $P\subset N_\R$ be a polyhedron.
  \begin{enum}
  \item A face $F\prec P$ is a polyhedron in $N_\R$.  If $F\neq P$ then
    $\dim(F)<\dim(P)$.
  \item If $F\prec P$ and $P$ is integral (resp. integral $\Gamma$-affine) then
    $F$ has the same property.
  \item If $F,F'\prec P$ and $F\subset F'$ then $F\prec F'$.  More
    generally, if $F\cap F'\neq\emptyset$ then $F\cap F'$ is a face of
    $F$ and of $F'$ (and of $P$).
  \item $P$ has finitely many faces.
  \item If $P$ is a polytope then $P = \conv(\vertices(P))$, and the
    convex hull of a finite set of points is a
    polytope~\cite[Corollary~2.4.3]{barvinok:convexity}. 

  \end{enum}
\end{prop}

\smallskip
Collections of polyhedra will also be of interest:

\begin{defn} \label{defn:polyhedral.complexes}
  \begin{enum}
  \item A \emph{polyhedral complex} is a finite collection $\Pi$ of polyhedra
    in $N_\R$, called the \emph{cells} or \emph{faces} of $\Pi$, satisfying
    \begin{itemize}
    \item[(PC1)] if $P,P'\in\Pi$ and $P\cap P'\neq\emptyset$ then $P\cap P'$ is a face of $P$ and of $P'$, and
    \item[(PC2)] if $P\in\Pi$ and $F\prec P$ then $F\in\Pi$.
    \end{itemize}
    The \emph{support} of $\Pi$ is the set $|\Pi| = \bigcup_{P\in\Pi} P$.
    The \emph{dimension} of $\Pi$ is the dimension of its
    highest-dimensional cell; $\Pi$ has \emph{pure dimension $d$} if
    every maximal cell has dimension $d$.
    We say that $\Pi$ is \emph{integral} (resp. \emph{integral $\Gamma$-affine})
    if all of its cells are integral (resp. integral $\Gamma$-affine).
    
  \item A \emph{polytopal complex} is a polyhedral complex whose cells are
    polytopes.   

  \item A \emph{refinement} of a polyhedral complex $\Pi$ is a polyhedral
    complex $\Pi'$ with the same support, and such that each cell of $\Pi$
    is a union of cells of $\Pi'$.

  \item Let $\Pi,\Pi'$ be polyhedral complexes.  We define 
    \[ \Pi\cap\Pi' = \{ \text{all faces of } P\cap P'~:~
    P\in\Pi,\,P'\in\Pi',\text{ and }P\cap P'\neq\emptyset \}. \]
    It is easy to show
    that $\Pi\cap\Pi'$ is a polyhedral complex, and that
    $|\Pi\cap\Pi'| = |\Pi|\cap|\Pi'|$.  In particular, if $|\Pi|=|\Pi'|$
    then $\Pi\cap\Pi'$ is a common refinement of $\Pi$ and $\Pi'$.

  \end{enum}
\end{defn}

\begin{eg} We \label{eg:complexes}
  let $N = M = \Z^2$ as in~\parref{eg:polytopes}.  Let
  \[\begin{split} P_1 &= \R_{\geq 0}(0,1) \quad
  P_2 = \R_{\geq 0}(1,0) \quad
  P_3 = \conv\{(-1,-1),(0,0)\}\\
  P_4 &= (-1,-1) + \R_{\geq 0}(-1,0) \quad
  P_5 = (-1,-1) + \R_{\geq 0}(0,-1) \end{split}\]
  as in Figure~\ref{fig:complexes}, and let
  \[ \Pi_1 = \{P_1,P_2,P_3,P_4,P_5,\{(0,0)\},\{(-1,-1)\}\}. \]
  Then $\Pi_1$ is an integral $\Z$-affine polyhedral complex of pure
  dimension $1$ in $\R^2$.  

  Let $Q_1$ denote the triangle $\conv\{(0,0),(0,1),(1,0)\}$ and let
  $Q_2 = \conv\{(1,1),(0,1),(1,0)\}$, as in Figure~\ref{fig:complexes}.
  These are integral $\Z$-affine polytopes with three vertices and three
  edges.  Let $\Pi_2$ be the collection of all faces of $Q_1$ and $Q_2$.
  This is an integral $\Z$-affine polytopal complex of pure dimension $2$
  in $\R^2$.  It contains four vertices, five edges, and two faces, and
  its support $\Pi_2$ is the unit square.
  
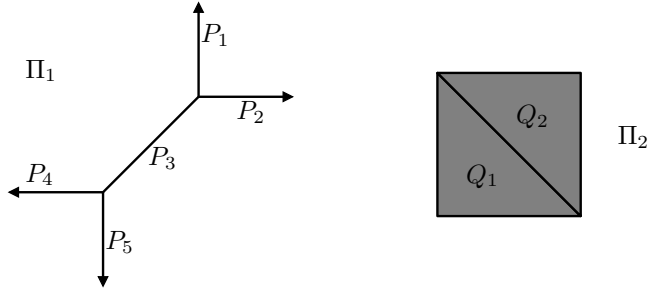
  \genericfig[ht]{complexes}{An integral $\Z$-affine polyhedral complex
    $\Pi_1$ of pure dimension $1$ and an integral $\Z$-affine polytopal
    complex $\Pi_2$ of pure dimension $2$ in $\R^2$.  Only the maximal
    cells are labeled.} 

\end{eg}

Intersections of \emph{linear} half-spaces are called cones:

\begin{genthm}{Definition/Proposition}{true}{}{paragraph}\rmfamily\upshape
  \label{prop:cone.properties} 
  \begin{enum}
  \item A (convex polyhedral) \emph{cone} (resp. \emph{integral cone}) in $N_\R$
    is an intersection 
    $\sigma$ of finitely many linear (resp integral linear)
    half-spaces in $N_\R$.  Any face of a cone is a cone.  We say that
    $\sigma$ is \emph{pointed} if $0$ 
    is a vertex of $\sigma$, or equivalently if $\sigma$ does not contain
    any nonzero linear space.  

  \item Let $v_1,\ldots,v_r\in N_\R$.  The subset 
    $\sigma = \sum_{i=1}^r \R_{\geq 0} v_i$
    is a cone in $N_\R$, and any cone can be written in this
    form~\cite[p.12]{fulton:toric_varieties}. 
    The cone $\sigma$ is integral if $v_i\in N$ for all $i$, and any
    integral cone can be written 
    $\sigma = \sum_{i=1}^r \R_{\geq 0} v_i$
    for $v_1,\ldots,v_r\in N$.

  \item Let $\sigma = \sum_{i=1}^r \R_{\geq 0} v_i \subset N_\R$ be a
    cone.  The (polar) \emph{dual cone} to $\sigma$ is the cone
    \[ \sigma^\vee = \{ u\in M_\R~:~\angles{u,v}\leq 0
    \text{ for all }v\in \sigma\}
    = \bigcap_{i=1}^r \{u\in M_\R~:~\angles{u,v_i}\leq 0\}. \]
    We have $\sigma = \sigma^{\vee\vee}$~\cite[(1.2.1)]{fulton:toric_varieties},
    and $\sigma$ is integral if and only if $\sigma^\vee$ is integral.

  \item The \emph{annihilator} of a cone $\sigma\subset N_\R$ is the
    annihilator of the vector space $\spn(\sigma)$: 
    \[ \sigma^\perp = \{ u\in M_\R~:~\angles{u,v} = 0
    \text{ for all }v\in \sigma\}. \]
    It is a linear space in $M_\R$.

  \item Let $\sigma\subset N_\R$ be a cone.  The map 
    $\tau\mapsto\tau^\perp\cap\sigma^\vee$ is an inclusion-reversing
    bijection between the 
    faces of $\sigma$ and the faces of $\sigma^\vee$, with inverse
    $\tau'\mapsto(\tau')^\perp\cap\sigma$
    \cite[(1.2.10)]{fulton:toric_varieties}.  We have
    $\dim(\tau) + \dim(\tau^\perp\cap\sigma^\vee) = n$.

  \item A \emph{fan} $\Delta$ in $N_\R$ is a polyhedral complex whose
    cells are cones (called the \emph{cones} of $\Delta$).  The fan
    $\Delta$ is \emph{complete} if  
    $|\Delta| = N_\R$.  The fan $\Delta$ is \emph{pointed} if $\{0\}\in\Delta$,
    or equivalently if all cells of $\Delta$ are pointed cones.

  \item Let $P\subset N_\R$ be a polyhedron.  The \emph{normal fan} to $P$
    is the fan $\cN(P)$ in $M_\R$ whose cells are the cones
    \[ \cN(P,F) \coloneq \{ u\in M~:~ F\subset\face_u(P) \}
    \sptxt{where} F\prec P. \]
    This fan is integral if $P$ is integral.
    See~\cite[p.26]{fulton:toric_varieties}. 

  \item Let $P$ be a polyhedron.  Its normal fan $\cN(P)$ is complete if
    and only if $P$ is a polytope, and $\cN(P)$ is pointed if and only if
    $\dim(P) = \dim_\R(N_\R)$.

  \end{enum}
\end{genthm}

\begin{eg} Let \label{eg:cones}
  $N = M = \Z^2$ as in~\parref{eg:polytopes}.
  Let $\sigma = \R_{\geq 0}(0,1) + \R_{\geq 0} (1,1)$.  This is an
  integral pointed cone in $\R^2$.  It has four faces: $\sigma$ itself,
  $\tau_1 = \R_{\geq 0}(0,1)$, $\tau_2 = \R_{\geq 0}(1,1)$, and
  $\{(0,0)\}$.  The dual cone is 
  $\sigma^\vee = \R_{\geq 0}(-1,0) + \R_{\geq 0}(1,-1)$; its faces are
  \[\begin{split} 
  \tau_1' &= \R_{\geq 0}(-1,0) = \tau_1^\perp\cap\sigma^\vee\quad
  \tau_2' = \R_{\geq 0}(1,-1) = \tau_2^\perp\cap\sigma^\vee\\
  \sigma^\vee &= \{(0,0)\}^\perp\cap\sigma^\vee\quad
  \{(0,0)\} = \sigma^\perp\cap\sigma^\vee. \end{split}\] 
  See Figure~\ref{fig:cones}.
  
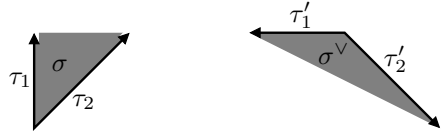
  \genericfig[ht]{cones}{An integral pointed cone $\sigma$ in $\R^2$ and
    its polar dual $\sigma^\vee$.  For $i=1,2$ we have
    $\tau_i' = \tau_i^\perp\cap\sigma^\vee$.}
\end{eg}

\begin{eg} Let \label{eg:fans}
  $N = M = \Z^2$ and let $S = [0,1]^2$, as in~\parref{eg:polytopes}.  
  Label the vertices $v_1,v_2,v_3,v_4$ and edges $F_1,F_2,F_3,F_4$ of $S$
  as in Figure~\ref{fig:normalfan}.  The normal fan to the polytope $S$ is
  drawn in Figure~\ref{fig:normalfan}; it is a complete integral pointed
  fan.  For $i=1,2,3,4$ the set of $u\in\R^2$ such that $v_i\in\face_u(S)$
  is the $i$th quadrant $\sigma_i$.  To say that $F_i\subset\face_u(S)$ is
  to say that both $v_i$ and $v_{i+1}$ are in $\face_u(S)$, so 
  $\cN(S,F_i) = \sigma_i\cap\sigma_{i+1} = \tau_i$, where the subscripts
  are taken modulo $4$.

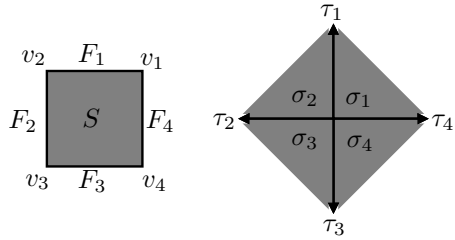
  \genericfig[ht]{normalfan}{The unit square and its normal fan.  For
    $i=1,2,3,4$ we have $\cN(S,v_i) = \sigma_i$ and
    $\cN(S,F_i) = \tau_i$}
  
\end{eg}

\section{Compactification procedures}
\label{sec:compactification}

\paragraph 
In this section we describe a procedure for constructing a partial
compactification $N_\R(\Delta)$ of $N_\R$ associated to a fan
$\Delta$.  This procedure is analogous to the construction of 
the toric variety $X(\Delta)$ associated to $\Delta$
(see~\S\ref{sec:extended.trop}); the space 
$N_\R(\Delta)$ will serve as the (extended) tropicalization of
$X(\Delta)$.  We then describe the closure $\bar P$ of a polyhedron $P$ in a
suitable partial compactification $N_\R(\Delta)$. 
The compactification $\bar P$ will 
correspond to a ``polyhedral subdomain'' of $X(\Delta)$; this 
generalizes~\cite[\S4]{gubler:tropical} and~\cite[\S3]{ekl:non_arch_amoebas}.

The construction of $N_\R(\Delta)$ is originally due to
Kajiwara~\cite{kajiwara:tropical_toric_geometry},  and was
later described by Payne~\cite[\S3]{payne:analytification}.  We follow Payne's
treatment.  

\begin{notn}
  We let $\bar\R$ be the additive monoid $\R\cup\{-\infty\}$, endowed with
  the topology which restricts to the standard topology on $\R$ and for
  which the sets of the form $[-\infty,a)$ for $a\in\R$ constitute
  a neighborhood basis of $-\infty$.
\end{notn}

\begin{defn} Let \label{defn:partial.compactification}
  $\sigma\subset N_\R$ be a cone.  The 
  \emph{partial compactification of $N_\R$ with respect to $\sigma$} is
  the space $N_\R(\sigma) = \Hom_{\R_{\geq 0}}(\sigma^\vee,\bar\R)$ of
  monoid homomorphisms respecting multiplication by $\R_{\geq 0}$,
  equipped with the topology of pointwise convergence.  We
  use $\angles{\cdot,\cdot}_\sigma$ to denote the the pairing
  $\sigma^\vee\times N_\R(\sigma)\to\bar\R$.
\end{defn}

See~\parref{eg:proj.plane} for an example.
Roughly, $N_\R(\sigma)$ is a space that compactifies $N_\R$ in the
directions of the faces of $\sigma$; this statement is made precise in the
following proposition.  By \emph{topological embedding} we mean an
injection of topological spaces that is a homeomorphism onto its image. 

\begin{prop}[{\cite[\S3]{payne:analytification}}] Let \label{prop:tropsigma}
  $\sigma\subset N_\R$ be a cone.
  \begin{enum}
  \item
    Let $\tau\prec\sigma$, let $v\in N_\R/\spn(\tau)$, and define
    $\iota(v)\in N_\R(\sigma)$ by
    \[ \angles{u,\iota(v)}_\sigma =
    \begin{cases}
      \angles{u,v} & \quad\text{ if } u\in\tau^\perp\cap\sigma^\vee \\
      -\infty      & \quad\text{ otherwise}
    \end{cases}
    \]
    for $u\in\sigma^\vee$.  Then $\iota(v)$ is a well-defined element of
    $N_\R(\sigma)$, and  
    \[ \iota~:~ \Djunion_{\tau\prec\sigma} N_\R/\spn(\tau) \isom N_\R(\sigma) \]
    is a bijection.  Furthermore, for each $\tau\prec\sigma$ the restriction
    $\iota|_{N_\R/\spn(\tau)}\inject N_\R(\sigma)$ is a
    topological embedding.

  \item 
    If $\sigma^\vee = \sum_{i=1}^r \R_{\geq 0} u_i$ then the map
    \[ v\mapsto(\angles{u_1,v}_\sigma,\ldots,\angles{u_r,v}_\sigma)~:~
    N_\R(\sigma)\inject\bar\R^r \]
    is a topological embedding with closed image. 

  \item 
    For $\tau\prec\sigma$ the inclusion $\sigma^\vee\subset\tau^\vee$
    induces a topological embedding $N_\R(\tau)\inject N_\R(\sigma)$ with
    open image.

  \end{enum}
\end{prop}

\pf
We will only prove~(i).
Since $\tau^\perp\cap\sigma^\vee$ is a face of $\sigma^\vee$
\parref{prop:cone.properties}, we have 
$u_1+u_2\in\tau^\perp\cap\sigma^\vee\iff
u_1,u_2\in\tau^\perp\cap\sigma^\vee$, which shows that $\iota(v)$ 
is a well-defined element of $N_\R(\sigma)$.  
We claim that $\iota$ is injective.
For $v\in N_\R/\spn(\tau)$ it is clear from the definition that we can
recover $\tau^\perp\cap\sigma^\vee$, and hence that we can recover
$\tau$, from the element $\iota(v)$, so it suffices to show that
$\iota|_{N_\R/\spn(\tau)}$ is injective for all $\tau\prec\sigma$.  This
follows from the fact that $\tau^\perp\cap\sigma^\vee$
spans $\tau^\perp$~(\ref{prop:cone.properties},v).  As for surjectivity: given any
$v_0\in N_\R(\sigma)$, the set $v_0\inv(\R)$ is a face of
$\sigma^\vee$, and is hence of the form $\tau^\perp\cap\sigma^\vee$;
by linear algebra, we conclude that $v_0 = \iota(v)$ for suitable
$v\in N_\R/\spn(\tau)$.

The topology on $N_\R/\spn(\tau)$ coincides with the the topology of
pointwise convergence, thinking of $N_\R/\spn(\tau)$ as the space of
linear functions on $\tau^\perp$.  It follows that $\iota|_{N_\R/\spn(\tau)}$
is a topological embedding. \qed

From this point on we will identify 
$\Djunion_{\tau\prec\sigma} N_\R/\spn(\tau)$ with $N_\R(\sigma)$ without
mentioning the map $\iota$.  See~\parref{eg:proj.plane} for an example.

\begin{rem} \label{rem:after.tropsigma}
  \begin{enum}
  \item    
    Later in this section we will give a ``local''
    description of the topology on $N_\R(\sigma)$; cf.
    \parref{rem:local.description}.

  \item
    If $\sigma$ is a pointed cone then $N_\R = N_\R(\{0\})$ naturally sits
    inside of $N_\R(\sigma)$ by~(\ref{prop:tropsigma}\,i).  
    In this case $N_\R$ is dense in $N_\R(\sigma)$.

  \item
    Let $\sigma$ be a cone and let $\tau\prec\sigma$.  Then for
    $v\in N_\R(\tau)$ and $u\in\sigma^\vee\subset\tau^\vee$, by definition we
    have $\angles{u,v}_\sigma = \angles{u,v}_\tau$ under the natural
    inclusion $N_\R(\tau)\inject N_\R(\sigma)$. 

  \item 
    We mentioned in the course of the proof of~\parref{prop:tropsigma}
    that $v\in N_\R/\spn(\tau)$ if and only if 
    $v\inv(\R) = \tau^\perp\cap\sigma^\vee$.
  \end{enum}
\end{rem}

\begin{defn} Let \label{defn:N.R.Delta}
  $\Delta$ be a pointed fan in $N_\R$.  The 
  \emph{partial compactification of $N_\R$ with respect to $\Delta$} is
  the space $N_\R(\Delta)$ obtained by gluing the spaces
  $N_\R(\sigma)$ for $\sigma\in\Delta$ using the open immersions
  $N_\R(\tau)\inject N_\R(\sigma)$ for $\tau\prec\sigma$.
\end{defn}

It follows from~(\ref{prop:tropsigma}\,i) that there is a canonical
bijection
\[ \Djunion_{\sigma\in\Delta} N_\R/\spn(\sigma) \isom N_\R(\Delta). \]
Moreover if $\Delta$ is the fan whose cones are the faces of a single
cone $\sigma$, then $N_\R(\sigma)$ is canonically identified with
$N_\R(\Delta)$.  See~\parref{eg:proj.plane}

\begin{eg}[The affine and projective planes] Let \label{eg:proj.plane}
  $\sigma_1$ be the first quadrant in $N_\R = \R^2$ (the toric variety
  associated to $\sigma_1$ is isomorphic to the affine plane).  The faces of
  $\sigma_1$ are $\sigma_1$ itself, $\tau_1 = \R_{\geq0}(1,0)$, 
  $\tau_2=\R_{\geq0}(0,1)$, and $\{0\}$.  Therefore
  \begin{equation} \label{eq:decomp.tropsigma} \begin{split}
    N_\R(\sigma_1) &= N_\R\djunion(N_\R/\spn(\tau_1))\djunion(N_\R/\spn(\tau_2))
    \djunion(N_\R/\spn(\sigma_1)) \\
    &=\R^2\djunion(\{+\infty\}\times\R)\djunion(\R\times\{+\infty\})\djunion\{(+\infty,+\infty)\}
    = (\R\djunion\{\infty\})^2. 
  \end{split}\end{equation}

  Let $\sigma_2 = \R_{\geq 0}(1,0) + \R_{\geq 0}(-1,-1)$ and
  $\sigma_3 = \R_{\geq 0}(0,1) + \R_{\geq 0}(-1,-1)$, let
  $\tau_3 = \R_{\geq 0}(-1,-1)$, and let
  $\Delta = \{ \sigma_1,\sigma_2,\sigma_3,\tau_1,\tau_2,\tau_3,\{0\}\}$.
  Then $\Delta$ is a complete integral pointed fan (its associated toric
  variety is the projective plane; cf.~\parref{eg:trop.projspace}).  By
  definition~\parref{defn:N.R.Delta}   we have
  $N_\R(\Delta) = N_\R(\sigma_1)\cup N_\R(\sigma_2)\cup N_\R(\sigma_3)$.
  See Figure~\ref{fig:tropdelta} for a picture.

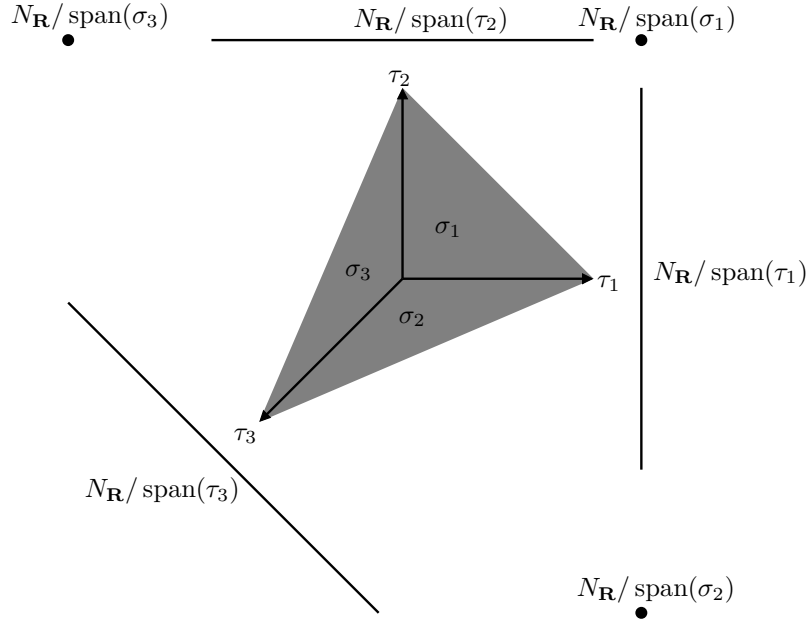
  \genericfig[ht]{tropdelta}{A picture of $N_\R(\Delta)$, where $\Delta$
    is the fan of~\parref{eg:proj.plane}.}   

\end{eg}

\paragraph The \label{par:functoriality.Delta}
construction of $N_\R(\Delta)$ is functorial in $\Delta$,
in the following sense.  Let $N_\R,N_\R'$ be finite-dimensional real
vector spaces and let $\sigma$ (resp. $\sigma'$) be a cone in $N_\R$
(resp. $N_\R'$).  Let $\phi: N'_\R\to N_\R$ be a linear map with
$\phi(\sigma')\subset\sigma$.  The dual map
$\phi^*:M_\R\to M'_\R = (N'_\R)^*$ induces a monoid homomorphism 
$\sigma^\vee\to(\sigma')^\vee$, and hence a continuous map
$\phi: N_\R'(\sigma')\to N_\R(\sigma)$ extending $\phi$.

Now let $\Delta$ (resp. $\Delta'$) be a pointed fan in $N_\R$ (resp. $N'_\R$), and let
$\phi:N'_\R\to N_\R$ be a linear map respecting the fans $\Delta',\Delta$, i.e.,
such that for every $\sigma'\in\Delta'$ there exists
$\sigma\in\Delta$ such that $\phi(\sigma')\subset\sigma$.  Then we can
glue the maps $N_\R'(\sigma')\to N_\R(\sigma)$ to give a 
continuous map $\phi: N_\R(\Delta')\to N_\R(\Delta)$ extending $\phi$.

\subparagraph
It is clear from the construction and
(\ref{prop:tropsigma},iii) that if $\Delta'$ is a subfan of
$\Delta$ (i.e.\ if every cone in $\Delta'$ is a cone in $\Delta$) then
$N_\R(\Delta')\to N_\R(\Delta)$ is an open immersion.   

\begin{rem}
  One can show that $N_\R(\Delta)$ is compact if and only if $\Delta$ is a
  complete fan.  More generally, if $\phi:N'\to N$ is a linear map
  respecting fans $\Delta'$ and $\Delta$ as above, then 
  the extending map $\phi: N_\R(\Delta')\to N_\R(\Delta)$ is proper if and
  only if $\phi\inv(|\Delta|) = |\Delta'|$.  This mirrors the situation
  for toric varieties.  See~\cite[\S3]{payne:analytification}.
\end{rem}

\begin{defn} In \label{defn:S.sigma}
  the case when $\sigma\subset N$ is an
  \emph{integral} pointed cone, we define
  \[ S_\sigma \coloneq \sigma^\vee\cap M \]
  and
  \[ N_\Gamma(\sigma) \coloneq
  \Hom(S_\sigma,\Gamma\cup\{-\infty\}) =
  \{
  v\in N_\R(\sigma)~:~v(S_\sigma)\subset\Gamma\cup\{-\infty\}
  \} \subset N_\R(\sigma) \] 
  with $\Gamma$ as in~\parref{notn:convex}.  If $\Delta$ is an integral
  pointed fan, we define 
  \[ N_\Gamma(\Delta) = \bigcup_{\sigma\in\Delta} N_\Gamma(\sigma)
  \subset N_\R(\Delta). \]
\end{defn}

As above we have
\[ \Djunion_{\tau\prec\sigma} N_\Gamma/(\spn(\tau)\cap N_\Gamma)
\isom N_\Gamma(\sigma) \sptxt{and}
\Djunion_{\sigma\in\Delta} N_\Gamma/(\spn(\sigma)\cap N_\Gamma)
\isom N_\Gamma(\Delta). \]
The constructions of $N_\Gamma(\sigma)$ and $N_\Gamma(\Delta)$ are
functorial with respect to linear maps $\phi:N'_\R\to N_\R$  as
in~\parref{par:functoriality.Delta} such that $\phi(N')\subset N$.

\paragraph We proceed with the compactification of a polyhedron in $N_\R$.  
More specifically, we will take the closure of a polyhedron $P$ inside of
a partial compactification $N_\R(\sigma)$, but this will only make sense when
$\sigma$ partially compactifies $N_\R$ in the directions in which $P$ is
infinite.  

\begin{defn}
  Let $P\subset N_\R$ be a polyhedron.  The 
  \emph{cone of unbounded directions} or \emph{recession cone} of $P$ is
  the cone $\cU(P)$ which is polar dual to the cone
  \[ \cU(P)^\vee \coloneq \{ u\in M_\R~:~\face_u(P)\neq\emptyset\} = |\cN(P)|. \] 
  We say that $P$ is \emph{pointed} if $\cU(P)$ is pointed, or
  equivalently if $P$ does not contain an affine space.
\end{defn}

If $P = \bigcap_{i=1}^r \{ v\in N_\R:\angles{u_i,v}\leq a_i \}$ 
then $\cU(P)^\vee = \sum_{i=1}^r \R_{\geq0}u_i$ and 
$\cU(P) = \bigcap_{i=1}^r \{v\in N_\R:\angles{u_i,v}\leq 0\}$.
It follows that $\cU(P)^\vee = M_\R$ if and only if $P$ is bounded, and
that $\cU(P)$ is integral when $P$ is integral.
See~\cite[\S2.16]{barvinok:convexity}. 

\begin{eg} Let \label{eg:unbounded}
  $N_\R = \R^2$ and let $P\subset\R^2$ be the polyhedron
  \[ P = \{ (x,y)\in\R^2~:~x\geq 1,\, y\geq 1,\, x+y\geq 3 \}. \]
  We have
  \[ \cU(P)^\vee = \R_{\geq 0}(-1,0)+\R_{\geq0}(0,-1)+\R_{\geq0}(-1,-1) 
  = \R_{\geq 0} (-1,0) + \R_{\geq 0} (0,-1) \]
  and therefore its cone of unbounded directions is the first quadrant.  See
  Figure~\ref{fig:unbounded}.  Note that $P = \conv\{(1,2),(2,1)\}+\cU(P)$.
  
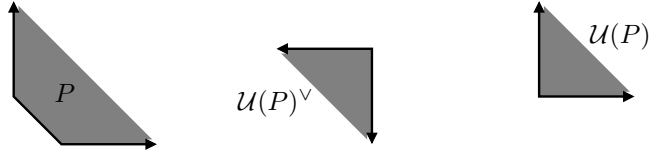
  \genericfig[ht]{unbounded}{A polyhedron $P$, the cone
    $\cU(P)^\vee$, and its cone of unbounded directions $\cU(P)$.}

\end{eg}

The following lemma is standard:

\begin{lem} Let \label{lem:bounded.faces}
  $P\subset N_\R$ be a pointed polyhedron and let $\sigma=\cU(P)$.  
  \begin{enum}
  \item $\displaystyle\cU(P) = \{ v'\in N_\R ~:~ v + \R_{\geq0} v'\subset P
    \text{ for all (resp. any) } v\in P \}$.
  \item If $F_b$ denotes the union of the bounded faces of $P$ then
    \[ P = F_b + \sigma = \{ u_1 + u_2 ~:~ u_1\in F_b,~ u_2\in\sigma \}. \]
  \end{enum}
\end{lem}

\pf We will only prove~(ii).
Write $P = \bigcap_{i=1}^r \{ v\in N_\R:\angles{u_i,v}\leq a_i \}$ for
some $u_i\in\sigma^\vee$ and $a_i\in\R$.
Let $v_1\in P$ and $v_2\in\sigma$.  For each $i$ we have
\[ \angles{u_i,v_1+v_2} = \angles{u_i,v_1} + \angles{u_i,v_2}
\leq a_i \]
since $\angles{u_i,v_2}\leq 0$ by definition.  This shows that
$P + \sigma\subset P$, so $F_b+\sigma\subset P$.

For the other inclusion, let $v\in P$ be arbitrary, and let $F$ be the
unique face of $P$ such that $v$ is contained in the relative interior
of $F$.  We will prove by induction on $\dim(F)$ that 
$v\in F_b+\sigma$.  If $\dim(F) = 0$ then we are done because $F$ is
bounded, so suppose that $\dim(F) > 0$.  Let $u\in\sigma^\vee$ be such that
$F = \face_u(P)$.  If $u$ is in the interior of $\sigma^\vee$ then $F$ is
bounded, and otherwise there exists nonzero $v_0\in\sigma$ such that
$\angles{v_0,u} = 0$.  Let $a_0 = \max\{a\in\R:v - av_0\in P\}$ ---
this is finite because $\angles{u_i,v_0}\neq 0$ for some $i$ ---
and let $v_1 = v - a_0v_0$.  By construction, $v_1$ is in the boundary
of $F$, and hence is contained in the relative interior of a face of
strictly smaller dimension.  This proves that $P = F_b + \sigma$.\qed

We omit the proofs of the following two lemmas, which follow more or less
immediately from \parref{lem:bounded.faces}.

\begin{lem} Let \label{lem:face.unbounded}
  $F$ be a face of a pointed polyhedron $P\subset N_\R$.  Then
  $\cU(F)\prec\cU(P)$.
\end{lem}

\begin{lem} Let \label{lem:unbounded.intersection}
  $P,P'\subset N_\R$ be pointed polyhedra such that
  $P\cap P'\neq\emptyset$.  Then $\cU(P\cap P')=\cU(P)\cap\cU(P')$.
\end{lem}

\paragraph Let
$P = \bigcap_{i=1}^r \{ v\in N_\R:\angles{u_i,v}\leq a_i \}$ be a
pointed polyhedron and let $\sigma = \cU(P)$.   Then we have
\begin{equation} \label{eq:P.inequality} 
P = \{ v\in N_\R~:~\angles{u,v}\leq\max_{v'\in P}\angles{u,v'}
\text{ for all } u\in\sigma^\vee \} \end{equation}
because $u_1,\ldots,u_r\in\sigma^\vee$.  More generally, let
$\tau\prec\sigma$ and let $\pi_\tau: N_\R\to N_\R/\spn(\tau)$ denote the
projection.  Then
\begin{equation} \label{eq:pi.tau.P} \pi_\tau(P)
= \{ v\in N_\R/\spn(\tau)~:~\angles{u,v}_\sigma\leq
\max_{v'\in P}\angles{u,v'}\text{ for all } u\in\sigma^\vee \}. \end{equation}
This can be seen as follows: one inclusion is clear, so suppose that 
$v\in N_\R/\spn(\tau)$ satisfies 
$\angles{u,v}\leq\max_{v'\in P}\angles{u,v'}$ for all 
$u\in\sigma^\vee\cap\tau^\perp$.  If $v_1\in N_\R$ lifts $v$ and satisfies
$\angles{u_i,v}\leq a_i$ for all $u_i\notin\tau^\perp$ then $v_1\in P$ by
the above, so $v\in\pi_\tau(P)$.

\begin{defn}
  Let $P\subset N_\R$ be a pointed polyhedron and let $\sigma = \cU(P)$.
  The \emph{compactification} $\bar P$ of $P$ is the closure of $P$ in 
  $N_\R(\sigma)$.  If $P$ is integral $\Gamma$-affine then we set 
  $\bar P_\Gamma = \bar P\cap N_\Gamma(\sigma)$.
\end{defn}

See~\parref{eg:P.compactification} for an example.

\subparagraph
Let $P = \bigcap_{i=1}^r \{ v\in N_\R:\angles{u_i,v}\leq a_i \}$ be a
pointed polyhedron, and define $f: N_\R(\sigma)\inject\bar\R^r$ by 
$f(v) = (\angles{u_1,v}_\sigma,\ldots,\angles{u_r,v}_\sigma)$ as
in~\parref{prop:tropsigma}. 
Then $f(\bar P)$ is a closed subset of the compact space
$\prod_{i=1}^r \bar\R_{\leq a_i}$, so $\bar P$ is compact.

\begin{prop} Let \label{prop:compactified.polyhedron}
  $P = \bigcap_{i=1}^r \{ v\in N_\R:\angles{u_i,v}\leq a_i \}$ be a
  pointed polyhedron with cone of unbounded directions $\sigma$.  Then
  \[ \begin{split}
    \bar P &= \Djunion_{\tau\prec\sigma} \pi_\tau(P) 
    = \Djunion_{\tau\prec\sigma} \{ v\in N_\R/\spn(\tau)~:~
    \angles{u,v}_\sigma\leq\max_{v'\in P}\angles{u,v'}
    \text{ for all } u\in\sigma^\vee \} \\
    &= \Djunion_{\tau\prec\sigma} \{ v\in N_\R/\spn(\tau)~:~
    \angles{u_i,v}\leq a_i\text{ for all } u_i\in\tau^\perp \} \\
    &= \{u:\sigma^\vee\to\bar\R~:~\angles{u,v}_\sigma\leq
    \max_{v'\in P}\angles{u,v'}\text{ for all } u\in\sigma^\vee \}.
  \end{split}\]
\end{prop}

\pf The second equality was proved in~\eqref{eq:P.inequality}, the
third equality follows from the fact that $\tau^\perp\cap\sigma^\vee$ is
spanned by the $u_i$ contained in $\tau^\perp$, and the last equality is
obvious.  The set
\[ \{u:\sigma^\vee\to\bar\R~:~\angles{u,v}_\sigma\leq
\max_{v'\in P}\angles{u,v'}\text{ for all } u\in\sigma^\vee \} \]
is closed since $\angles{u,v}_\sigma\leq\max_{v'\in P}\angles{u,v'}$ is a
closed condition for fixed $v$, so $\bar P$ is contained in the right-hand
side by~\eqref{eq:pi.tau.P}.  Conversely, let $\tau\prec\sigma$ be a face
of positive dimension, let $v\in N_\R/\spn(\tau)$ be
such that $\angles{u,v}\leq\max_{v'\in P}\angles{u,v'}$ for all
$u\in\sigma^\vee\cap\tau^\perp$, and let $v_1\in P$ be a lift of $v$.  Let
$v_2$ be in the relative interior of $\tau$.  Then 
$\pi_\tau(v_1 + av_2) = v$ for all $a\in\R$, 
$v_1+av_2\in P$ for all $a\in\R_{\geq 0}$, and $v_1+av_2\to v$ as
$a\to\infty$.\qed

\begin{eg} Let \label{eg:P.compactification}
  $P\subset\R^2$ be the polyhedron of~\parref{eg:unbounded} and let
  $\sigma = \cU(P)$, the first quadrant.  Then 
  \[ N_\R(\sigma) =
  \R^2\djunion(\{\infty\}\times\R)\djunion(\R\times\{\infty\})\djunion\{(\infty,\infty)\} \] 
  as in~\eqref{eq:decomp.tropsigma}.  According
  to~\parref{prop:compactified.polyhedron}, under this identification we have
  \[ \bar P = P \djunion(\{\infty\}\times[1,\infty))
  \djunion([1,\infty)\times\{\infty\})\djunion\{(\infty,\infty)\}. \]
  See Figure~\ref{fig:barp}.

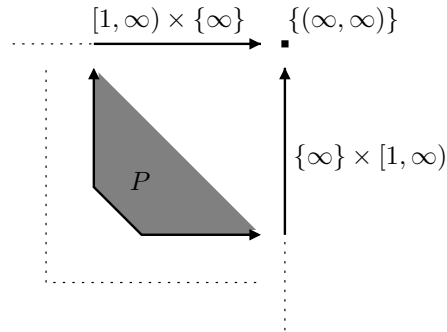
  \genericfig[ht]{barp}{The compactification of the polyhedron $P$
    of~\parref{eg:P.compactification}.} 
  
\end{eg}

\begin{rem} Let \label{rem:barP.in.NDelta}
  $\Delta$ be a pointed fan in $N_\R$.  Let $P\subset N$ be a pointed
  polyhedron, and suppose that its cone of unbounded directions $\sigma$
  is a cone of $\Delta$.  Then $N_\R(\sigma)\subset\Delta$ so $\bar P$ is
  naturally a subspace of $N_\R(\Delta)$; as $\bar P$ is compact,  
  $\bar P$ agrees with the closure of $P$ in $N_\R(\Delta)$.
  If $F\prec P$ then it follows from~\parref{lem:face.unbounded} that
  $\bar F\subset\bar P\subset N_\R(\Delta)$.
\end{rem}

\begin{rem} Let \label{rem:local.description}
  $\Delta$ be a pointed fan in $N_\R$.  A base for the topology of
  $N_\R(\Delta)$ is given by the interiors of the compactifications of the
  integral $\Gamma$-affine polyhedra whose cone of unbounded directions is a cone
  of $\Delta$.  See~\cite[Remark~3.4]{payne:analytification}.
\end{rem}

\paragraph This \label{par:barP.cap.Nmodsigma}
is a convenient place to mention the following construction, which will
come up later.  Let $\Delta$ be a pointed fan in $N_\R$ and fix
$\sigma\in\Delta$.  Let $N'_\R = N_\R/\spn(\sigma)$ and let 
$\Delta_\sigma$ be the pointed fan in $N'_\R$ whose cones are the images
of the cones $\tau\in\Delta$ such that $\sigma\prec\tau$.  Then
\[ N'_\R(\Delta_\sigma) = \Djunion_{\sigma\prec\tau} N_\R/\spn(\tau) \]
and so we have a natural inclusion 
$N'_\R(\Delta_\sigma)\inject N_\R(\Delta)$.

Let $P\subset N_\R$ be a pointed polyhedron with cone of unbounded
directions $\sigma'\in\Delta$, and suppose that $\sigma\prec\sigma'$.
Let $P' = \pi_\sigma(P)\subset N_\R'$, a
polyhedron in $N'_\R$.  It follows immediately
from~\parref{prop:compactified.polyhedron} that the compactification 
$\bar P'$ of $P'$ inside of 
$N'_\R(\pi_\sigma(\sigma'))\subset N'_\R(\Delta_\sigma)$ is equal to
$\bar P\cap N'_\R(\Delta_\sigma)$.

\section{A review of affinoid algebras}
\label{sec:affinoids}

\paragraph 
In this section we give a brief introduction of the theory of affinoid algebras
for the benefit of the reader who is not familiar with the language of rigid
analytic spaces.  We will only briefly mention the
global theory of rigid analytic spaces as it will not play a major role in
the sequel. 
Our main reference for all thing rigid-analytic is~\cite{bgr:nonarch},
although we refer the reader to~\cite{aws:2007} for an introduction
to the theory.  See also~\cite[\S3]{ekl:non_arch_amoebas} for an
introduction to the subject in the context of tropical geometry.

We fix the following notation for the rest of this paper:

\begin{notn} \label{notn:rigid}

  \begin{tabular}{ll}
    $K$ & is a field that is complete with respect to \\
    $\val: K\to\R\cup\{\infty\}$ & a nontrivial non-Archimedean valuation \\
    $|\cdot| = \exp(-\val(\cdot))$ & is the associated absolute value \\
    $\sO_K$  & is the valuation ring of $K$ \\
    $\fm_K\subset\sO_K$ & is the maximal ideal \\
    $k = \sO_K/\fm_K$   & is the residue field \\
    $\Gamma_K = \val(K^\times)$ & is the value group of $K$ \\
    $\Gamma = \sqrt{\Gamma_K} = \val(\bar K{}^\times)$ & 
    is the saturation of the value group.
  \end{tabular}\\
\end{notn}

Note that $\Gamma$ is divisible and hence dense in $\R$.  The base of the
exponential used in the definition of $|\cdot|$ can be any number greater
than $1$; we will use the natural exponential for concreteness.

\paragraph
The theory of rigid analytic spaces was invented by Tate in order to give
more structure to his non-Archimedean uniformization of elliptic curves
with split multiplicative reduction.  It closely parallels the theory of
complex analytic spaces, in that it exhibits many of the rigidity
characteristics of algebraic geometry while carrying a finer, analytic
topology.  We will try to emphasize the analogy with the theory of
varieties over a field.

\paragraph[Tate algebras]
Rigid spaces are modeled on closed subspaces of the
$p$-adic closed unit ball (or polydisc)
\[ \B^n_K(\bar K) = \{ (x_1,\ldots,x_n)\in\bar K{}^n~:~|x_i|\leq 1
\text{ for all } i \}, \]
which plays the same role as affine $n$-space in algebraic geometry.
(We use the \emph{closed} unit ball because the ring of analytic functions
on a ``compact'' space is well-behaved; in any case, $\B^n_K(\bar K)$ is
still open in the $p$-adic topology.)  
An infinite sum of elements in a complete non-Archimedean field converges
if and only if the absolute values of the summands approaches zero, so
one might expect that the holomorphic functions converging on this set
would correspond to the formal power series 
$\sum_\nu a_\nu x^\nu\in K\ps{x_1,\ldots,x_n}$ such that
$|a_\nu|\to 0$ as $|\nu|\to\infty$, where $\nu = (\nu_1,\ldots,\nu_n)$ and
$|\nu| = \nu_1+\cdots+\nu_n$.  This leads to the definition of the Tate
algebra, which plays the same role as a polynomial ring in algebraic geometry.

\begin{defn}
  The \emph{Tate algebra} in $n$ variables is the $K$-algebra
  \[ K\angles{x_1,\ldots,x_n} = 
  \bigg\{ \sum_\nu a_\nu x^\nu\in K\ps{x_1,\ldots,x_n}~:~
    |a_\nu| \to 0 \text{ as } |\nu|\to\infty \bigg\}. \]
\end{defn}

\begin{thm} The \label{thm:tate.algebras}
  Tate algebra $T_n = K\angles{x_1,\ldots,x_n}$ satisfies the
  following properties:
  \begin{enum}
  \item $T_n$ is an integral domain.  Moreover it is noetherian, regular,
    and a unique factorization domain.  For every maximal ideal $\fm$ of
    $T_n$ the local ring $(T_n)_\fm$ has dimension $n$ and its residue
    field $T_n/\fm$ is a finite extension of $K$.
  \item $T_n$ is a Jacobson ring: every prime ideal of $T_n$ is the
    intersection of the maximal ideals containing it.  In particular, if
    $\fa$ is an ideal of $T_n$ then an element of $T_n/\fa$ is nilpotent
    if and only if it is contained in every maximal ideal of $T_n/\fa$.
  \end{enum}
\end{thm}

Generally one expects $T_n$ to enjoy any ring-theoretic property
satisfied by $K[x_1,\ldots,x_n]$ that is not explicitly related to being of
finite type over a field.

\paragraph
The $K$-algebra homomorphisms from $K\angles{x_1,\ldots,x_n}$ to
$\bar K$ are in bijective correspondence with
$\B^n_K(\bar K)$ via $f\mapsto (f(x_1),\ldots,f(x_n))$.
Theorem~(\ref{thm:tate.algebras},ii) then allows us to view 
$K\angles{x_1,\ldots,x_n}$ as a function algebra on $\B^n_K(\bar K)$.
If we set $\B^n_K = \Max(K\angles{x_1,\ldots,x_n})$, the maximal spectrum
of the Tate algebra, then the map
\[ f\mapsto\ker(f)~:~ \B^n_K(\bar K)\To\B^n_K \]
is a surjection whose fibers are the $\Gal(K^\sep/K)$-orbits.

\begin{defn}
  A \emph{$K$-affinoid algebra} is a $K$-algebra that is isomorphic to a
  quotient of a Tate algebra.
\end{defn}

The maximal spectrum of an affinoid algebra is therefore a
\emph{Zariski}-closed subspace of a unit ball $\B^n_K$, defined by some
ideal $\fa\subset K\angles{x_1,\ldots,x_n}$.  (In general, a closed
analytic subspace of a rigid space should be thought of as being Zariski-closed.)
By~\parref{thm:tate.algebras}, an affinoid algebra is a Jacobson ring, and
therefore a reduced affinoid algebra $A$ is a function algebra on
the space $(\Max A)(\bar K) \coloneq \Hom_K(A,\bar K)$.  An affinoid
algebra is equipped with a canonical semi-norm%
\footnote{\label{fnt:norms}A \emph{semi-norm} $|\cdot|$ on a ring $A$ is a function
  $A\to\R_{\geq0}$ satisfying the ultrametric triangle inequality and
  such that $|1|=1$ and $|fg|\leq|f|\,|g|$ for all $f,g\in A$.  A
  semi-norm is called a \emph{norm} if $|f|=0$ if and only if $f=0$.}
$|\cdot|_{\sup}$, called the
\emph{supremum semi-norm}, defined by
\[ |f|_{\sup} = \sup_{\xi\in\Max A} |f(\xi)|
= \sup_{\xi\in(\Max A)(\bar K)} |f(\xi)|. \]
(Recall that there is a unique absolute value on any finite extension of
$K$ extending $|\cdot|$.)
If $A = K\angles{x_1,\ldots,x_n}$ and $f = \sum a_\nu x^\nu$ then
$|f|_{\sup} = \max_\nu |a_\nu|$; in this case $|\cdot|_{\sup}$ is called
the Gauss norm. A form of Gauss' lemma states that $|\cdot|_{\sup}$ is
multiplicative on $K\angles{x_1,\ldots,x_n}$
(i.e.\ $|fg|_{\sup}=|f|_{\sup}|g|_{\sup}$ for all  
$f,g\in K\angles{x_1,\ldots,x_n}$).
In general the supremum semi-norm may not be
multiplicative, but it is always power-multiplicative
(i.e.\ $|f^m|_{\sup} = |f|_{\sup}^m$ for all $m\geq 0$).  

\begin{thm}[Maximum Modulus Principle]
  Let $A$ be a $K$-affinoid algebra and let $f\in A$.  Then there exists
  $\xi\in\Max(A)$ such that $|f(\xi)| = |f|_{\sup}$.  In particular, 
  $\xi\mapsto |f(\xi)|$ is bounded and attains a maximum value on $\Max(A)$.
\end{thm}

It is clear that $|f|_{\sup} = 0$ if and only if $f$ is nilpotent.
If $A$ has no nilpotents then $A$ is complete and separated with
respect to $|\cdot|_{\sup}$~\cite[Theorem~6.2.4/1]{bgr:nonarch}.

\begin{rem} The \label{rem:tate.universal.prop}
  Tate algebra satisfies the following universal property (analogous
  to the universal property satisfied by a polynomial ring): if $A$ is a
  $K$-affinoid algebra (resp any $K$-Banach algebra%
  \footnote{A $K$-algebra that is complete and separated with respect to a
  norm extending the absolute value on $K$.}), then to give a
  $K$-algebra homomorphism $f: K\angles{x_1,\ldots,x_n}\to A$ is
  equivalent to choosing $a_1,\ldots,a_n\in A$ with $|a_i|_{\sup}\leq 1$
  (resp.\ such that $\{a_i^m\}_{m\geq0}$ is bounded).  That is, there exists a 
  unique homomorphism $f$ such $f(x_i) = a_i$.  
  See~\cite[Propositions~1.4.3/1 and~6.2.3/2]{bgr:nonarch}.
\end{rem}

\begin{eg}[Annuli] Let \label{eg:annuli}
  $r\in K^\times$ and let $\rho = |r|$.  Suppose that $\rho\leq1$.  Consider
  the affinoid algebra
  \[ A = K\angles{x,y}/(xy-r). \]
  Let $X\subset\B^2_K$ be its maximal spectrum and let
  $p_1: X\to\B^1_K$ be the projection onto the first factor.
  Then $p_1$ maps $X(\bar K)$ isomorphically onto 
  $\{\xi\in\B^1_K(\bar K)~:~|\xi|\geq\rho \}$.
  We call $X$ the \emph{annulus} of inner radius $\rho$ and outer radius
  $1$.  This is an example of a Laurent domain;
  see~\parref{par:laurent.domains}. 
\end{eg}

\paragraph
There is a notion of cofiber (tensor) product in the category of
$K$-affinoid algebras.  It is constructed as a completion of an ordinary
tensor product, but may be described more concretely as follows.  If
$A = K\angles{x_1,\ldots,x_n}/\fa$ and $B = K\angles{y_1,\ldots,y_m}/\fb$
are affinoid algebras then we set
\[ A\hat\tensor_K B = 
K\angles{x_1,\ldots,x_n,y_1,\ldots,y_m}/(\fa + \fb). \]
This $K$-algebra is visibly affinoid, and satisfies the universal
property of the cofiber product in the category of $K$-affinoid algebras.

If $A$ is $K$-affinoid then we set 
\[ A\angles{x_1,\ldots,x_n} = A\hat\tensor_K K\angles{x_1,\ldots,x_n}. \]

\paragraph
In order to put a sheaf of rings on the maximal spectrum $\Max(A)$ of an
affinoid algebra $A$, one has to understand the analogue of a
distinguished affine open subset.  As these will be quite a bit more
general than the complement of the zero locus of a regular function, it is
convenient to define an affinoid open subset by universal property:

\begin{defn*}
  Let $A$ be an affinoid algebra and let $U\subset\Max(A)$.  If there
  exists a homomorphism of affinoid algebras $f:A\to B$ such that $f^*$
  identifies $\Max(B)$ with $U$, and such that a homomorphism $g:A\to C$
  taking $\Max(C)$ into $U$ extends uniquely to 
  a homomorphism $B\to C$, then we say that $U = \Max(B)$ is a
  \emph{affinoid subdomain} of $\Max(A)$.
\end{defn*}

\paragraph If \label{par:laurent.domains}
the topology on $\Max(A)$ is going
to be ``analytic'', one would certainly hope that a subset of the form
$\{\xi:|f(\xi)|\leq 1\}$ would be an affinoid open for any $f\in A$.  In
fact we will want to consider the more general kind of analytic
subsets: 

\begin{defn*}
  Let $A$ be an affinoid algebra.  A \emph{Laurent domain}
  is a subset of $\Max(A)$ of the form
  \[ D(\bb f,\bg\inv) = \{ \xi\in\Max(A)~:~ |f_1(\xi)|,\ldots,|f_n(\xi)|\leq 1,~
  |g_1(\xi)|,\ldots,|g_m(\xi)|\geq 1 \} \]
  for some $\bb f = f_1,\ldots,f_n,\bg = g_1,\ldots,g_m\in A$. 
  If $m = 0$ we call $D(\bb f)$ a \emph{Weierstrass domain}.
  We set
  \[ A\angles{\bb f,\bg\inv} = A\angles{x_1,\ldots,x_n,y_1,\ldots,y_m}
  / (x_1-f_1,\ldots,x_n-f_n,~y_1g_1-1,\ldots,y_mg_m-1). \]
\end{defn*}

The following proposition follows almost immediately from the universal
properties of the Tate algebra and the completed tensor product. 

\begin{prop}
  The natural map $A\to A\angles{\bb f,\bg\inv}$ induces a bijection 
  \[ \Max(A\angles{\bb f,\bg\inv})\isom D(\bb f,\bg\inv)\subset\Max(A) \]
  that exhibits $D(\bb f,\bg\inv)$ as an affinoid subdomain of $\Max(A)$.
\end{prop}

\begin{eg} The \label{eg:laurent.domains}
  annulus of~\parref{eg:annuli} is by definition the Laurent domain
  $D((x/r)\inv)$ inside $\B^1_K$.  For $m\geq 0$ the Laurent domain
  $D((x^m/r)\inv)\subset\B^1_K$ is the set 
  $\{\xi:|\xi|\leq \rho^{1/m}\}$.  For any 
  $\mu\in|\bar K{}^\times|$ with $|\mu|<1$ we can find $m\geq 1$ such that
  $\mu^m\in|K^\times|$; thus we can define the annulus of
  inner radius $\mu$ and outer radius $1$ as above. 
  One can identify the coordinate ring of this annulus with the algebra
  \[ \bigg\{\sum_{i\in\Z} a_i x^i~:~|a_i|\to 0\text{ as } i\to\infty
  \text{ and }
  |a_i|\mu^{i}\to 0\text{ as } i\to -\infty\bigg\} \]
  The supremum norm is 
  $|\sum a_i x^i|_{\sup} = \sup\{|a_i|,\,|a_i|\mu^i~:~i\in\Z\}$.

  Consider the Weierstrass domain in $\B^1_K$
  \[ \B^1_K(\mu) \coloneq D(x^m/r) = \{ \xi\in\B^1_K~:~|\xi|\leq\mu \}. \]
  This is the ball of radius $\mu$; it is defined for every 
  $\mu\in|\bar K{}^\times|$ with $\mu\leq 1$.
  The coordinate ring of $\B^1_K(\mu)$ is naturally identified with
  the modified Tate algebra
  \[ T_{1,\mu}\coloneq \bigg\{\sum_{i\geq 0} a_i x^i ~:~ |a_i| \mu^i\to 0\bigg\}, \] 
  and the supremum  norm 
  is $|\sum a_i x^i|_{\sup} = \sup |a_i|\mu^i$.  In fact, for arbitrary
  $\mu\in|\bar K{}^\times|$ the algebra $T_{1,\mu}$ is $K$-affinoid with
  supremum norm $|\sum a_i x^i|_{\sup} = \sup |a_i|\mu^i$;
  we define $\B^1_K(\mu) = \Sp(T_{1,\mu})$ for any $\mu\in|\bar K{}^\times|$.

  The constructions above extend in an evident manner to define $n$-balls
  $\prod_{i=1}^n\B^1_K(\mu_i)$ and polyannuli of different radii, and to
  characterize their affinoid algebras and $\sup$ norms;
  see~\parref{eg:trop.B1} and~\parref{eg:trop.annulus}.
  There is a caveat however: if $\rho\in\R_{>0}$ but 
  $\rho\notin|\bar K{}^\times|$ then  
  $\{\xi\in\B^1_K:|\xi|\leq\rho\}$ is \emph{not} an affinoid subdomain of
  $\B^1_K$. 
\end{eg}

\paragraph
Here we give a brief sketch of the globalization procedure for rigid
spaces.  Let $A$ be an affinoid algebra and let $X = \Max(A)$.  A subset
$U\subset X$ is an \emph{admissible open subset} if it has a set-theoretic
covering $\{U_i\}$ by affinoid subdomains such that for any map of
affinoids $f: A\to B$ with $f^*(\Max(B))\subset U$ the cover 
$\{(f^*)\inv(U_i)\}$ of $\Max(B)$ has a finite subcover.  A set-theoretic covering
$\{U_i\}$ of an admissible open subset $U$ is an \emph{admissible cover}
provided that for any map of affinoids $f: A\to B$ such that
$f^*(\Max(B))\subset U$ the covering $\{(f^*)\inv(U_i)\}$ of $\Max(B)$ has
a refinement consisting of finitely many affinoid subdomains.  In
particular, any affinoid subdomain is an admissible open, and any cover by
finitely many affinoid subdomains is an admissible cover.

The admissible open subsets of $X$ form a Grothendieck topology whose
covers are the admissible covers.  Therefore $X$ has the structure of a
\emph{G-topological space}, i.e.\ a set endowed with a Grothendieck topology
on a collection of subsets.  (The point-set topology generated by
the affinoid opens induces the  $p$-adic topology on
$X$, which is totally disconnected and therefore too fine --- we want
$\B^n_K$ to be connected.)  The Tate acyclicity theorem roughly states
that there is a sheaf of rings $\sO_X$ on $X$ such that $\sO_X(\Max(B)) = B$
for every (admissible) affinoid open subset $\Max(B)\subset X$.  The
locally ringed G-topological space $\Max(A)$ is
called an \emph{affinoid space} and is denoted $\Sp(A)$.  A morphism of 
affinoid spaces is a morphism as locally ringed G-topological spaces.  Any
morphism $\Sp(B)\to\Sp(A)$ arises from a unique 
homomorphism $A\to B$. 

A \emph{rigid-analytic space} is a G-topological space (satisfying some
technical hypotheses) which admits an
admissible cover by affinoid spaces, and a morphism of rigid analytic 
spaces is a morphism in the category of locally ringed G-topological
spaces.

\begin{eg} The \label{eg:open.ball}
  \emph{rigid-analytic open unit ball} is the rigid space
  $\D^1_K = \bigcup_{\rho\in|\bar K{}^\times|,\,\rho < 1} \B^1_K(\rho)$,
  where $\B^1_K(\rho)$ is the ball of radius $\rho$ defined
  in~\parref{eg:laurent.domains}.  This cover is admissible by the
  maximum modulus principle.  More generally, we define
  $\D^1_K(\rho)$ for $\rho\in|\bar K{}^\times|$ in an evident manner.
\end{eg}

\begin{eg} Rigid-analytic \label{eg:affine.space}
  affine $m$-space is the rigid space
  $\A^{m,\an}_K = \bigcup_{\rho\in|\bar K{}^\times|} \B^m_K(\rho)$.
  Again this cover is admissible by the maximum modulus principle.
  We can define the rigid-analytic projective space
  $\bP^{m,\an}_K$ by gluing $m+1$ copies of $\A^{m+1,\an}_K$, but in fact
  $\bP^{m,\an}_K$ is covered by the $m+1$ closed unit balls 
  $\B^{m+1}_K\subset\A^{m+1,\an}_K$ since we can always normalize
  $[x_0:\cdots:x_m]$ so that $\max |x_i| = 1$.
\end{eg}

\paragraph
There is an analytification functor $X\mapsto X^\an$ from the category
of $K$-schemes locally of finite type to the category of rigid analytic
spaces.  This functor respects most notions common to both categories,
such as  open and closed immersions, finite, proper
(see~\ref{par:properness}), and 
projective morphisms, fiber products, etc.  Furthermore the set underlying
$X^\an$ is canonically identified with the set of closed points $|X|$
of $X$, and the completed local ring of $X^\an$ at $\xi\in|X|$ agrees
with $\hat\sO_{X,\xi}$.  The analytification can be defined by universal
property (as is the case over $\C$), but can also be described concretely
as follows.  The analytification of the algebraic affine space $\A^n_K$ is
the analytic affine space $\A^{n,\an}_K$ defined
above~\parref{eg:affine.space}.  If 
$X\subset\A^n_K$ is the closed subscheme cut out by a collection of
polynomials $f_1,\ldots,f_m\in K[x_1,\ldots,x_n]$ then 
$X^\an\subset\A^{n,\an}_K$ is the closed subspace defined by the same
polynomials $f_1,\ldots,f_m\in\Gamma(\A^{n,\an}_K,\sO)$.  Finally, if 
$X$ is an arbitrary locally-finite-type $K$-scheme covered by the affine
open subsets $\{U_i\}$ then $X^\an$ is obtained by pasting the
analytifications $U_i^\an$. 

The analogues of Serre's GAGA theorems hold in this
context~\cite{kiehl:gaga}.  In particular, any projective rigid space
(including any proper curve) has a unique algebraization.

\section{Kajiwara-Payne extended tropicalizations}  
\label{sec:extended.trop}

\paragraph In this section we set our notation regarding toric varieties
and review Kajiwara-Payne's construction of the tropicalization of a toric
variety over a 
non-Archimedean field \cite{kajiwara:tropical_toric_geometry,payne:analytification}.  
We refer the reader to~\cite{fulton:toric_varieties}
for a general reference for toric varieties.

\paragraph
Affine toric varieties are associated to integral pointed cones $\sigma$
in $N_\R$ as follows.  Recall~\parref{defn:S.sigma} that
$S_\sigma = \sigma^\vee\cap M$.  This is a finitely-generated monoid by
Gordan's Lemma~\cite[Proposition~1.2.1]{fulton:toric_varieties},
and $S_\sigma$ spans $\sigma^\vee$.  

\begin{notn*}
  Let $\sigma\subset N_\R$ be an integral pointed cone and let
  $K[S_\sigma]$ be the semigroup ring associated to $S_\sigma$.  For 
  $u\in S_\sigma$ we let $x^u$ denote the corresponding element of
  $K[S_\sigma]$, so 
  \[ K[S_\sigma] = \bigg\{ \sum_{u\in S_\sigma} a_u x^u~:~ a_u\in K
    \text{ and } a_u = 0\text{ for almost all } u \bigg\}. \]
  The affine toric variety over $K$ associated to
  $\sigma$ is denoted $X(\sigma)=\Spec(K[S_\sigma])$.
\end{notn*}

\begin{defn} Let \label{defn:trop.map}
  $\sigma\subset N_\R$ be an integral pointed cone.  We define the
  \emph{tropicalization map}
  \[ \trop:~ |X(\sigma)|\surject N_\Gamma(\sigma) \sptxt{by}
  \angles{u,\trop(\xi)}_\sigma = -\val(x^u(\xi)) \]
  with the notation in~\parref{defn:partial.compactification}
  and~\parref{defn:S.sigma}.  We also define 
  \[ \trop:~ X(\sigma)(\bar K) \surject N_\Gamma(\sigma) 
  \sptxt{by} \angles{u,\trop(\xi)}_\sigma = -\val(x^u(\xi)) \]
  and we define 
  \[ \trop:~ X(\sigma)^\an\surject N_\Gamma(\sigma) \]
  by identifying the set underlying $X(\sigma)^\an$ with $|X(\sigma)|$.
\end{defn}

The above definition makes sense because the residue field
$\kappa(\xi)$ of a closed point $\xi\in|X(\sigma)|$ is a finite extension of
$K$, and therefore inherits a unique valuation extending the one on $K$.
Note that $\trop: X(\sigma)(\bar K)\to N_\Gamma(\sigma)$ agrees with the
composition $X(\sigma)(\bar K)\to |X(\sigma)|\to N_\Gamma(\sigma)$.
See also~\cite[\S3]{payne:analytification}.

\begin{rem*}
  It may seem that we have changed our convention regarding the sign of $\trop(\xi)$
  from the definition given in the introduction.  However the latter definition
  rests on a choice of basis for $M$; choosing an appropriate basis, we
  recover the tropicalization map of the introduction.
  See~\parref{eg:trop.torus} below.
\end{rem*}

\begin{rem} Let \label{rem:tau.from.xi}
$\sigma\subset N_\R$ be a pointed cone and
let $\xi\in|X(\sigma)|$. 
According to~\parref{defn:trop.map}, we have
$\angles{u,\trop(\xi)}_\sigma\in\R$ if and 
only if $x^u(\xi)\neq 0$ for $u\in S_\sigma$. 
Remark~(\ref{rem:after.tropsigma},iv) then implies
that the set $\{u\in S_\sigma:x^u(\xi)\neq 0\}$ is equal to 
$S_\sigma\cap\tau^\perp$ where $\tau\prec\sigma$ is the face such that 
$\trop(\xi)\in N_\R/\spn(\tau)\subset N_\R(\sigma)$.
\end{rem}

\begin{eg} Let \label{eg:trop.torus}
  $e_1,\ldots,e_n$ be a basis for $M$ and let $e_1',\ldots,e_n'$ be
  the dual basis for $N$.
  Let $\sigma=\{0\}$, so $\sigma^\vee = \Z^n$ and
  \[ T = X(\sigma) = \Spec(K[x_1^{\pm 1},\ldots,x_n^{\pm 1}]) \cong\G_m^n \]
  is a torus, where $x_i \coloneq x^{-e_i}\in K[M]$.  
  Let $\xi = (\xi_1,\ldots,\xi_n)\in T(\bar K)$
  and let $\trop(\xi) = \sum_{i=1}^n v_i e_i'\in\R^n$.  According
  to~\parref{defn:trop.map} we have
  \[ v_i = \angles{e_i,\trop(\xi)} = \val(x_i(\xi)) = \val(\xi_i). \]
  Hence $\trop: T(\bar K)\to\R^n$ is simply
  \begin{equation} \label{eq:trop.map}
    (\xi_1,\ldots,\xi_n)\mapsto (\val(\xi_1),\ldots,\val(\xi_n)). 
  \end{equation}
  For general $\sigma$, the tropicalization map restricted to the dense
  torus in $X(\sigma)$ is of the above form: it is simply the
  vector of valuations of the coordinates of the point.
\end{eg}

\begin{eg} Choose \label{eg:trop.An}
  bases as in~\parref{eg:trop.torus}, and let
  $\sigma = \sum_{i=1}^n\R_{\geq 0} e_i'$.  Then
  $S_\sigma = \sum_{i=1}^n \Z_{\leq 0} e_i$, so 
  \[ X(\sigma) = \Spec(K[x_1,\ldots,x_n]) = \A^n_K \]
  is affine $n$-space.  The partial compactification $N_\R(\sigma)$ is
  identified with $(\R\cup\{\infty\})^n$ as in~\parref{eg:proj.plane}, and the
  tropicalization map 
  $\trop:\A^n_K(\bar K) = \bar K{}^n\to(\R\cup\{\infty\})^n$ is given
  by~\eqref{eq:trop.map} again, where now we allow
  $\val(\xi_i)$ to be $+\infty$ when $x_i(\xi)=0$

\end{eg}

\paragraph The \label{par:affine.functoriality} 
definition of the tropicalization map is functorial with respect to
equivariant morphisms of affine toric varieties, in the following sense.
Let $\phi: N'_\R\to N_\R$ be a homomorphism respecting a choice of
integral structures and carrying one integral pointed cone $\sigma'$ into
another integral pointed cone $\sigma$, as in~\parref{defn:S.sigma}.  
Then $\phi^*: M_\R\to M'_\R$
maps $S_\sigma$ into $S_{\sigma'}$, and therefore induces maps
$K[S_\sigma]\to K[S_{\sigma'}]$ and $X(\sigma')\to X(\sigma)$ making the
following diagram commute:
\begin{equation} \label{eq:affine.functoriality}
\xymatrix{
  {X(\sigma')(\bar K)} \ar[r] \ar[d] & {|X(\sigma')|} \ar[r]^\trop \ar[d] &
  {N_\Gamma(\sigma')} \ar[d]^\phi \\
  {X(\sigma)(\bar K)} \ar[r] & {|X(\sigma)|} \ar[r]^\trop &
  {N_\Gamma(\sigma).} 
}\end{equation}

\begin{notn}
  Let $\Delta$ be an integral pointed fan in $N_\R$.  
  We let $X(\Delta)$ denote the toric variety obtained by gluing the
  affine toric varieties $X(\sigma)$ along the open immersions
  $X(\tau)\inject X(\sigma)$ for $\tau\prec\sigma$. 
\end{notn}

\begin{defn}
  Let $\Delta$ be an integral pointed fan in $N_\R$.  We define the
  \emph{tropicalization map} 
  \[ \trop:~ |X(\Delta)| \surject N_\Gamma(\Delta) \sptxt{or}
  \trop:~ X(\Delta)(\bar K) \surject N_\Gamma(\Delta) \]
  by gluing the maps $\trop: |X(\sigma)|\to N_\Gamma(\sigma)$ 
  for $\sigma\in\Delta$ using the
  diagram~\eqref{eq:affine.functoriality}.  We define
  \[ \trop:~ X(\Delta)^\an\surject N_\Gamma(\Delta) \]
  by identifying the set underlying $X(\Delta)^\an$ with
  $|X(\Delta)|$. 
\end{defn}

\begin{eg} Let \label{eg:trop.projspace}
  $\Delta$ be the fan of~\parref{eg:proj.plane}.  The associated toric variety
  $X(\Delta)$ is isomorphic to the projective plane, and we can identify 
  $N_\R(\Delta)$ with 
  \[ \big((\R\cup\{\infty\})^3\setminus\{(\infty,\infty,\infty)\}\big)/\R \]
  where $\R$ acts by translations.  The tropicalization map is then given by
  \[ \trop\big[\xi_1\,:\,\xi_2\,:\,\xi_3\big] 
  = \big[\val(\xi_1)\,:\,\val(\xi_2)\,:\,\val(\xi_3)\big] \]
  with the sign conventions as in~\parref{eg:trop.torus}.
  A similar construction works in higher dimensions.
\end{eg}

The above definition of the tropicalization map is functorial with respect
to equivariant morphisms of toric varieties; 
cf.~\parref{par:affine.functoriality} and~\parref{par:functoriality.Delta}.
See~\cite[\S3]{payne:analytification} for more details.

\paragraph Let \label{par:torus.orbits}
$\Delta$ be an integral pointed fan in $N_\R$, and 
let $T = \Spec(K[M])\cong\bG_{m}^{ r}$ be the dense torus inside 
$X(\Delta)$.  There is a bijective correspondence~\cite[\S3.1]{fulton:toric_varieties}
\[ \sigma\mapsto T_\sigma~:~\Delta\isom
\{ \text{the $T$-orbits of }|X(\Delta)| \}, \]
defined as follows.  For $\sigma\in\Delta$ we define
$T_\sigma = \Spec(K[\sigma^\perp\cap M])$, with 
the inclusion $T_\sigma\inject X(\sigma)\subset X(\Delta)$ being given by
the homomorphism
\begin{equation} \label{eq:T.sigma.in.X.tau}
  K[S_\sigma]\To K[\sigma^\perp\cap M]\quad\text{ defined by }\quad
x^u \mapsto
\begin{cases}
  x^u \quad&\text{ if } u\in\sigma^\perp\cap S_\sigma \\
  0   \quad&\text{ otherwise.}
\end{cases} 
\end{equation}
In particular, if $\xi\in|T_\sigma|$ and $u\in S_\sigma$ then
$x^u(\xi)\neq 0$ if and only if $u\in\sigma^\perp\cap S_\sigma$,
so by~\parref{rem:tau.from.xi} we have 
$\trop(\xi)\in N_\R/\spn(\sigma)$.

In fact  more is true: $T_\sigma$ is a torus with lattice of characters
$\sigma^\perp\cap M$, and the dual of 
$(\sigma^\perp\cap M)\tensor_\Z\R = \sigma^\perp$ is 
$N_\R/\spn(\sigma)$, so (replacing $M_\R$ with $\sigma^\perp$ and
$N_\R$ with $N_\R/\spn(\sigma)$) we have a tropicalization map
$\trop: |T_\sigma|\to N_\R/\spn(\sigma)$.  An elementary compatibility
check shows that we have a commutative square  
\begin{equation} \label{eq:djunion.trop} \xymatrix{
  {\Djunion_{\sigma\in\Delta} |T_\sigma|} \ar[d]^\sim \ar[r]^(.4){\djunion\trop} &
  {\Djunion_{\sigma\in\Delta} N_\R/\spn(\sigma)} \ar[d]^\sim \\
  {|X(\Delta)|} \ar[r]^{\trop} & {N_\R(\Delta).} 
}\end{equation}
See~\parref{eg:orbits.trop.compat} for an example, and
see~\cite[\S3]{payne:analytification} for more details.

\subparagraph The \label{par:bar.T.sigma}
closure $\bar T_\sigma$ of $T_\sigma$ in $X(\Delta)$ is a $T$-equivariant
closed subvariety, and the map $\sigma\mapsto\bar T_\sigma$ is a bijection
between the cones of $\Delta$ and the $T$-equivariant closed subvarieties
of $X(\Delta)$.  The scheme $\bar T_\sigma$ is the toric variety with dense
torus $T_\sigma$ given by the fan $\Delta_\sigma$ in $N_\R'=N_\R/\spn(\sigma)$
defined in~\parref{par:barP.cap.Nmodsigma}.  If $\tau$ is a cone
of $\Delta$ such that $\sigma\prec\tau$ then the inclusion 
$\bar T_\sigma\cap X(\tau)\inject X(\tau)$ is explicitly given by the map
\begin{equation}\label{eq:V.sigma.in.X.tau}
K[S_\tau]\surject K[S_\tau\cap\sigma^\perp]\qquad\text{ defined by }
x^u \mapsto
\begin{cases}
  x^u &\quad\text{ if } u\in\sigma^\perp \\
  0   &\quad\text{ otherwise}.
\end{cases}
\end{equation}
We have $N_\R'(\Delta_\sigma) = \Djunion_{\sigma\prec\tau} N_\R/\spn(\tau)$
and the following square commutes:
\begin{equation}\label{eq:T.tau.compat}
\xymatrix{
  {|\bar T_\sigma|} \ar[d] \ar[r]^(.3){\trop} &
  {\Djunion_{\sigma\prec\tau} N_\R/\spn(\tau)} \ar[d] \\
  {|X(\Delta)|} \ar[r]^(.35){\trop} & {\Djunion_{\tau\in\Delta} N_\R/\spn(\tau).} 
}
\end{equation}

\begin{egsub} Let \label{eg:orbits.trop.compat}
  $\sigma_1$ be the fan of~\parref{eg:proj.plane}, so
  $X(\sigma)\cong\A^2_K$.  The decomposition~\eqref{eq:decomp.tropsigma}
  corresponds to the decomposition of $|\A^2_K|$ into $\G_m^2$-orbits
  \[ |\A^2_K| = |\G_m^2|\djunion|(\{0\}\times\G_m)|
  \djunion|(\G_m\times\{0\})|\djunion\{(0,0)\} \]
  under~\eqref{eq:djunion.trop}.  The invariant subvariety
  $\{0\}\times\A^1_K$ corresponds to the cone $\tau_1\prec\sigma_1$,
  and~\eqref{eq:T.tau.compat} expresses the compatibility
  of the tropicalization 
  $\trop: |\{0\}\times\A^1_K|\to\{\infty\}\times(\R\cup\{\infty\})$
  with $\trop:|\A^2_K|\to(\R\cup\{\infty\})^2$.

\end{egsub}

\section{Polyhedral subdomains of toric varieties}
\label{sec:polyhedral.subdomains}

\paragraph In this section we introduce a class of admissible affinoid
open subdomains of toric varieties which correspond
to polyhedral data inside its tropicalization.  These so-called polyhedral
subdomains are generalizations of the polytopal subdomains of affinoid
algebras introduced in~\cite[\S3]{ekl:non_arch_amoebas} and studied
in~\cite[\S4]{gubler:tropical}.  They enable a local study of the
tropicalization of a subvariety of a torus.

\begin{defn}
  Let $P\subset N_\R$ be an integral $\Gamma$-affine pointed polyhedron with cone of
  unbounded directions $\sigma = \cU(P)$.  The
  \emph{polyhedral subdomain} associated to $P$ is the set
  $U_P \coloneq \trop\inv(\bar P)\subset X(\sigma)^\an$.
\end{defn}

\begin{remsub}
  We will show~\parref{prop:polyhedral.affinoids} that $U_P$ is an
  affinoid open subdomain.  However for this to be
  true it is necessary that $P$ be integral $\Gamma$-affine:
  see~\parref{eg:trop.annulus} and the
  remark at the end of~\parref{eg:laurent.domains}.
\end{remsub}

\begin{remsub}
  If $P$ is a \emph{polytope} in $N_\R$ then $U_P$ is
  a polytopal subdomain as defined
  in~\cite[\S3]{ekl:non_arch_amoebas} and~\cite[Proposition~4.4]{gubler:tropical}.
  More accurately, Gubler's polytopal subdomain $U_P$ is the Berkovich
  space associated to the affinoid subdomain $U_P$
  (see~\ref{prop:polyhedral.affinoids}). 
  We choose to work with classical rigid spaces as opposed to Berkovich
  spaces simply because rigid spaces are more accessible and they suffice
  for our purposes.
\end{remsub}

The subset $U_P$ of $X(\sigma)^\an$ is in fact an affinoid
subdomain whose coordinate ring is the following affinoid algebra
(see~\ref{prop:polyhedral.affinoids}). 

\begin{defn} Let \label{defn:polyhedral.subdomain}
  $P\subset N_\R$ be an integral $\Gamma$-affine pointed polyhedron with cone of
  unbounded directions $\sigma = \cU(P)$.  We define
  \[ K\angles{U_P} = \bigg\{\sum_{u\in S_\sigma} a_u x^u ~:~
    a_u\in K,~ \val(a_u) - \angles{u,v}\to\infty\text{ for all } v\in P
  \bigg\} \]
  where the convergence (as always) is taken on the complements of finite subsets of
  $S_\sigma$.  If $A$ is any $K$-affinoid algebra we set
  $A\angles{U_P} = A\hat\tensor_K K\angles{U_P}$.  For 
  $f = \sum a_u x^u\in K\angles{U_P}$ we define
  \[ |f|_P = \sup_{\substack{u\in S_\sigma\\v\in P}} 
  |a_u| \exp(\angles{u,v}). \]
\end{defn}

\begin{rem} Let
  $\xi\in|U_P|$ and let $v = \trop(\xi)\in\bar P$.  Then 
  $\val(x^u(\xi))$ is by definition $-\angles{u,v}_\sigma$, so if
  $f = \sum a_u x^u\in K\angles{U_P}$ then $|a_u x^u(\xi)|\to 0$.  In other
  words, $K\angles{U_P}$ is precisely the ring of power series that appear to
  converge on all points of $U_P$.  This is made precise
  in~\parref{prop:polyhedral.affinoids}. 
\end{rem}

\begin{rem} Let \label{rem:check.on.vertices}
  $u\in S_\sigma$ and let $\face_u(P)\subset P$ be the associated
  face.  By definition of $\face_u(P)$, for any $v\in \face_u(P)$ we have
  $\angles{u,v} = \sup_{v'\in P}\angles{u,v'}$.
  Since any face contains a vertex, it follows that
  $f = \sum a_u x^u$ is in $K\angles{U_P}$ if and only if
  $\ord(a_u) - \angles{u,v}\to\infty$ for all $v\in \vertices(P)$.
  Moreover, by~\parref{prop:compactified.polyhedron} for any 
  $u\in S_\sigma$ the function $v\mapsto\angles{u,v}_\sigma$ on $\bar P$
  takes its   maximum on a vertex of $P$.  Therefore
  \begin{equation} \label{eq:norm.P.alt}
    |f|_P = \max_{\substack{u\in S_\sigma\\v\in\vertices(P)}}
    |a_u|\exp(\angles{u,v}) 
    = \sup_{\substack{u\in S_\sigma\\v\in\bar P}} 
    |a_u| \exp(\angles{u,v}_\sigma) < \infty. 
  \end{equation}
\end{rem}

\begin{rem} The \label{rem:K.UP.completion}
  function $|\cdot|_P$ defines a $K$-algebra norm (see
  footnote~\ref{fnt:norms}) on $K[S_\sigma]$ such that
  $K\angles{U_P}$ is the completion of $K[S_\sigma]$ with respect to this 
  norm.  In other words, $(K\angles{U_P},|\cdot|_P)$ is a $K$-Banach algebra.
\end{rem}

\begin{eg} Choose \label{eg:trop.B1}
  bases $e_1,\ldots,e_n$ for $M$ and $e_1',\ldots,e_n'$ for $N$ and let
  let $x_i = x^{-e_i}\in K[M]$ as in~\parref{eg:trop.torus}.  Let 
  $P = [0,\infty)^n\subset N_\R\cong\R^n$.
  Then $\sigma\coloneq\cU(P)=P$, $S_\sigma=\Z_{\leq0}^n$, 
  $\bar P = [0,\infty]^n\subset N_\R(\sigma)$, and
  $X(\sigma) = \A^n_K$, as in~\parref{eg:trop.An}
  and~\parref{eg:P.compactification}.  Hence 
  $U_P = \trop\inv(\bar P) = 
  \{ (\xi_1,\ldots,\xi_n)\in|\A^n_K|:\val(\xi_i)\geq 0\} = \B^n_K$.  
  This agrees with the fact that
  \[ K\angles{U_P} = \bigg\{\sum_{\nu\in\Z_{\geq 0}^n} a_\nu x^\nu~:~|a_\nu|\to 0\bigg\}
  = K\angles{x_1,\ldots,x_n} \]
  is a Tate algebra.
  More generally, if we take $P = \prod_{i=1}^n[r_i,\infty)$ for
  $r_1,\ldots,r_n\in\Gamma$ then 
  $U_P = \trop\inv\big(\prod_{i=1}^n[r_i,\infty]\big) =
  \prod_{i=1}^n\B^1_K(\exp(-r_i))\subset\A^n_K$ and 
  \[ K\angles{U_P} = \bigg\{\sum_{\nu\in\Z_{\geq 0}^n} 
  a_\nu x^\nu~:~|a_\nu|\exp(r_1\nu_1+\cdots+r_n\nu_n)\inv\to 0 \bigg\}. \]
  See~\parref{eg:laurent.domains}.  
\end{eg}

\begin{eg} With \label{eg:trop.annulus} 
  the notation in~\parref{eg:trop.B1}, let 
  $r_1,\ldots,r_n,s_1,\ldots,s_n\in\Gamma$ with $r_i\leq s_i$ and let
  $P = \prod_{i=1}^n[r_i,s_i]$.  This $P$ is a polytope, so 
  $\cU(P)=\{0\}$, $S_\sigma=M$, and $\bar P = P$.  The polytopal subdomain 
  $U_P = \trop\inv(P)$ is the polyannulus
  $\{(\xi_1,\ldots,\xi_n)\in|\G_m^{ n}|~:~r_i\leq
  \val(\xi_i)\leq s_i\}$.  The associated affinoid algebra is
  \[ K\angles{U_P} = \bigg\{ \sum_{\nu\in\Z_{\geq 0}^n} a_\nu x^\nu~:~
  |a_\nu| \mu^\nu\to 0\text{ for all } \mu\in\prod_{i=1}^n
  \{\exp(-r_i),\exp(-s_i)\} \bigg\} \]
  by~\parref{rem:check.on.vertices}.  See~\parref{eg:laurent.domains}.
\end{eg}

The following proposition is due to
Einsiedler,  Kapranov, and Lind~\cite[Proposition~3.1.8]{ekl:non_arch_amoebas}
and also to Gubler~\cite[Proposition~4.1]{gubler:tropical} when $P$ is a
polytope.  The general case is more difficult since $U_P$ is
not a Laurent domain in an easily identifiable affinoid subdomain
when $P$ is unbounded. 

\begin{prop} Let \label{prop:polyhedral.affinoids}
  $P\subset N_\R$ be an integral $\Gamma$-affine pointed polyhedron and let $\sigma$ be
  its cone of unbounded directions.
  \begin{enum}
  \item The ring $K\angles{U_P}$ is a $K$-affinoid algebra.
  \item The inclusion $K[S_\sigma]\inject K\angles{U_P}$ induces an open
    immersion $\Sp(K\angles{U_P})\inject X(\sigma)^\an$, and
  \item the image of this open immersion is equal to $U_P$.
    In particular, $U_P$ is an admissible affinoid open subset of $X(\sigma)^\an$.
  \item The supremum norm on $K\angles{U_P}$ agrees with $|\cdot|_P$, i.e.,
    for $f\in K\angles{U_P}$ we have
    \[ |f|_P = \sup_{\xi\in|U_P|} |f(\xi)| = \max_{\xi\in|U_P|} |f(\xi)|. \]
  \item The ring $K\angles{U_P}$ is a Cohen-Macaulay ring of dimension $g$. 
  \end{enum} 
\end{prop}

\pf 
Write $P = \bigcap_{i=1}^r\{ v\in N_\R:\angles{u_i,v}\leq b_i \}$ where
$b_i\in\Gamma$ and $\exp(b_i) = |x^{u_i}|_P$, and assume that
$u_1,\ldots,u_r$ generate $S_\sigma$.  Let
$\phi:\Z_{\geq0}^r\surject S_\sigma$ be the map 
$(\nu_1,\ldots,\nu_r)\mapsto\sum_{i=1}^r \nu_i u_i$; this induces a
surjective map $\phi:K[y_1,\ldots,y_r]\surject K[S_\sigma]$ given by 
$\phi(y^\nu) = x^{\phi(\nu)}$.  We identify $X(\sigma)$ with the image of
the associated closed immersion $X(\sigma)\inject\A^r_K$.  Letting
$\beta_i = \exp(b_i)$, we have 
$U_P = X(\sigma)^\an\cap\prod_{i=1}^r\B^1_K(\beta_i)$
because $\xi\in|U_P|$ if and only if 
\[ \val(y_i(\xi)) = \val(x^{u_i}(\xi)) = -\angles{u_i,\trop(\xi)}_\sigma \geq -b_i \]
(see~\ref{eg:trop.B1}). 
This proves that $U_P$ is an affinoid subdomain
of $X(\sigma)^\an$, and furthermore that $U_P$ is a closed subspace of
$\prod_{i=1}^r\B^1_K(\beta_i)$.  Let $\bb b = (b_1,\ldots,b_r)$ and let
\[ T_{r,\bb b} = \bigg\{\sum_{\nu\in\Z_{\geq 0}^r} a_\nu y^\nu
\in K\ps{y_1,\ldots,y_r}~:~ \val(a_\nu) - \nu\cdot\bb b\to\infty\bigg\}, \]
so $T_{r,\bb b}$ is an affinoid algebra with supremum norm
$|\sum a_\nu y^\nu|_{\sup} = \max |a_\nu|\beta^\nu$, and
$\Sp(T_{r,\bb b}) = \prod_{i=1}^r \B^1_K(\beta_i)\subset\A^{r,\an}_K$
(see~\ref{eg:trop.B1} and~\cite[\S6.1.5]{bgr:nonarch}).
The ideal defining $U_P\subset\prod_{i=1}^r\B^1_K(\beta_i)$ is the extension of 
$\fa = \ker(\phi)$; let $A = T_{r,\bb b}/\fa T_{r,\bb b}$, so 
$U_P = \Sp(A)$.  Since $|x^{mu_i}|_P = \beta_i^m$ for all $m\geq 0$, the
homomorphism $\phi$ extends uniquely to a homomorphism
$\phi: T_{r,\bb b}\to K\angles{U_P}$
(see~\ref{rem:tate.universal.prop}).  This homomorphism kills $\fa$ and
therefore descends to $\bar\phi: A\to K\angles{U_P}$.  We claim that
$\bar\phi$ is an isomorphism.

First we show that $\bar\phi$ is injective, i.e.\ 
$\ker(\phi) \subset \fa T_{r,\bb b}$.  Let $f = \sum a_\nu y^\nu\in\ker(\phi)$,
so for $u\in S_\sigma$ we have $\sum_{\phi(\nu)=u} a_\nu = 0$ (note
$\lim_{\nu\in\phi\inv(u)} |a_\nu| = 0$).  Setting
$f^u = \sum_{\phi(\nu)=u} a_\nu y^\nu$ we have 
$f = \sum_{u\in S_\sigma} f^u$ and $\phi(f^u) = 0$; since every
ideal in $T_{r,\bb b}$ is closed~\cite[Proposition~6.1.1/3]{bgr:nonarch} it
suffices to show that $f^u\in\fa T_{r,\bb b}$ for all $u$.  Thus we may
assume that $f = f^u$ for some $u\in S_\sigma$.   The sum
$f = \sum_{\phi(\nu)=u} a_\nu(y^\nu - y^{\nu_0})$
converges for fixed $\nu_0\in\phi\inv(u)$, so since
$y^\nu-y^{\nu_0}\in\fa$ for all $\nu$, we have $f\in\fa T_{r,\bb b}$ (again
since $\fa T_{r,\bb b}$ is closed).

Therefore $A\subset K\angles{U_P}$.  Next we claim that $|\cdot|_P$ restricts
to the supremum norm $|\cdot|_{\sup}$ on $A$.  For any vertex $v$ of $P$
the supremum norm on 
\[ K\angles{U_{\{v\}}} = 
\bigg\{\sum_{u\in M} a_u x^u~:~\val(a_u) - \angles{u,v}\to\infty\bigg\} \]
is $|\sum a_u x^u|_{\{v\}} = \sup_{u\in S_\sigma}(|a_u|\exp(\angles{u,v}))$
by~\cite[Proposition~4.1]{gubler:tropical} or using the fact that
$U_{\{v\}} = \Sp(K\angles{U_{\{v\}}})$ is a polyannulus;
see~\parref{eg:laurent.domains} and~\parref{eg:trop.annulus}.  Since
$U_{\{v\}}$ is an affinoid 
subdomain of $U_P$, for $f\in K\angles{U_P}$ we have
\[ |f|_{\sup} \geq \max_{v\in\vertices(P)}\sup_{\xi\in|U_{\{v\}}|} |f(\xi)| 
= \max_{v\in\vertices(P)} |f|_{\{v\}} = |f|_P  \]
where the last equality holds by~\eqref{eq:norm.P.alt}.  To prove the
inequality $|f|_{\sup}\leq|f|_P$ we must show that $|f(\xi)|\leq|f|_P$
for all $\xi\in|U_P|$.  For 
$f = \sum_{u\in S_\sigma} a_u x^u\in K\angles{U_P}$ we have
\[ |f(\xi)| \leq\sup_{u\in S_\sigma} |a_u|\,|x^u(\xi)|
= \sup_{u\in S_\sigma} |a_u|\exp(\angles{u,\trop(\xi)}_\sigma)
\leq |f|_P, \]
where the last inequality comes from~\eqref{eq:norm.P.alt}.

The reduced affinoid algebra $A$ is complete and separated with respect to 
$|\cdot|_{\sup} = |\cdot|_P$ by~\cite[Theorem~6.2.4/1]{bgr:nonarch}.  But
$A$ contains $K[S_\sigma]$ which is dense in
$K\angles{U_P}$ as noted in~\parref{rem:K.UP.completion}, so 
$A = K\angles{U_P}$.  This proves (i)--(iv).
By Hochster's Theorem \cite[Theorem~2.1]{cox:toric_varieties},
$X(\sigma)$ is Cohen-Macaulay of dimension $g$.  Assertion~(v) follows
because the completed local rings of $X(\sigma)$ and $X(\sigma)^\an$
agree; see~\cite[Appendix~A]{conrad:irredcomps} for details.\qed

\begin{rem} It \label{rem:U.P.immersion}
  follows from the proof of~\parref{prop:polyhedral.affinoids} that if
  $u_1,\ldots,u_r$ is a set of generators for $S_\sigma$ such that
  $P = \bigcap_{i=1}^r\{v\in N_\R:\angles{u_i,v}\leq b_i\}$ and
  $\beta_i\coloneq \exp(b_i) =|x^{u_i}|_P$ then we have a closed immersion
  \[ U_P\inject \B^1_K(\beta_1)\times\cdots\times\B^1_K(\beta_r) \]
  with the parameter on $\B^1_K(\beta_i)$ mapping to 
  $x^{u_i}\in K\angles{U_P}$.
\end{rem}

\paragraph
Let $P$ be an integral $\Gamma$-affine pointed polygon in $N_\R$.
The tropicalization map $\trop: X(\sigma)^\an\to N_\Gamma(\sigma)$
restricts to a surjective map $\trop: U_P\surject\bar P_\Gamma$.
If $\Delta$ is an integral pointed fan containing $\sigma=\cU(P)$  then
$X(\sigma)\subset X(\Delta)$ and hence we may identify $U_P$ with the
admissible affinoid open subset $\trop\inv(\bar P)$ in $X(\Delta)^\an$.

\section{Tropicalizations of embedded subspaces}
\label{sec:tropicalizations}

\paragraph In this section we define the tropicalizations of analytic and
algebraic subspaces of toric varieties.  The definitions are
self-contained and illustrated by some examples, but the reader may want to
consult~\cite{gathmann:tropical_geometry}, for instance, for an
introduction to the subject of tropical geometry.

\begin{defn} Fix \label{defn:tropicalization}
  an integral $\Gamma$-affine pointed polyhedron $P \subset N_\R$ with
  cone of unbounded 
  directions $\sigma$, and let $Y\subset U_P$ be the closed analytic
  subspace defined by an ideal $\fa\subset K\angles{U_P}$.  Define 
  \[ \Trop_\Gamma(Y) = \trop(Y) \subset \bar P_\Gamma, \]
  where $\trop: U_P\to\bar P_\Gamma$ is the tropicalization
  map~\parref{defn:trop.map}, and let $\Trop(Y)\subset\bar P$ be the
  closure of $\Trop_\Gamma(Y)$.  The set $\Trop(Y)$ is called the
  \emph{tropicalization} of $Y$ (as a subspace of $U_P$), and the map 
  $\trop = \trop|_{|Y|}: |Y|\to\Trop(Y)$ is again called the 
  \emph{tropicalization map}.

  If the ambient space is not clear from context we will write
  \[ \Trop_\Gamma(Y,\bar P_\Gamma)\sptxt{and}
  \Trop(Y, \bar P). \]

\end{defn}

\begin{rem}
  It is more natural to define $\Trop(Y)$ as the image of 
  the Berkovich analytic space $Y^\berk$ associated to $Y$ under the
  natural map $\trop: Y^\berk\to N_\R$, as in \cite[\S5]{gubler:tropical}
  or~\cite[Definition~4.1]{draisma:tropical}.  This approach has several
  advantages: for instance, there is no need to take closures, 
  the tropicalization inherits topological properties of the Berkovich
  space (e.g.\ connectedness), and there is no problem in the case of a trivial
  valuation.  We have chosen the above definition simply in order to avoid
  discussing Berkovich spaces. 
\end{rem}

\subparagraph With \label{par:remark.tropicalization}
the above definition it is clear that the tropicalization satisfies
the following functoriality property: let $\phi: N'\to N$ be a homomorphism
carrying an integral $\Gamma$-affine pointed polyhedron $P'\subset N'_\R$ into
another $P\subset N_\R$, so $\phi$ extends to a map
$\phi: N'_\R(\cU(P'))\to N_\R(\cU(P))$ taking $\bar P'$ into $\bar P$. 
If $Y'\subset U_{P'}$ and $Y\subset U_P$ are closed analytic
subvarieties such that the induced map $U_{P'}\to U_P$ takes
$Y'$ into $Y$, then $\phi(\Trop(Y'))\subset\Trop(Y)$.

For example, if $N' = N$, $P'\subset P$, and $Y' = Y\cap U_{P'}$, then
$\Trop(Y') = \Trop(Y)\cap\bar P{}'$.

It is also clear that the definition of $\Trop(Y)$ is insensitive to
finite extension of the base field $K$.  

\begin{rem} The \label{rem:tropicalization.valuations}
  definition of the tropicalization given above agrees with
  Gubler's tropicalization \cite{gubler:tropical}.
  The point of this section is to show that $\Trop(Y)$ is determined
  by the valuations of the coefficients of the power
  series vanishing on $Y$, thus showing that $\Trop(Y)$ can be effectively
  calculated and (in certain cases anyway) that it is a well-behaved
  convex-geometric object.  See~\parref{rem:trop.surj}. 
\end{rem}

\begin{eg} In \label{eg:dim1.trop}
  order to illustrate~\parref{rem:tropicalization.valuations}, we begin
  with the simplest example.  Let $N = M = \Z$ and let $P = [0,\infty)$, so
  $U_P = \B^1_K$ as in~\parref{eg:trop.B1}.  Let
  $x = x^{(-1)}\in K[M]$, the character corresponding to $-1\in M$, so the
  coordinate ring of $U_P$ is $K\angles x$ with the conventions
  in~\parref{eg:trop.B1}. Let
  $f = \sum_{u=0}^\infty a_u x^u\in K\angles x$ be nonzero and let
  $Y = V(f)\subset\B^1_K$ be the subspace defined by $f$.  Then
  $\Trop(Y)\subset\bar P=[0,\infty]$, and $v\in\Trop(Y)$ if and only if there
  exists $\xi\in\bar K$ such that $f(\xi)=0$ and $\val(\xi)=v$.
  For a particular choice of $\xi\in\bar K$, if there were some $u$ such
  that $\val(a_u\xi^u) < \val(a_{u'}\xi^{u'})$ for all $u'\neq u$ then by
  the ultrametric inequality, $\val(f(\xi)) = \val(a_u\xi^u)$, so
  $f(\xi)\neq 0$.  Writing $\val(\xi)=v$ and $\val(\xi')=v'$, this says that
  if there exists $u\geq 0$ such that $\val(a_u) + uv < \val(a_{u'})+u'v$ for
  all $u'\neq u$ then $v\notin\Trop(f)$.  In other words, a necessary
  condition for $v\in\Trop(f)$ is that there must exist at least
  \emph{two} numbers $u,u'$ such that 
  $\val(a_u)+uv = \val(a_{u'})+u'v\leq\val(a_{u''})+u''v$ for all
  $u''\geq0$.  (By the theorem of the Newton 
  polygon, or by~\parref{thm:main.tropical.thm.affinoid} below, this
  condition is also sufficient.  See~\parref{eg:newtonpoly1}.)

\end{eg}

\paragraph Let \label{par:height.graph}
$f = \sum_{u\in S_\sigma} a_u x^u\in K\angles{U_P}$ be nonzero and let
$\tau\prec\sigma$.  Define the 
\emph{height graph of $f$ with respect to $\tau$} to be
\[ H(f,\tau) = \{(u,\val(a_u))~:~u\in S_\sigma\cap\tau^\perp,~a_u\neq 0\}
\subset(S_\sigma\cap\tau^\perp)\times\R. \]
For $v\in N_\R/\spn(\tau)$ let
\begin{equation} \label{eq:vert.p.f}
\vertices_v(f) = \minset((-v,1),~H(f,\tau)) \subset H(f,\tau),
\end{equation}
where we regard $(N_\R/\spn(\tau))\times\R$ as a space of linear functionals on 
$\tau^\perp\times\R$.  This is a nonempty finite set by the definition of
$K\angles{U_P}$.  Define the \emph{initial form} of $f$ with respect
to $v\in N_\R/\spn(\tau)$ to be
\[ \inn_v(f) = \sum_{(u,\val(a_u))\in\vertices_v(f)} a_u\,x^u. \]
In other words, $\inn_v(f)$ is the (finite) sum of those monomials $a_u x^u$ such
that $u\in S_\sigma\cap\tau^\perp$ and 
\begin{equation} \label{eq:cuts.out.Trop}
  \val(a_u) - \angles{u,v} = 
  \min\{\val(a_{u'})-\angles{u',v}~:~u'\in S_\sigma\cap\tau^\perp \}. 
\end{equation}

\begin{eg} Continuing \label{eg:dim1.trop2}
  with~\parref{eg:dim1.trop}, we have
  $H(f,\{0\}) = \{(-u,\val(a_u)):a_u\neq 0\}\subset\Z_{\geq 0}\times\R$,
  and for $v\in[0,\infty)$ we have
  \[ \vertices_v(f) = \{ (-u,\val(a_u))~:~ \val(a_u)+uv
  \text{ is minimal among } \{\val(a_{u'})+u'v~:~u'\geq0 \}\}. \]
  Hence by the reasoning in~\parref{eg:dim1.trop}, if
  $v\in\Trop(Y)$ then $\#\vertices_v(f)\geq 2$, or equivalently
  $\inn_v(f)$ is not a monomial.  This is true in general:
\end{eg}

\paragraph Let \label{par:ultrametric}
$\xi\in|U_P|$, let $v = \trop(\xi)$, and suppose that $v\in N_\R/\spn(\tau)$,
i.e.\ that $\xi\in|T_\tau|$~\parref{par:torus.orbits}.
For $u\in S_\sigma\cap\tau^\perp$ and $a_u\in K$ we have
\[ \val(a_u\,x^u(\xi)) = \val(a_u) + \val(x^u(\xi)) = \val(a_u) - \angles{u,v}, \]
and for $u\in S_\sigma$ with $u\notin\tau^\perp$ we have $x^u(\xi) = 0$
by~\parref{rem:tau.from.xi}.
Therefore, the initial form $\inn_v(f)$ is the sum of those monomials
$a_ux^u$ with \emph{minimal} valuation when evaluated on $\xi$.  If 
$\inn_v(f) = a_u x^u$ is a \emph{monomial}
then $\val(f(\xi)) = \val(a_u x^u(\xi)) \neq \infty$ by the ultrametric
triangle inequality, so $f(\xi)\neq 0$.
Therefore, if $f(\xi) = 0$ then $\inn_v(f)$ is not a monomial.  It is a
fundamental fact that in an appropriate sense, the preceding condition is
\emph{sufficient} for there to exist a zero $\xi$ of $f$ with
$\trop(\xi)=v$.

\begin{thm} Let \label{thm:main.tropical.thm.affinoid}
  $P \subset N_\R$ be an integral $\Gamma$-affine pointed polyhedron with cone
  of unbounded directions $\sigma$ and let
  $Y\subset U_P$ be the closed analytic subspace defined by an ideal
  $\fa\subset K\angles{U_P}$.  Then
  \begin{enum}
  \item $\Trop(Y) = \{ v\in\bar P~:~\inn_v(f)\text{ is not a monomial for
      any } f\in\fa \}$.
  \item $\Trop_\Gamma(Y) = \Trop(Y)\cap\bar P_\Gamma$.
  \item For $\tau\prec\sigma$ we have 
    $\Trop(Y)\cap(N_\R/\spn(\tau)) = \Trop(Y\cap\bar T_\tau^\an)$.
  \end{enum}
\end{thm}

\subparagraph Part~(iii) \label{par:explain.iii}
requires some explanation.  Recall~\parref{par:bar.T.sigma}
that $\bar T_\tau$ is the affine toric variety defined by the image
$\sigma'$ of $\sigma$ in $N'_\R = N_\R/\spn(\tau)$, and that 
$\bar P\cap N'_\R(\sigma')$ is the compactification of 
the integral $\Gamma$-affine pointed polyhedron $P' = P\cap N'_\R$
by~\parref{par:barP.cap.Nmodsigma}.  Therefore
we may consider $Y\cap\bar T_\tau^\an$ as a closed analytic subspace of
$U_{P'} = U_P\cap\bar T_\tau^\an$, and consider its tropicalization inside
$N_\R'(\sigma')\subset N_\R(\sigma)$.  The statement of~(iii) is thus an
important compatibility that allows us to compute tropicalizations inside
toric varieties by reducing to the case of a torus.
See~\cite[Corollary~3.8]{payne:analytification}. 

\begin{remsub}
  Gubler has pointed out to us
  that~(\ref{thm:main.tropical.thm.affinoid},i) can be strengthened to
  show that the local structure of $\Trop(Y)$ at a point $v\in P^\circ$
  agrees with the tropicalization of the initial degeneration of $Y$ at
  $v$, as in the algebraic case.
\end{remsub}

\pf[of~\parref{thm:main.tropical.thm.affinoid}]
Let $C = \{ v\in\bar P:\inn_v(f)\text{ is not a monomial for any } f\in\fa \}$.
Since the condition for $\inn_v(f)$ not to be monomial is a closed
condition for fixed $f$, and since $\Trop_\Gamma(Y)\subset C$
by~\parref{par:ultrametric}, 
we have $\Trop(Y)\subset C$.  Since $\Q\subset\Gamma$ and
$C$ is defined by equations of the form~\eqref{eq:cuts.out.Trop}, the
closure of $C\cap\bar P_\Gamma$ is $C$.  Hence it suffices 
to show that if
\begin{equation} \tag{*}
  v\in\bar P_\Gamma \text{ is such that }
  \inn_v(f)\text{ is not a monomial for any } f\in \fa
\end{equation}
then $v\in\Trop_\Gamma(Y)$.  Let $v$ satisfy (*), and suppose that 
$v\in P$ (i.e.\ $v\in N_\R$).
After possibly making a finite extension of the ground field, we may translate
the problem by $-v$ to assume that $v = 0$.  Let $A$ be the ring
$K\angles{U_{\{0\}}} = \{ \sum_{u\in M} a_u x^u:|a_u|\to 0 \}$
and let $|\cdot|_0$ be the supremum norm on $A$, so
$|\sum a_u x^u|_0 = \max |a_u|$.
To show $0\in\Trop_\Gamma(Y)$, we must show that
$U_{\{0\}}\cap Y\neq\emptyset$, i.e. that $\fa$ does not generate the
unit ideal in $A$.  Suppose to the contrary that $1\in\fa A$.  Then there
exist $f_1,\ldots,f_r\in A$ and $g_1,\ldots,g_r\in\fa$ such that
$1 = \sum_{i=1}^r f_i g_i$.  Since $A$ is the completion of $K[M]$ under
$|\cdot|_0$, there exist sequences $\{f_{i,j}\}_{j\geq 1}\subset K[M]$
such that $\lim_{j\to\infty} f_{i,j} = f_i$.  Since 
$K[M]\subset K\angles{U_P}$ we have 
$h_j = \sum_{i=1}^r f_{i,j} g_i\in\fa$ and $\lim_{j\to\infty} h_j = 1$,
the limit always taken with respect to $|\cdot|_0$.  Writing
$h_j = \sum_{u\in S_\sigma} a_{j,u} x^u$ we have
\[ |h_j - 1|_0 = \max\{|a_{j,0}-1|,~|a_{j,u}|~:~u\in S_\sigma\setminus\{0\}\}. \]
Therefore $\val(a_{j,0}) = 0$ for $j\gg 0$ and $\val(a_{j,u})\to\infty$
uniformly for $u\neq 0$, so $\inn_0(h_j) = a_{j,0}$ for $j\gg 0$.  But
$h_j\in\fa$ and $\inn_0(h_j)$ is a monomial, a contradiction.

Now suppose that $v\in\bar P\cap(N_\R/\spn(\tau))$ for some $\tau\prec\sigma$.  
Let $N_\R' = N_\R/\spn(\tau)$ and $P' = P\cap N_\R'$,
and let $\sigma'$ be the image of $\sigma$ in $N_\R'$, so
$\bar P' = \bar P\cap N_\R'(\sigma')$ as in~\parref{par:explain.iii}.
The inclusion $U_{P'}\inject U_P$ corresponds to the surjection
$K\angles{U_P}\surject K\angles{U_{P'}}$ defined using the rule
\eqref{eq:V.sigma.in.X.tau}.  For $f\in K\angles{U_P}$ let $f'$ be its
image under this map.  By construction, 
we have $\inn_v(f) = \inn_v(f')$.  Therefore, the
above argument (as applied to $N_\R'$, $P'$, and $Y\cap U_{P'}$) shows
that there exists $\xi\in Y\cap U_{P'}$ such that $\ord(\xi) = v$.\qed

See~\parref{rem:trop.surj} for some remarks on the above proof.

\begin{eg} In \label{eg:newtonpoly1}
  this example we explain how~(\ref{thm:main.tropical.thm.affinoid},i)
  implies a large part of the theorem of the Newton polygon.  Let
  $f = \sum_{u=0}^\infty a_ux^u\in K\angles x$ as
  in~\parref{eg:dim1.trop}, where $x = x^{(-1)}$ still.  By definition the
  Newton polygon $\NP(f)$ is the lower convex hull of 
  $\{(u,\val(a_u)):a_u\neq 0\}$.  In order to maintain our sign
  conventions we let $\NP'(f)$ be the lower convex hull of $H(f,\{0\})$;
  this is the Newton polygon of $f$ flipped over the $y$-axis. It is an
  elementary exercise to   
  show that a line segment $\conv\{(-u,\val(a_u)),(-u',\val(a_{u'}))\}$ is
  contained in $\NP'(f)$ if and only if there exists $v\geq 0$ such that 
  \[ \{(-u,\val(a_u)),\,(-u',\val(a_{u'}))\}\subset\vertices_v(f), \]
  in which case $v$ is the slope of the line segment.  In particular, the
  line segments in $\NP'(f)$ are exactly the sets of the form 
  $\conv(\vertices_v(f))$ for $v\geq 0$.  See Figure~\ref{fig:newtonpoly}.
  Hence the elementary reasoning
  of~\parref{eg:dim1.trop2} translates into the easy direction of the
  theorem of the Newton polygon: if $f(\xi) = 0$ then 
  $\inn_v(f)$ is not a monomial, so $\#\vertices_v(f)\geq 2$, so
  $\conv(\vertices_v(f))$ is a line segment and hence $\val(\xi)$ is
  a slope of $\NP'(f)$.  Theorem~(\ref{thm:main.tropical.thm.affinoid},i)
  provides part of the hard direction: if $v$ is a slope of $\NP'(f)$
  then $\conv(\vertices_v(f))$ is a line segment, so
  $\#\vertices_v(f)\geq 2$, so $\inn_v(f)$ is not a monomial and hence
  there is at least one zero $\xi$ of $f$ such that $\val(\xi)=v$. 

  The full theorem of the Newton polygon (including information about
  multiplicities) is the one-dimensional case of the 
  intersection multiplicity formula~\parref{thm:multiplicity.ps};
  see~\parref{eg:multiplicity.ps.dim1}.  The multiplicity information is
  encoded in the Newton complex of $f$~\parref{par:newton.complex}.

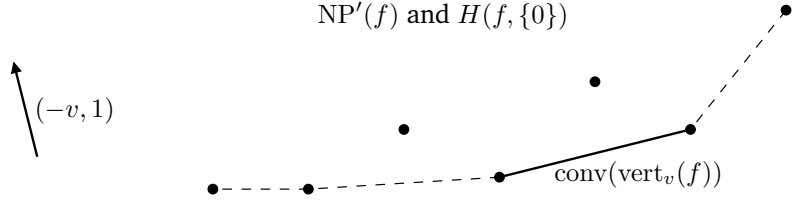
  \genericfig[ht]{newtonpoly}{This figure illustrates~\parref{eg:newtonpoly1}.  
    The dots represent $H(f,\{0\})$ and the dotted line represents
    $\NP'(f)$.  A number $v$ is a slope of $\NP'(f)$ if and 
    only if $\vertices_v(f)$ contains at least two elements, i.e.\ 
    $\conv(\vertices_v(f))$ is a line segment, in which case the slope of
    the line segment is $v$.}

\end{eg}

For another example see~\S\ref{sec:hypersurfaces}.
We now consider tropicalizations of (algebraic) subschemes of toric
varieties.  

\begin{defn}
  Let $\Delta$ be an integral pointed fan in $N_\R$ and let 
  $Y\subset X(\Delta)$ be a closed subscheme.  Define
  \[ \Trop_\Gamma(Y) = \trop(Y) \subset N_\Gamma(\Delta), \]
  and let $\Trop(Y)\subset N_\R(\Delta)$ be the closure  of
  $\Trop_\Gamma(Y)$.  The set $\Trop(Y)$ is called the
  \emph{tropicalization} of $Y$ (as a subscheme of $X(\Delta)$), and the
  map $\trop = \trop|_{|Y|}: |Y|\to\Trop(Y)$ is again called the
  \emph{tropicalization map}.

  As before if the ambient space is not clear from context we will write
  \[ \Trop_\Gamma(Y,N_\Gamma(\Delta)) \sptxt{and} 
  \Trop(Y, N_\R(\Delta)). \]
\end{defn}

\subparagraph It follows from the compatibility properties of the
tropicalization noted in~\parref{par:remark.tropicalization} that 
for any integral $\Gamma$-affine pointed polyhedron $P\subset N_\R$ whose cone
of unbounded directions $\sigma$ is contained in $\Delta$, we have
\[ \Trop(Y^\an\cap U_P) = \Trop(Y) \cap \bar P \sptxt{and} 
\Trop_\Gamma(Y^\an\cap U_P) = \Trop_\Gamma(Y) \cap \bar P_\Gamma. \]

While~\parref{thm:main.tropical.thm} does not follow formally
from~\parref{thm:main.tropical.thm.affinoid}, the proof carries over
verbatim.  

\begin{thm} Let \label{thm:main.tropical.thm}
  $\Delta$ be an integral pointed fan in $N_\R$ and let
  $Y\subset X(\Delta)$ be a closed subscheme.  Suppose that for
  $\sigma\in\Delta$ the closed subscheme 
  $Y\cap X(\sigma)\subset X(\sigma)$ is defined by the ideal 
  $\fa_\sigma\subset K[S_\sigma]$.  Then
  \begin{enum}
  \item $\Trop(Y) = \bigcup_{\sigma\in\Delta} 
    \{ v\in N_\R(\sigma)~:~\inn_v(f)\text{ is not a monomial for any } 
    f\in\fa_\sigma \}$.
  \item $\Trop_\Gamma(Y) = \Trop(Y)\cap\bar P_\Gamma$.
  \item For $\tau\prec\sigma$ we have 
    $\Trop(Y)\cap(N_\R/\spn(\tau)) = \Trop(Y\cap\bar T_\tau)$.
  \end{enum}
\end{thm}

\begin{rem} Theorem~\parref{thm:main.tropical.thm} \label{rem:trop.surj}
  is well-known
  (see~\cite{payne:fibers,speyer_sturmfels:tropical_grassmannian,ekl:non_arch_amoebas} for instance).
  The characterization of the tropicalization of an \emph{analytic}
  subspace of a torus~\parref{thm:main.tropical.thm.affinoid} has not
  appeared previously, although it is well-known to the experts.  
  The proof of~\parref{thm:main.tropical.thm.affinoid} closely resembles
  the proof of~\cite[Theorem~2.2.5]{ekl:non_arch_amoebas}, which 
  relates the Bieri-Groves set of a subvariety $Y$ of a torus
  with its tropicalization $\Trop(Y)$;
  in this sense the proof of~\cite[Theorem~2.2.5]{ekl:non_arch_amoebas} is
  the ``valuation-theoretic'' version of the proof
  of~\parref{thm:main.tropical.thm.affinoid}.  

  The main piece of machinery that is used in the proof
  of~\parref{thm:main.tropical.thm.affinoid} 
  is the interpretation~\parref{prop:polyhedral.affinoids} of the polyannulus 
  \[ U_{\{0\}} = \{(\xi_1,\ldots,\xi_n)\in|\G_m^{ n}| ~:~
  |\xi_1| = \cdots = |\xi_n| \} \]
  as the set of maximal ideals of the affinoid algebra
  $K\angles{U_{\{0\}}}$; then to show that there is a point $\xi$ of a
  subvariety $Y$ of $\G_m^{ n}$ inside $U_{\{0\}}$, i.e.\ such that
  $\trop(\xi)=0$, reduces to the \emph{algebraic} problem of showing that an
  ideal in $K\angles{U_{\{0\}}}$ is not the unit ideal.  This approach is
  quite standard once one is familiar with the theory of affinoid
  algebras, and is a compelling first application of the theory of
  rigid spaces to tropical geometry; in fact the author would argue
  that~\parref{thm:main.tropical.thm} is at heart a theorem in rigid
  analysis.  Another significant advantage of the rigid-analytic approach
  to this and other tropical problems is that the theory of rigid spaces
  has \emph{already been set up} to work over fields endowed with a
  non-discrete valuation, i.e., whose valuation ring is not noetherian.

  For a brief history of~\parref{thm:main.tropical.thm} as well as a
  stronger version,
  see~\cite{payne:fibers}.  See also~\cite[Theorem~4.2]{draisma:tropical}
  for a (different) proof of~\parref{thm:main.tropical.thm}
  that uses affinoid algebras.
\end{rem}

\section{Tropical hypersurfaces and the Newton complex}
\label{sec:hypersurfaces}

\paragraph 
When a closed analytic subspace $Y$ of a polyhedral subdomain of a toric
variety is defined by a single equation $f$, its tropicalization comes
equipped with extra combinatorial structures (as is well-known in the
algebraic case): the set $\Trop(Y)$ is the
support of a polyhedral complex, which is ``dual'' to the
so-called Newton complex $\New(f)$ also naturally associated to
$f$.  The Newton complex should be regarded as recording the multiplicity
information missing from $\Trop(Y)$.  These extra structures render
$\Trop(Y)$ easily computable in terms of $f$, and will later be used to
compute a local intersection multiplicity formula for rigid-analytic
complete intersections~\parref{thm:multiplicity.ps}.  The difficulty in setting up the
theory is showing that these complexes are in fact
\emph{finite}, so we begin with the key finiteness result.

\begin{notn*}
  Let $P\subset N_\R$ be an integral $\Gamma$-affine pointed polyhedron and let
  $f\in K\angles{U_P}$ be nonzero.  For any subset $\Sigma\subset P$ we define
  \[ \vertices_\Sigma(f) = \bigcup_{v\in\Sigma} \vertices_v(f), \]
  where $\vertices_v(f)$ is defined in~\eqref{eq:vert.p.f}.
\end{notn*}

\begin{lem} Let \label{lem:inu.finiteness}
  $P\subset N_\R$ be an integral $\Gamma$-affine pointed polyhedron and let
  $f\in K\angles{U_P}$ be nonzero.
  \begin{enum}
  \item The set $\vertices_P(f)$ is finite.
  \item There exists $\epsilon > 0$ such that for all $f'\in K\angles{U_P}$
    with $|f-f'|_P<\epsilon$ we have $\vertices_P(f) = \vertices_P(f')$.
  \end{enum}
\end{lem}

\pf Let $\sigma = \cU(P)$ and write $f = \sum_{u\in S_\sigma} a_u x^u$.
For fixed $v\in P$ we have $\val(a_u) - \angles{u,v}\to\infty$ by
definition; let 
$m(v) = \min_{u\in S_\sigma}\{\val(a_u)-\angles{u,v}\}$, so
\[ \vertices_v(f) = \{(u,\val(a_u))~:~\val(a_u)-\angles{u,v} = m(v) \} \]
by~\eqref{eq:cuts.out.Trop}, and hence
\[ \vertices_P(f) = \{(u,\val(a_u))~:~\val(a_u)-\angles{u,v} = m(v)
\text{ for some } v\in P \}. \]
Let $F_b$ be the union of the bounded faces of $P$, so $P = F_b + \sigma$
by~\parref{lem:bounded.faces}.  Let $\fa\subset K[S_\sigma]$ be the ideal
generated by $\{x^u:a_u\neq 0\}$, so since $K[S_\sigma]$ is noetherian,
there exist $u_1,\ldots,u_r\in S_\sigma$ such that
$\fa = (x^{u_1},\ldots,x^{u_r})$.  Let
\[ \alpha = \max\{\ord(a_{u_i})-\angles{u_i,v}~:~i=1,\ldots,r,~v\in F_b\}. \]
Let $v\in P$ and write
$v = v' + v''$ for $v'\in F_b$ and $v''\in\sigma$.  Note that for any
$u_0\in S_\sigma$ we have
\begin{equation} \label{eq:m.v.bound} \begin{split}
  m(v) &= \min_{u\in S_\sigma} \{\val(a_u) - \angles{u,v'} - \angles{u,v''}\}\\
  &\leq \min_{u\in S_\sigma}\{\val(a_u)-\angles{u,v'}\} + 
  \min_{u\in S_\sigma} \{-\angles{u,v''}\} \\
  &\leq m(v') - \angles{u_0,v''} \leq \alpha - \angles{u_0,v''}. 
\end{split}\end{equation}

Let $v_1,\ldots,v_s$ be the vertices of $P$, so
$F_b\subset\conv\{v_1,\ldots,v_s\}$.  Let
\[ \Psi = \{ u\in S_\sigma~:~\ord(a_u) - \angles{u,v_i}\leq\alpha
\text{ for some } i = 1,\ldots,s \}, \]
so $\Psi$ is a finite set.  We will show that
$\vertices_P(f)\subset\Psi$.  Fix $u\in S_\sigma\setminus\Psi$ and assume
that $a_u\neq 0$ (since otherwise $(u,\val(a_u))\notin\vertices_P(f)$ by
definition), so $x^u\in\fa$.  Fix $i_0\in\{1,\ldots,s\}$ such that 
$\angles{u,v_{i_0}} = \max_{i=1,\ldots,s} \angles{u,v_i}$.  
Let $v\in F_b$ and write
$v = \sum_{i=1}^s t_i v_i$ with $0\leq t_i\leq 1$ and 
$\sum_{i=1}^s t_i=1$.  Then we have 
\[ \val(a_u) - \angles{u,v} = \val(a_u) - \sum_{i=1}^s t_i\angles{u,v_i} 
\geq \val(a_u) - \sum_{i=1}^s t_i\angles{u,v_{i_0}} 
= \val(a_u) - \angles{u,v_{i_0}} > \alpha \]
where the final inequality holds because $u\notin\Psi$.  Now let
$v\in P$ be arbitrary, and write $v = v' + v''$ for $v'\in F_b$
and $v''\in\sigma$.  Since $x^u$ is contained in the monomial ideal $\fa$
we can write $u = u_{j_0} + u'$ for some $j_0=1,\ldots,r$ and 
$u'\in S_\sigma$.  We calculate
\[\begin{split}
  \val(a_u) - \angles{u,v} &= (\val(a_u) - \angles{u,v'}) - \angles{u,v''} \\
  &> \alpha - \angles{u,v''} = \alpha - \angles{u_{j_0},v''} - \angles{u',v''} 
  \geq \alpha - \angles{u_{j_0},v''}
\end{split}\]
since $v'\in F_b$ and $\angles{u',v''}\leq 0$.  But
$m(v)\leq\alpha-\angles{u_{j_0},v''}$ by~\eqref{eq:m.v.bound}, so
$\val(a_u)-\angles{u,v} > m(v)$ for all $v\in P$ and hence
$u\notin\vertices_P(f)$.  This proves~(i).

For $\epsilon > 0$ we let 
$f' = \sum_{u\in S_\sigma} a_u' x^u\in K\angles{U_P}$ denote a generic
power series satisfying $|f-f'|_P < \epsilon$; for such an $f'$ we define
$m'(v)$ as above.  Note that $|f-f'|_P < \epsilon$ if and only if
\[ \min\{\val(a_u-a_u')-\angles{u,v}~:~u\in S_\sigma,~v\in P\} >
-\log(\epsilon). \]

\textbf{Step 1.} First we choose $\epsilon$ small enough that
$\val(a_u) = \val(a_u')$ for all $u\in\Psi$.  Since
$\vertices_P(f)\subset\Psi$ we have 
$m(v) = \min_{u\in\Psi}\{\val(a_u)-\angles{u,v}\}$ for all $v\in P$, and
hence $m'(v)\leq m(v)$ for all $v\in P$.

\textbf{Step 2.} Decreasing $\epsilon$ if necessary we may assume that 
$-\log(\epsilon) > \alpha$.  We claim that $\vertices_P(f')\subset\Psi$.
Fix $u\in S_\sigma\setminus\Psi$.  For all $i=1,\ldots,s$ we have
\[ \val(a_u') - \angles{u,v_i} \geq
\min\{\val(a_u) - \angles{u,v_i},~\val(a_u-a_u')-\angles{u,v_i}\} >
\alpha. \]
It follows that $\val(a_u')-\angles{u,v} > \alpha$ for all $v\in F_b$ as
above.  Now let $v\in P$ be arbitrary, and write $v = v'+v''$ for
$v'\in F_b$ and $v''\in\sigma$, so
\[ m'(v)\leq m(v) \leq \alpha - \angles{u,v''}
< (\val(a_u') - \angles{u,v'}) - \angles{u,v''} = \val(a_u') -
\angles{u,v}, \]
where we used~\eqref{eq:m.v.bound} for the second inequality.  It follows
that $u\notin\vertices_P(f')$.

\textbf{Step 3.} We claim that we can decrease $\epsilon$ further so that
$\vertices_P(f)\subset\vertices_P(f')$.  Choose $w_1,\ldots,w_t\in P$ such
that $\vertices_P(f) = \bigcup_{i=1}^t \vertices_{w_i}(f)$, and suppose
for the moment that $m(w_i) = m'(w_i)$ for $i=1,\ldots,t$.  Let
$u\in\vertices_P(f)$, and suppose that
$(u,\val(a_u))\in\vertices_{w_{i_0}}(f)$, so
\[ m'(w_{i_0}) = m(w_{i_0}) = \val(a_u) - \angles{u,w_{i_0}} = 
\val(a_u') - \angles{u,w_{i_0}}, \]
and hence $u\in\vertices_P(f')$.  Therefore it suffices to show that 
$m(w_i)=m'(w_i)$ for $i=1,\ldots,t$ after potentially shrinking $\epsilon$
again.  In fact, if $-\log(\epsilon) > \max\{m(w_1),\ldots,m(w_t)\}$ then
for $u\in S_\sigma$ and $i=1,\ldots,t$ we have
\[ \val(a_u') - \angles{u,w_i} \geq 
\min\{\val(a_u)-\angles{u,w_i},~\val(a_u-a_u')-\angles{u,w_i}\} \geq
m(w_i), \]
which shows that $m'(w_i)\geq m(w_i)$.

\textbf{Step 4.} Finally we claim that
$\vertices_P(f')\subset\vertices_P(f)$ with the above conditions on
$\epsilon$.  We are done if $\Psi = \vertices_P(f)$, so
assume that there exists $u_0\in\Psi\setminus\vertices_P(f)$.  Let 
$v\in P$, and choose $u\in\vertices_P(f)$ such that 
$\val(a_u)-\angles{u,v} = m(v) < \val(a_{u_0})-\angles{u_0,v}$.  Since
$\val(a_u) = \val(a_u')$ and $\val(a_{u_0}) = \val(a_{u_0}')$ we have
$m'(v)\leq \val(a_u') -  \angles{u,v} < \val(a_{u_0}') - \angles{u_0,v}$.
Since $v$ was arbitrary, this proves that $u_0\notin\vertices_P(f')$.\qed

\paragraph We move on to tropicalizations of hypersurfaces.  For
convenience we use the following piece of

\begin{notn*}
  Let $P\subset N_\R$ be an integral $\Gamma$-affine pointed polyhedron (resp. let
  $\sigma\subset N_\R$ be an integral pointed cone) and let
  $f\in K\angles{U_P}$ (resp. $f\in K[S_\sigma]$).  We denote the closed
  analytic subspace of $U_P$ (resp. closed subscheme of $X(\sigma)$) defined by $f$ by
  $V(f)$, and we set 
  \[ \Trop_\Gamma(f) \coloneq \Trop_\Gamma(V(f)) \sptxt{and}
  \Trop(f) \coloneq \Trop(V(f)). \]
  As before if the ambient space is not clear from context we write
  \[ \Trop_\Gamma(f,\bar P_\Gamma),\qquad\Trop(f,\bar P),\qquad
  \Trop_\Gamma(f,N_\Gamma(\sigma)),\sptxt{and}\Trop(f,N_\R(\sigma)). \]
\end{notn*}

It is clear that if $\sigma = \cU(P)$ and $f\in K[S_\sigma]$ then 
$V(f)^\an\cap U_P = V(f|_{U_P})$ and 
$\Trop(f)\cap\bar P = \Trop(f|_{U_P})$, so the ambiguity in the notation
should not cause confusion.  

We note that $\Trop(f)$ is determined by $f$ in the way one might expect: 

\begin{lem} Let \label{lem:trop.hypersurface}
  $P\subset N_\R$ be an integral $\Gamma$-affine pointed polyhedron (resp. let
  $\sigma\subset N_\R$ be an integral pointed cone) and let 
  $f\in K\angles{U_P}$ (resp. $f\in K[S_\sigma]$) be nonzero.  Then
  \[ \Trop(f) = \{ v\in\bar P\text{ (resp. $v\in N_\R(\sigma)$) }~:~
  \inn_v(f)\text{ is not a monomial }\}. \]
\end{lem}

\pf The algebraic version follows from the rigid-analytic version, so
assume $f\in K\angles{U_P}$.  We must show that if $v\in\bar P_\Gamma$ and
$\inn_v(f)$ is not a monomial then $v\in\trop(V(f))$.  Reducing to the
case $v=0$ as in the proof of~\parref{thm:main.tropical.thm.affinoid}, we
would like to show that $f$ is a unit in $A = K\angles{U_{\{0\}}}$ if and only
if $\inn_0(f)$ is a monomial.  Let 
$\mathring A = \{ g\in A:|g|_0\leq 1\}$ and
$\check A = \{ g\in\mathring A:|g|_0 < 1\}$, and let
$\td A = \mathring A/\check A\cong K[M]$.  By scaling we may assume
$|f|_0=1$, so its residue $\td f\in\td A$ is nonzero.  If 
$\td f\in\td A{}^\times$ then there exists $g\in\mathring A$ such that
$\td f\td g = 1$, so $fg = 1-h$ for $h\in\check A$.  Since
$|h|_0 < 1$ we have $\lim_{m\to\infty} h^m = 0$, so
$fg\sum_{m=0}^\infty h^m = 1$ and hence $f\in\mathring A{}^\times$.%
\footnote{This is a general fact about affinoid algebras.
  See~\cite[Proposition~1.2.5/8]{bgr:nonarch}.} 
But $\inn_0(f)$ is a monomial if and only if $\td f$ is a monomial, in
which case $\td f\in\td A{}^\times$, so $f\in A^\times$.\qed

\paragraph[The polyhedral complex structure on $\Trop(f)$]
Let \label{par:trop.poly.cplx.1}
$\sigma\subset N_\R$ be an integral pointed cone and let 
$f\in K[S_\sigma]$ be nonzero.  Let $\tau\prec\sigma$, let 
$N_\R' = N_\R/\spn(\tau)$, and let $H(f,\tau)$ be the height graph of
$f$, where we are using the notation of~\parref{par:height.graph}.  Assume
that the image of $f$ in $K[S_\sigma\cap\tau^\perp]$ is nonzero. 
By~\parref{lem:trop.hypersurface} we have
$\#\vertices_v(f)\geq 2$ (i.e.\ $\inn_v(f)$ is not a monomial) if and only
if $v\in\Trop(f)$.  For $v\in \Trop(f)\cap N_\R'$ define 
\[ \gamma_v = \{ v'\in N_\R'~:~\vertices_{v'}(f)\supset\vertices_v(f) \}. \]
It is standard (see for
instance~\cite[Theorem~2.1.1]{ekl:non_arch_amoebas}) 
that $\{\gamma_v:v\in\Trop(f)\cap N_\R'\}$ is an 
integral $\Gamma$-affine polyhedral complex in $N_\R'$ 
of pure codimension $1$ (that is, all maximal cells have dimension 
$\dim_\R(N_\R')-1$),
and since $v\in\gamma_v$ the support of this complex is exactly
$\Trop(f)\cap N_\R'$.   We will write $\Trop(f)\cap N_\R'$ to denote the
polyhedral complex as well as its support.  To summarize:

\begin{prop*}
  Let $\sigma\subset N_\R$ be an integral pointed cone and let 
  $f\in K[S_\sigma]$; let $\tau\prec\sigma$ and assume
  that the image of $f$ in $K[S_\sigma\cap\tau^\perp]$ is nonzero.  Then
  $\Trop(f)\cap(N_\R/\spn(\tau))$ is the support of a natural polyhedral
  complex in $N_\R/\spn(\tau)$ of pure codimension $1$.
\end{prop*}

\begin{egsub} To \label{eg:tropical.line.1}
  illustrate, let $N = M = \Z^2$, let $\sigma=\{0\}$, and let
  $x = x^{(-1,0)},y = x^{(0,-1)}\in K[M]$.
  Let $\lambda\in K$ have valuation $1$ and let
  $f(x,y) = x+y+\lambda\in K[M]$.  Then
  \[\begin{split} 
    R_1\coloneq (1,1) + \R_{>0}(-1,-1) &= \{ v\in\R^2~:~\inn_v(f) = x + y \}  \\ 
    R_2\coloneq (1,1) + \R_{>0}(0,1) &= \{ v\in\R^2~:~\inn_v(f) = x + \lambda \} \\
    R_3\coloneq (1,1) + \R_{>0}(1,0) &= \{ v\in\R^2~:~\inn_v(f) = y + \lambda \} \\
    \{(1,1)\} &= \{ v\in\R^2~:~\inn_v(f) = x + y + \lambda \}. 
  \end{split}\]
  Each $R_i$ is an open ray in $\R^2$, and 
  $\bar R_i = R_i\cup\{(1,1)\} = \gamma_v$ for any
  $v\in R_i = \relint(\bar R_i)$.  The vertex $(1,1)$ is equal to
  $\gamma_{\{(1,1)\}}$. 
  Hence $\Trop(f) = \{\bar R_1,\bar R_2,\bar R_3,\{(1,1)\}\}$ as a
  polyhedral complex.  See Figure~\ref{fig:trop_new}.
\end{egsub}

\subparagraph Now \label{par:trop.poly.cplx.2}
let $P\subset N_\R$ be an integral $\Gamma$-affine pointed polyhedron with cone of
unbounded directions $\sigma$, let $\tau\prec\sigma$, let 
$N_\R' = N_\R/\spn(\tau)$, and let $P' = \bar P\cap N_\R'$.  Let
$f = \sum a_u x^u\in K\angles{U_P}$ have nonzero image in
$K\angles{U_{P'}}$.  As above we have  
$\vertices_{P'}(f)\subset H(f,\tau)$ and
by~\parref{lem:inu.finiteness} (as applied to the image of
$f$ in $K\angles{U_{P'}}$) the set $\vertices_{P'}(f)$ is finite.  
Define 
\[ f' = \sum \{a_u x^u~:~(u,\val(a_u))\in\vertices_{P'}(f)\} \in
K[S_\sigma\cap\tau^\perp], \]
so $f'$ is a Laurent polynomial such that 
$\vertices_v(f') = \vertices_v(f)$ for all $v\in P'$.  Again since
$\inn_v(f)$ is a monomial if and only if $\#\vertices_v(f)=1$,
by~\parref{lem:trop.hypersurface} we have 
$\Trop(f)\cap N_\R' = \Trop(f')\cap P'$.  Therefore,

\begin{prop*}
  Let $P\subset N_\R$ be an integral $\Gamma$-affine pointed polyhedron with cone of
  unbounded directions $\sigma$, let $\tau\prec\sigma$, and let 
  $P' = \bar P\cap(N_\R/\spn(\tau))$.  Let
  $f = \sum a_u x^u\in K\angles{U_P}$ have nonzero image in
  $K\angles{U_{P'}}$.  Then $\Trop(f)\cap(N_\R/\spn(\tau))$ is the
  intersection of the support of a pure-codimension-$1$ polyhedral complex
  in $N_\R/\spn(\tau)$ with the polyhedron $P'$.
\end{prop*}

\begin{remsub}
  In this case $\Trop(f)\cap N_\R'$ is the support of the
  polyhedral-complex-theoretic
  intersection~(\ref{defn:polyhedral.complexes},iv) of the complex  
  $\Trop(f')\cap N_\R'$ with the complex whose cells are the faces of 
  $P'$, but this extra structure does not seem very useful.
\end{remsub}

\paragraph[The Newton complex] We \label{par:newton.complex}
use the notation in~\parref{par:trop.poly.cplx.1}.
Let $\pi: M_\R\times\R\to M_\R$ denote the projection onto the first
factor.  For $v\in N_\R'$ we define
\[ \check\gamma_v = \pi(\conv(\vertices_v(f))). \]
This is an integral $\Z$-affine 
polytope in $\spn(\tau)\subset M_\R$.  Again 
it is standard~\cite[Corollary~2.1.2]{ekl:non_arch_amoebas}
that $\New(f,\tau) \coloneq \{\check\gamma_v:v\in N_\R'\}$ is an
(integral $\Z$-affine) polytopal complex in $\spn(\tau)$, called the
\emph{Newton complex} of $f$.  When $\tau=\{0\}$ we omit it and simply
write $\New(f)$.  It is clear that the support of $\New(f,\tau)$ is
\[ |\New(f,\tau)| = \conv\{u\in S_\sigma\cap\tau^\perp~:~a_u\neq 0\}; \]
this is the \emph{Newton polytope} of $f$.  

\begin{egsub} Continuing  \label{eg:tropical.line.2}
  with~\parref{eg:tropical.line.1}, we have
  \[\begin{split}
    \check\gamma_{\{(1,1)\}} &= \conv\{(0,0),(-1,0),(0,-1)\} \\
    v\in R_1\implies \check\gamma_v &= \conv\{(-1,0),(0,-1)\} \\
    v\in R_2\implies \check\gamma_v &= \conv\{(0,0),(-1,0)\} \\
    v\in R_3\implies \check\gamma_v &= \conv\{(0,0),(0,-1)\}.
  \end{split}\]
  If $v$ is in one of the connected components of $\R^2\setminus\Trop(f)$
  then $\check\gamma_v$ is one of the vertices
  $\{(0,0)\}$, $\{(-1,0)\}$, $\{(0,-1)\}$, so $\inn_v(f)$ is a monomial.
  See Figure~\ref{fig:trop_new}.

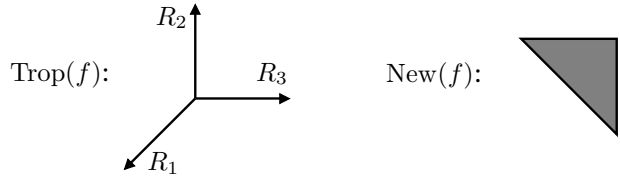
  \genericfig[ht]{trop_new}{The tropicalization and Newton complex of
    $f=x+y+\lambda$.} 

\end{egsub}

The complexes $\Trop(f)\cap N_\R'$ and $\New(f,\tau)$ are dual to each
other in the following sense:

\begin{propsub}
  We use the notation of~\parref{par:newton.complex}.
  \begin{enum}
  \item For $v,v'\in\Trop(f)\cap N_\R'$ we have $\gamma_v\prec\gamma_{v'}$ if
    and only if $\check\gamma_v\succ\check\gamma_{v'}$.  In particular, 
    $\gamma_v=\gamma_{v'}\iff\check\gamma_v=\check\gamma_{v'}$.

  \item For $v\in\Trop(f)\cap N_\R'$ the cells $\gamma_v$ and
    $\gamma_{v'}$ are orthogonal to each other in the sense that
    the linear subspace of $N_\R'$ associated to the affine span of
    $\gamma_v$ is orthogonal to the linear subspace of $\spn(\tau)$
    associated to the affine span of $\check\gamma_v$.

  \item For $v\in\Trop(f)\cap N_\R'$ we have
    $\dim(\gamma_w)+\dim(\check\gamma_w)=\dim_\R(N_\R')$. 
  \end{enum}
\end{propsub}

For $v\in\Trop(f)\cap N_\R'$ we call $\check\gamma_v$ the \emph{dual cell}
to $\gamma_v$.  This establishes a bijection between the cells of 
$\Trop(f)\cap N_\R'$ and the positive-dimensional cells of $\New(f,\tau)$
(the zero-dimensional cells correspond to the connected components of 
$N_\R'\setminus\Trop(f)$).  The ``duality'' between $\Trop(f)$ and
$\New(f,\tau)$ is not intrinsic (indeed, $\New(f)$ contains multiplicity
information missing from $\Trop(f)$); rather, they are related
manifestations of the combinatorial structure of the power series $f$
living in dual vector spaces.

\subparagraph
We resume the notation of~\parref{par:trop.poly.cplx.2}, so
$f\in K\angles{U_P}$ and $f'\in K[S_\sigma\cap\tau^\perp]$.
For $v\in P'$ we have $\vertices_v(f) = \vertices_v(f')$, so
\[ \check\gamma_v = \pi(\conv(\vertices_v(f))) = \pi(\conv(\vertices_v(f'))). \]
We define 
\[ \New(f,\tau) \coloneq 
\{\check\gamma_v~:~v\in P'\}\subset\New(f',\tau). \]
This is not in general a polyhedral complex as it may well happen that 
there exist $v\in P'$ and $v'\in N_\R'$ such that 
$\check\gamma_{v'}\prec\check\gamma_v$ but $\check\gamma_{v'}$ is not a
cell of $\New(f,\tau)$ (i.e.\ the corresponding cell $\gamma_{v'}$ is not
contained in $P'$).  We will only use the fact that there is a polytope
$\check\gamma_v = \pi(\conv(\vertices_v(f)))\in\New(f,\tau)$ associated to
every $v\in\Trop(f)\cap N_\R'$.

\begin{remsub} If \label{rem:new.trop.agree}
  $|f-g|_P\ll 1$ then $\vertices_{P'}(f) = \vertices_{P'}(g)$
  by~(\ref{lem:inu.finiteness},ii) and therefore  
  $\Trop(f) = \Trop(g)$ and $\New(f,\tau) = \New(g,\tau)$.
\end{remsub}

\begin{rem}
  Let $\sigma\subset N_\R$ be a pointed cone and let 
  $U\subset X(\sigma)^\an$
  be an admissible open subset that can be written as a union
  of polyhedral subdomains $\{U_{P_i}\}$ associated to polyhedra $P_i$
  with cone of unbounded directions 
  $\sigma$.  For instance we can take $U$ to be the rigid-analytic open unit ball 
  $\D^n_K = \bigcup_{r>0}\trop\inv([r,\infty]^n)$ inside of
  $\A^{n,\an}_K$, or we can take $U$ to be the analytic torus
  $T^\an = X(\{0\})^\an = \bigcup_{r>0}\trop\inv([-r,r]^n)$.  There is an
  evident tropicalization map $\trop: |U|\to\bigcup P_i$.  Let
  $f$ be an analytic function on $U$ and define $\Trop(f)$ to be the
  closure of $\trop(|V(f)|)$.  The finiteness
  lemma~\parref{lem:inu.finiteness} implies that $\Trop(f)$ is a ``locally
  finite polyhedral complex''.  This complex is not in general finite but is
  still interesting to study.
\end{rem}

\begin{rem}
  Let $P\subset N_\R$ be an integral $\Gamma$-affine pointed polyhedron and let 
  $Y\subset U_P$ be the closed analytic subspace defined by some ideal
  $\fa\subset K\angles{U_P}$.  If $P$ is a polytope then
  Gubler~\cite[Proposition~5.2]{gubler:tropical} has shown using the theory
  of semistable alterations of rigid spaces that $\Trop(Y)$
  is a finite union of (non-canonical) integral $\Gamma$-affine polytopes (among
  other things), as is the case for subschemes of a torus.  Such a result
  would follow from~\parref{par:trop.poly.cplx.2} for a pointed polyhedron $P$ 
  if one knew that 
  $\Trop(Y) = \bigcap_{i=1}^r \Trop(f_i)$ for some finite list of elements
  $f_1,\ldots,f_r\in\fa$.  While $\fa$ is certainly finitely generated it is
  \emph{not} necessarily the case (even for Laurent polynomials) that
  the intersection of the tropicalizations of a set of generators is equal
  to $\Trop(Y)$ (see~\ref{eg:non.proper.int}).  What one needs is a
  theorem that there exists a   ``universal Gr\"obner basis''
  of the ideal $\fa$ in $K\angles{U_P}$; see for
  instance~\cite[\S2]{speyer_sturmfels:tropical_grassmannian}.  The author
  would guess that such a theorem, suitably formulated, would be true.
  This issue is certainly deserving of further study as such a theorem
  would form an important part of the foundations of a theory of tropical
  analytic geometry.  

\end{rem}

\section{Continuity of roots I: the global version}
\label{sec:controots1}

\paragraph In this section we give a tropical
criterion~\parref{thm:controots1} for a family 
of $n$-tuples of power series in $n$ variables parametrized by a
one-dimensional base $S$
to define a rigid space that is finite and flat over $S$, so that
the number of common zeros of any member of the family is
independent of the parameter.  This will be a key ingredient
in~\S\ref{sec:non.transverse}.  A weaker version of this result has
appeared in~\cite{jdr:thesis}, 
where it was useful in explicitly counting the number of zeros of a
complicated system of power series by deforming the problem to a
much simpler one. 

The main rigid-analytic ingredient used in this section
is the direct image theorem for rigid spaces.  The statement is
exactly the same as the direct image theorem for algebraic geometry; the
subtlety is in the \emph{definition} of properness for morphisms of rigid spaces,
which we review below.

\paragraph The \label{par:controots1.intuition}
intuitive idea behind the continuity of roots theorem is as
follows.  Suppose for this paragraph that $K=\bar K$ (for simplicity).
Let $f_1,\ldots,f_n\in K\angles{x_1,\ldots,x_n,t}$ and let
$Y\subset\B^n_K\times\B^1_K$ be the closed analytic subspace defined by
the ideal $(f_1,\ldots,f_n)$.  For $t_0\in|\B^1_K|$ let
$f_{i,t_0}$ be the image of $f_i$ in $K\angles{x_1,\ldots,x_n}$ and
let $Y_{t_0}$ be the fiber of $Y$ over $t_0$, so $Y_{t_0}$ is the space of common
zeros of $(f_{1,t_0},\ldots,f_{n,t_0})$.  Let $\rho\in |K^\times|$ with $\rho<1$
and suppose that $Y$ is in fact contained in the Weierstrass subdomain
$\B^n_K(\rho)\times\B^1_K$ of $\B^n_K\times\B^1_K$
(cf.~\parref{eg:laurent.domains}): that is, $Y$ is a 
closed subscheme of $\B^n_K(\rho)\times\B^1_K$ that is simultaneously a
closed subscheme of $\B^n_K\times\B^1_K$.  Tropically, if $P$ is the
polyhedron $\R^n_{\geq 0}$ then our condition is equivalent to 
$\Trop(Y_{t_0})$ being contained in the closure of a polyhedron 
$\R_{\geq r}^n$ for some $r > 0$ and all $t_0$.

Roughly,
points of $Y_{t_0}$ are ``trapped'' inside of the smaller ball $\B^n_K(\rho)$
since they cannot escape to the boundary of $\B^n_K$ --- that is, no points
of $Y_{t_0}$ can enter or leave $\B^n_K(\rho)$ as the parameter $t_0$ varies
since otherwise we would have points ``jumping over'' the annulus
$\B^n_K\setminus\B^n_K(\rho)$.  Hence all of the finite rigid spaces
(equivalently, finite schemes) $Y_{t_0}$ must have the same length.

\paragraph To say that a ball is contained in the ``interior'' of a
larger ball is basically the notion of relative compactness:

\begin{defn*}[{\cite[\S9.6.2]{bgr:nonarch}}]
  Let $X = \Sp(A)$ and $Y = \Sp(B)$ be $K$-affinoid spaces, let 
  $f: X\to Y$ be a morphism, and let $U\subset X$ be an affinoid
  subdomain.  We say that \emph{$U$ is relatively compact in $X$ over $Y$}
  and we write $U\Subset_Y X$ provided that we can find a closed immersion
  $X\inject \B^n_K\times Y$ over $Y$ such that 
  $U\subset\B^n_K(\rho)\times Y$ for some $\rho\in|\bar K{}^\times|$ with
  $\rho<1$. 
\end{defn*}

In the above example~\parref{par:controots1.intuition}, we have
$\B^n_K(\rho)\times\B^1_K\Subset_{\B^1_K}\B^n_K\times\B^1_K$.

\paragraph Kiehl's notion of  properness \label{par:properness}
for morphisms of rigid spaces is defined in terms of relative compactness.  

\begin{defn*}[{\cite[Definition~9.6.2/2]{bgr:nonarch}}]
  Let $f:X\to Y$ be a morphism of rigid spaces, and for simplicity assume
  that $Y$ is affinoid.  We say that $f$ is
  \emph{proper} if it is separated (i.e.\ the diagonal is closed) and if it
  satisfies the following condition: there exist two admissible affinoid
  coverings $\{U_i\}_{i=1}^n$ and $\{V_i\}_{i=1}^n$ of $X$ such that 
  $U_i\Subset_Y V_i$ for all $i=1,\ldots,n$.
\end{defn*}

Properness over a general base is defined in such a way as to be
local on the base.  

\begin{thmsub}[{\cite[Theorem~9.6.3/1]{bgr:nonarch}}]
  Let $f:X\to Y$ be a proper morphism of rigid spaces and let $\sF$ be a
  coherent sheaf of $\sO_X$-modules.  Then $f_*\sF$ is a coherent
  sheaf of $\sO_Y$-modules. 
\end{thmsub}

The definition of a coherent sheaf of modules on a rigid space is similar
to the analogous definition in algebraic geometry, but the precise
definition is not important for our purposes since we will only use the
following simple consequence: 

\begin{corsub} Let \label{cor:proper.affinoid.finite}
  $X = \Sp(A)$ and $Y = \Sp(B)$ be affinoid spaces and let
  $f:X\to Y$ be a proper morphism.  Then $B$ is finite as an $A$-module. 
\end{corsub}

\begin{egsub} In \label{eg:finite.example}
  the situation of~\parref{par:controots1.intuition}, the spaces $Y$
  and $\B^1_K$ are affinoid and the map $Y\to\B^1_K$ is proper: in fact
  $Y\Subset_{\B^1_K} Y$ since 
  \[ Y\subset\B^n_K(\rho)\times\B^1_K\Subset_{\B^1_K}\B^n_K\times\B^1_K. \]
  Therefore $Y\to\B^1_K$ is finite; flatness follows from an
  unmixedness argument as in the proof of~\parref{thm:controots1}.
\end{egsub}

The following generalizes the fact~\parref{par:controots1.intuition}
that $\B^n_K(\rho)\times\B^1_K\Subset_{\B^1_K}\B^n_K\times\B^1_K$.

\begin{lem} Let \label{lem:relatively.compact}
  $P'\subset P\subset N_\R$ be integral $\Gamma$-affine pointed polyhedra such that
  $\tau = \cU(P')$ is a face of $\sigma = \cU(P)$ (so
  $\bar P{}'\subset\bar P\subset N_\R(\sigma)$).  If $P'$ is contained in the
  (topological) interior of $P$ then $U_{P'}\Subset_K U_P$.
\end{lem}

\pf First note that $U_{P'}$ is an admissible affinoid open subset of 
$X(\tau)^\an\subset X(\sigma)^\an$ and is 
therefore an affinoid subdomain of $U_P$~\cite[Proposition~9.3.1/3]{bgr:nonarch}.
Choose generators $u_1,\ldots,u_r$ for $S_\sigma$ such that we can write 
$P = \bigcap_{i=1}^r\{v\in N_\R:\angles{u_i,v}\leq b_i\}$, so
$U_P$ is a closed subspace of $\prod_{i=1}^r \B^1_K(\rho_i)$ as
in~\parref{rem:U.P.immersion} where $\rho_i = \exp(b_i)$.
Since $P'$ is contained in the interior of $P$, we can find 
$c_i\in\Gamma$ with $c_i < b_i$ such that 
$P'\subset\bigcap_{i=1}^r\{v\in N_\R:\angles{u_i,v}\leq c_i\}$.  
Setting $\mu_i = \exp(c_i)$, we have 
$U_{P'}\subset\prod_{i=1}^r\B^1_K(\mu_i)\Subset_K\prod_{i=1}^r\B^1_K(\rho_i)$.\qed

Generalizing~\parref{eg:finite.example} we have the following consequence
of~\parref{lem:relatively.compact}, which is a tropical criterion for an
affinoid space to be finite.

\begin{prop} Let \label{prop:automatically.finite}
  $P\subset N_\R$ be an integral $\Gamma$-affine pointed polyhedron and let 
  $Y$ be a closed analytic subspace of $U_P$ such that $\Trop(Y)$ is
  contained in the interior $\bar P{}^\circ$ of $\bar P$.  Then $Y$ is finite.
\end{prop}

\pf Let $\fa\subset K\angles{U_P}$ be the ideal defining $Y$ and let
$A = K\angles{U_P}/\fa$, so $Y = \Sp(A)$.
Write $P = \bigcap_{i=1}^r\{v\in N_\R:\angles{u_i,v}\leq b_i\}$ as in
the proof of~\parref{lem:relatively.compact}.  
By the maximum modulus principle as applied to the image of $x^{u_i}$ in
$A$, there exists $c_i\in\Gamma$ with $c_i < b_i$ such that
$\Trop(Y)$ is contained in the polyhedron
$P' = \bigcap_{i=1}^r\{v\in N_\R:\angles{u_i,v}\leq c_i\}$.  Since
$Y\subset U_{P'}\Subset_K U_P$ this shows that $Y\Subset_K Y$, so 
$Y\to\Sp(K)$ is proper and hence $Y$ is
finite by~\parref{cor:proper.affinoid.finite}.\qed

Note that~\parref{prop:automatically.finite} applies equally well to a
closed subscheme $Y$ of $X(\sigma)$ such that $\Trop(Y)$ is contained in
the interior of a compactified polyhedron $\bar P$.  This is also a
consequence of the balancing condition for tropical varieties: a
positive-dimensional tropical variety is ``infinite in all directions''.

\paragraph The \label{par:controots1.setup}
following is the setup for the continuity of roots theorem.
Let $\Delta$ be an integral pointed fan in $N_\R$ and let
$P_1,\ldots,P_r\subset N_\R$ be integral $\Gamma$-affine pointed polyhedra such
that $\sigma_i = \cU(P_i)\in\Delta$.  Let $\cP = \{P_1,\ldots,P_r\}$ and
define
\[ U_\cP = U_{P_1}\cup\cdots\cup U_{P_r} 
= \trop\inv(\bar P_1\cup\cdots\cup\bar P_r)
\subset X(\Delta)^\an. \]
This is an admissible open subset of $X(\Delta)^\an$.
By a \emph{family of subspaces of $U_\cP$} parametrized by a rigid space
$S$ we mean a closed analytic subspace $Y$ of $U_\cP\times S$.  Letting
$\pi:Y\to S$ be the projection onto the second factor, for
$s\in S$ we set $Y_s = \pi\inv(s)$; this is a closed analytic subspace of
$\kappa(s)\hat\tensor_K U_\cP$.  We say that $Y$ is a 
\emph{relative complete intersection} if for all $s\in S$ there is an
affinoid neighborhood $\Sp(A)\subset S$ such that for all $i=1,\ldots,r$
the closed subspace $Y\cap(U_{P_i}\times\Sp(A))$ of $U_{P_i}\times\Sp(A)$
is defined by $n$ 
equations $f_1,\ldots,f_n\in A\angles{U_P}$, where 
$A\angles{U_P} = A\hat\tensor_K K\angles{U_P}$ as 
in~\parref{defn:polyhedral.subdomain}

\begin{thm}[Continuity of roots I] We fix: \label{thm:controots1}
  \begin{enum}
  \item $S$ a normal connected rigid space of dimension one.
  \item $\Delta$ an integral pointed fan in $N_\R$.
  \item $\cP = \{P_1,\ldots,P_r\}$ a collection of integral $\Gamma$-affine
    polyhedra in $N_\R$ such that $\cU(P_i)\in\Delta$ for all $i$.
  \item $\cP' = \{P_1',\ldots,P_r'\}$ a second collection of integral
    $\Gamma$-affine polyhedra in $N_\R$ such that
    $\cU(P_i') = \cU(P_i)$ and $P_i'\subset P_i^\circ$ for all $i$.
  \item $Y$ a family of subspaces of $U_\cP$ parametrized by $S$ which is
    a relative complete intersection.
  \end{enum}
  Suppose that $\Trop(Y_s)$ is contained in 
  $\bigcup_{i=1}^r\bar P{}_i'$ for all $s\in|S|$.  Then (a) if each fiber
  $Y_s$ is a finite set then $\pi: Y\to S$ is finite and flat, and (b)
  every \emph{finite} fiber $Y_s$ has the same length even if $\pi$ has
  non-finite fibers.  
\end{thm}

Before giving the proof we mention the following (easier to formulate)
special case, which follows from~\parref{thm:controots1}
combined with~\parref{prop:automatically.finite}. 

\begin{cor} Let \label{cor:controots1}
  $A$ be an affinoid algebra that is a Dedekind domain, let 
  $P\subset N_\R$ be an integral $\Gamma$-affine pointed polyhedron, let
  $f_1,\ldots,f_n\in A\angles{U_P}$, and let
  $Y\subset U_P\times\Sp(A)$ be the subspace cut out by $f_1,\ldots,f_n$.
  If $\Trop(Y_s)$ is contained in the interior of $\bar P$ for all 
  $s\in|\Sp(A)|$ then $Y\to\Sp(A)$ is a finite, flat map.
\end{cor}

\begin{eg} Let \label{eg:trop.deform}
  $N = M = \Z^2$ and let $x = x^{(-1,0)}$ and $y = x^{(0,-1)}$ as
  in~\parref{eg:tropical.line.1}.  Let $\tau = \R_{\geq 0}(1,0)$.
  Let $\lambda\in K^\times$ have 
  valuation $1$, let $\mu\in K^\times$ have valuation $2$, and define 
  \[ f_1(x,y,t) = x + t y + \lambda,\quad
  f_2(x,y,t) = \mu x + y + \lambda\quad \in \quad K[S_\tau][t] 
  = K[x^{\pm 1},y,t]. \]
  The tropicalizations of the specializations $f_{1,t_0}$ and $f_{2,t_0}$ in
  $N_\R(\tau)$ for a specific 
  value of $t_0$ are drawn in Figure~\ref{fig:deform}.  Let
  \[ P = \{ (u_1,u_2)~:~ u_1\geq 0,~u_2\in[0,2]\}, \]
  so $\cU(P) = \tau$ and $\bar P\subset N_\R(\tau)$.  Let 
  $Y\subset X(\tau)\times\A^1_K$ be the subscheme defined by $(f_1,f_2)$,
  and for $\epsilon\in\Gamma$ positive let
  $Y_\epsilon = Y^\an\cap(U_P\times S_\epsilon)$, where
  $S_\epsilon$ is the annulus of inner radius $\exp(-\epsilon)$ and outer
  radius $\exp(\epsilon)$ as in~\parref{eg:trop.annulus}.  It is clear
  from the picture that the family $Y_\epsilon$ satisfies the hypotheses
  of~\parref{cor:controots1} for small enough $\epsilon$ since
  $\Trop(f_{1,t_0})\cap\Trop(f_{2,t_0})\subset\bar P{}^\circ$ for
  $\val(t_0)$ near zero.  Therefore $Y_\epsilon\to S_\epsilon$ is finite and
  flat, so in particular every fiber has the same length (of $1$),
  \emph{including} the fiber over $t_0=1$ where the intersection is a
  completed ray.  

  This picture underlies much of~\S\ref{sec:non.transverse}, where it is
  essential that we work with families of translations parametrized by
  rigid-analytic annuli $S_\epsilon$.

  \genericfig[ht]{deform}{Pictures of $\Trop(f_{1,t_0})$ and $\Trop(f_{2,t_0})$ in
    $N_\R(\tau)$ from~\parref{eg:trop.deform} evaluated at a generic value of 
    $t_0\in K^\times$ with $0<\val(t_0)\ll 1$.  The dotted vertical line on the
    right is $N_\R/\spn(\tau)$ and the solid line segment is $\bar P$; the
    dots in $\bar P$ are included in the tropicalizations.}

\end{eg}

\pf[of~\parref{thm:controots1}]
By hypothesis $Y\subset U_{\cP'}\times S$, so since 
$U_{P_i'}\times S\Subset_S U_{P_i}\times S$ by~\parref{lem:relatively.compact}
we have $Y\cap(U_{P_i'}\times S)\Subset_S Y\cap(U_{P_i}\times S)$ for all
$i$ and hence $\pi: Y\to S$ is proper.  Suppose that 
$Y_s$ is a finite set for all $s\in|S|$.  Then $Y$ is finite over $S$
by~\cite[Corollary~9.6.3/6]{bgr:nonarch}, so it
suffices to prove that $Y$ is $S$-flat.  The assertion is local on $Y$ and
$S$, so we may assume that $S = \Sp(A)$ is affinoid, where $A$ is an
affinoid algebra that is a Dedekind domain and therefore
Cohen-Macaulay.  Let $f_1,\ldots,f_n\in A\angles{U_{P_i}}$ be a collection
of power series defining $Y\cap(U_{P_i}\times\Sp(A))$, let 
$\fa = (f_1,\ldots,f_n)$, and let $B = A\angles{U_{P_i}}/\fa$, so
$Y\cap(U_{P_i}\times\Sp(A)) = \Sp(B)$.  Since $A\angles{U_{P_i}}$ is a
flat $A$-algebra~\cite[Theorem~A.1.5]{conrad:relampleness} with
Cohen-Macaulay fiber rings over maximal
ideals~(\ref{prop:polyhedral.affinoids},v), it follows 
from~\cite[Theorem~23.9]{matsumura:crt} that $A\angles{U_{P_i}}$ is itself
Cohen-Macaulay.  Thus $A\angles{U_{P_i}}$ is catenary of dimension $n+1$,
so by Krull's principal ideal theorem, if $\fp$ is a minimal prime of $B$
then $\dim(B/\fp)\geq 1$.  But the fibers of $\pi$ have dimension zero,
so $\dim(B/\fp)=1$ and hence by the unmixedness
theorem~\cite[Theorem~17.6]{matsumura:crt}, $\fa$ has no embedded prime
ideals.  Thus every associated prime of $B$ contracts to the zero ideal of
$A$, so since $A$ is a Dedekind domain, $B$ is a flat $A$-module.

In the general case, the theorem on semicontinuity of fiber
dimension~\cite[Theorem~4.9]{ducros:relative_dimension} implies that the set 
\[ Z = \{\eta\in|Y| ~:~ \dim_\eta(Y_{\pi(\eta)})\geq 1 \} \]
is (Zariski-)closed in $Y$, so the proper mapping
theorem~\cite[Proposition~9.6.3/3]{bgr:nonarch} implies that 
$\pi(Z)$ is a closed subset of $S$, which has dimension zero
if $\pi$ has any finite fibers.  Deleting $\pi(Z)$ from $S$ does
not affect its connectedness, so we are reduced to the case treated above.\qed

\begin{rem} It \label{rem:weaken.hypotheses}
  may be possible to weaken the hypotheses of~\parref{thm:controots1}
  to only require that $\Trop(Y_s)$ be contained in 
  $\bigcup_{i=1}^r\bar P{}_i^\circ$ for each $s$, or even in the interior of
  $\bigcup_{i=1}^r\bar P_i$, but it is not immediately obvious how one
  would do so.
\end{rem}

\section{Continuity of roots II: the local version}
\label{sec:controots2}

\paragraph The purpose of this section is to show that if
$f_1,\ldots,f_n$ is any family of $n$-tuples of power series in $n$
variables parametrized by a one-dimensional rigid space $S$, and if
$t\in|S|$ is a point
such that the specializations $f_{1,t},\ldots,f_{n,t}$ at $t$ have finitely
many common zeros, then $f_1\ldots,f_n$ defines a finite and flat rigid
space over a small affinoid neighborhood of $t$ in $S$.  This is the 
rigid-analytic fact that allows us to use a polynomial approximation
argument in order to derive the local intersection multiplicity formula
for rigid spaces from the analogous theorem for schemes.  The proof
of~\parref{thm:controots2} is more technical than~\parref{thm:controots1},
and we will assume more familiarity with rigid analytic spaces in it.  In
particular, we will assume that the reader has some 
knowledge of Raynaud's theory of formal models, which we briefly review
in~\parref{par:formal.models}. 

We begin with the statement of the theorem we will prove:

\begin{thm}[Continuity of roots II] Let \label{thm:controots2}
  $A$ be a $K$-affinoid algebra that is a Dedekind domain and let 
  $S = \Sp(A)$.  Let $X = \Sp(B)$ be a Cohen-Macaulay affinoid space of
  dimension $n+1$, let $f_1,\ldots,f_n\in B$, and let $Y\subset X$ be the
  subspace defined by the ideal $\fa = (f_1,\ldots,f_n)$.  Suppose that we
  are given a morphism $\pi: X\to S$ and a point $t\in|S|$ such that 
  the fiber $Y_t = \pi\inv(t)\cap Y$ has dimension zero.  Then there is an
  affinoid subdomain $U\subset S$ containing $t$ such that 
  $\pi\inv(U)\to U$ is finite and flat.  
\end{thm}

In particular, the rigid space $Y_s = \pi\inv(s)\cap Y$ is finite for all
$s\in|U|$ and has the same length as $Y_t$.

\begin{eg}
  The following special case makes~\parref{thm:controots2} look very much
  like a theorem of continuity of roots.
  Let $X = \B^n_K\times\B^1_K$ and $S = \B^1_K$, with $\pi:X\to S$ the projection
  onto the second factor.  Let 
  $f_1,\ldots,f_n\in K\angles{x_1,\ldots,x_n,t}$.  If the
  specializations $f_{1,0},\ldots,f_{n,0}$ at $0$ have only finitely many zeros
  in $\B^n_K$ then there exists $\epsilon>0$ such that 
  $f_{1,s},\ldots,f_{n,s}$ have the same number of zeros (counted with
  multiplicity) in $\B^n_{\kappa(s)}$ as $f_{1,0},\ldots,f_{n,0}$ when
  $|s| < \epsilon$.
\end{eg}

\paragraph Here \label{par:formal.models} we recall some notions used in
Raynaud's theory of formal 
models.  The main reference is Bosch, L\"utkebohmert, and Raynaud's series
of papers~\cite{bl:fmI,bl:fmII,blr:fmIII,blr:fmIV}.  The
\emph{ring of restricted power series in $n$ variables} over $\sO_K$ is
\[ \sO_K\angles{x_1,\ldots,x_n} = 
\bigg\{\sum_{\nu}a_\nu x^\nu\in \sO_K\ps{x_1,\ldots,x_n}~:~
|a_\nu|\to 0\bigg\} 
= \{ f\in K\angles{x_1,\ldots,x_n}~:~ |f|_{\sup}\leq 1 \}. \]
An $\sO_K$ algebra $A$ admitting a surjective homomorphism
$\phi: \sO_K\angles{x_1,\ldots,x_n}\surject A$ for some $n$ is called
\emph{topologically of finite type} or \emph{tf type}; if we can choose
$\phi$ such that $\ker(\phi)$ is a finitely generated ideal, we say that
$A$ is \emph{topologically of finite presentation} or 
\emph{tf presentation}.  If $A$ is tf  type and is
$\sO_K$-flat we say that $A$ is an \emph{admissible $\sO_K$-algebra}; in
this case $A$ is automatically tf
presentation and is complete and separated in the $\varpi$-adic 
topology for any nonzero $\varpi\in\fm_K$~\cite[Proposition~1.1]{bl:fmI}.
Note that $A$ is $\sO_K$-flat if and only if it has no $\varpi$-torsion.
An \emph{admissible formal $\sO_K$-scheme} is a formal
$\Spf(\sO_K)$-scheme that is locally isomorphic to the formal spectrum of
an admissible $\sO_K$-algebra (equipped with the $\varpi$-adic
topology).  

There is a \emph{rigid generic fiber} functor $\fX\mapsto \fX_\rig$ from
the category of 
quasi-compact admissible formal $\sO_K$-schemes to the category of
quasi-compact and quasi-separated rigid spaces over $K$; it becomes an 
equivalence after inverting so-called admissible formal blow-ups in the
source category.  If $\fX = \Spf(\sO_K\angles{x_1,\ldots,x_n}/\fa)$ is the
formal spectrum of an admissible $\sO_K$-algebra then 
$\fX_\rig = \Sp(K\angles{x_1,\ldots,x_n}/K\fa)$.  
The rigid generic fiber functor satisfies many
compatibility properties 
including respecting open immersions and fiber products.  If 
$X$ is a rigid space, an admissible formal scheme $\fX$ such that
$\fX_\rig\cong X$ is called a \emph{formal model} for $X$; 
such a model always exists when $X$ is quasi-compact and
quasi-separated.  If $X = \Sp(K\angles{x_1,\ldots,x_n}/\fa)$ is an
affinoid space then 
$\Spf(\sO_K\angles{x_1,\ldots,x_n}/\fa\cap\sO_K\angles{x_1,\ldots,x_n})$ is a
formal model for $X$; however, most formal models for $X$ will not be affine.  

Let $f:X\to Y$ be a morphism of quasi-compact and quasi-separated rigid
spaces.  The power of Raynaud's theory lies in the ability to choose
formal models $\fX$ and $\fY$ for $X$ and $Y$, respectively, along with
a morphism $\phi:\fX\to\fY$ such that $\phi_\rig = f$, in such a way that
$\phi$ retains any ``nice'' properties of $f$ (e.g.\ flatness).  This
allows one to use algebraic geometry to prove statements about rigid
spaces.

\begin{notn}
  We fix a nonzero element $\varpi\in\fm_K$.  For $m\geq 0$ we let
  $\sO_{K,m} = \sO_K/\varpi^{m+1}\sO_K$, and if $\fX$ is a formal
  $\Spf(\sO_K)$-scheme we let 
  $\fX_m = \sO_{K,m}\tensor_{\sO_K}\fX$.
\end{notn}

If $\fX$ is an admissible formal $\sO_K$-scheme then each $\fX_m$ is a flat
$\sO_{K,m}$-scheme of finite type (having the same underlying topological
space as $\fX$).  The following converse statement is well-known;
see~\cite[\S1]{bl:fmI}. 

\begin{lem} Let \label{lem:fibral.admissibility}
  $\{A_m\}_{m\geq 0}$ be an inverse system of $\sO_K$-algebras
  such that for all $m\geq 0$ the map $A_{m+1}\to A_m$ identifies $A_m$
  with $\sO_{K,m}\tensor_{\sO_K} A_{m+1}$.  If $A_0$ is an $\sO_{K,0}$-algebra of
  finite type and $A_m$ is a flat $\sO_{K,m}$-algebra for every $m\geq 0$ then
  $A = \varprojlim_m A_m$ is an admissible $\sO_K$-algebra and the
  natural maps $\sO_{K,m}\tensor_{\sO_K} A\to A_m$ are isomorphisms.
\end{lem}

\paragraph Let \label{par:reduction.map}
$\fX$ be an admissible formal $\sO_K$-scheme.  There is a functorial
\emph{reduction map} $\red:|\fX_\rig|\to|\fX|$, defined as follows.
Let $\xi\in|\fX_\rig|$ and let $\fU = \Spf(A)\subset\fX$ be a formal affine
such that $\xi$ is a point of
$\fU_\rig = \Sp(K\tensor_{\sO_K} A)$.  Then $\xi$ corresponds
to a surjective homomorphism $\phi: K\tensor_{\sO_K} A\surject K'$, where
$K'=\kappa(\xi)$ is a finite extension of $K$.  For boundedness reasons we
have $\phi(A)\subset\sO_{K'}$; the point $\red(\xi)$ corresponds to the
contraction of $\fm_{K'}$ in $A$.

In the above situation the ring $R = \phi(A)\subset\sO_{K'}$ is
a finite admissible local $\sO_K$-algebra of dimension one and the closed
immersion $\Spf(R)\inject\fX$ is called a \emph{rig-point} of $\fX$;
see~\cite[\S3.1]{bl:fmI}.  In this way the rig-points of $\fX$ correspond
naturally to the points of $\fX_\rig$~\cite[Lemma~3.4]{bl:fmI}.

\begin{lem} Let \label{lem:etale.invariance}
  $\fX$ be a quasi-compact admissible formal $\sO_K$-scheme and let
  $g: \sY\to\fX_0$ be an \'etale morphism of finite-type
  $\sO_{K,0}$-schemes.  There is 
  a unique (up to unique isomorphism) admissible formal $\sO_K$-scheme
  $\fY$ equipped with a morphism $f: \fY\to\fX$ such that 
  $\fY_0 \cong\sY$ with $f_0: \fY_0\to \fX$ identified with $g$, and such
  that $f_m: \fY_m\to\fX_m$ is \'etale for all $m\geq 0$.  Moreover, $f$
  is flat and $f_\rig:\fY_\rig\to\fX_\rig$ is an \'etale morphism of rigid
  spaces. 
\end{lem}

\pf The existence and uniqueness of $f:\fY\to\fX$ is a consequence of the
infinitesimal invariance of the \'etale site~\cite[18.1.2]{egaIV_4},
along with~\parref{lem:fibral.admissibility}.  The flatness of $f$ follows
from the fibral flatness criterion over general valuation
rings~\cite[Lemma~1.6]{bl:fmI}.  That $f$ is \'etale is a special case
of~\cite[Corollary~3.10]{blr:fmIII} --- or one can prove it directly by
reducing to the case of standard \'etale morphisms and using the Jacobi
criterion.\qed

The following proposition is a translation of the structure theorem for
separated finite-type morphisms with a finite
fiber~\cite[18.12.3]{egaIV_4} to rigid spaces, using formal models.

\begin{prop} Let \label{prop:str.thm.finite.fiber}
  $f:X\to Y$ be a separated morphism of quasi-compact quasi-separated
  rigid spaces, and suppose that the fibers of $f$ are finite.  For any
  $\eta\in|Y|$ there is an \'etale morphism $g:Y'\to Y$ and a point
  $\eta'\in g\inv(y)$ such that the product
  $X' = X\times_Y Y'$ decomposes into a disjoint union  of rigid spaces
  $X' = X_1'\djunion X_2'$ in such a way that
  $X_1'\to Y'$ is finite and $X_2'\to Y'$ has empty $\eta'$-fiber.
\end{prop}

\pf By~\cite[Corollary~5.10(b)]{bl:fmII} there exist formal models
$\fX$ and $\fY$ for $X$ and $Y$, respectively, along with a morphism
$\phi:\fX\to\fY$ with dimension-zero fibers such that $\phi_\rig = f$.
The morphism $\phi_0:\fX_0\to\fY_0$ is separated by
\cite[Proposition~4.7]{bl:fmI}, so by 
\cite[18.12.3]{egaIV_4} there is an \'etale morphism 
$\psi_0:\sY'\to\fY_0$ and a point $\eta_0'\in|\sY'|$ over 
$\eta_0 = \red(y)$ such that $\sX' = \sY'\times_{\fY_0}\fX_0$ breaks up into a
disjoint union $\sX_1'\djunion\sX_2'$ with $\sX_1'$ 
finite over $\sY'$ and $\sX_2'\to\sY'$ having empty $\eta_0'$-fiber.  Let
$\psi:\fY'\to\fY$ be the unique lift of 
$\psi_0:\sY'\to\fY_0$ as in~\parref{lem:etale.invariance}, and let 
$\fX' = \fY'\times_{\fY}\fX$ be the fiber product (in the
category of formal $\Spf(\sO_K)$-schemes; $\fX'$ is then admissible because 
$\fY'\to\fY$ is flat), so  
$\fX'_0 = \sX'$.  Since the topological spaces underlying
$\fX'$ and $\sX'$ are the same, we have
$\fX' = \fX_1'\djunion\fX_2'$, where $\fX_i'$ is an admissible
formal $\sO_K$-scheme lifting $\sX_i'$ for $i=1,2$.  It follows from
\cite[Lemma~1.5]{bl:fmI} that $\fX_1'\to\fY'$ is finite, and
certainly $\fX_2'\to\fY'$ has empty $\eta_0'$-fiber.  
Let $X' = \fX'_\rig$, $Y' = \fY'_\rig$, 
$g = \psi_\rig:Y'\to Y$, and
$X_i' = (\fX_i')_\rig$ for $i=1,2$, so 
$X' = Y'\times_Y X = X_1'\djunion X_2'$.
Then $X_1'\to Y'$ is finite, and if $\eta'\in|Y'|$ is any point that
reduces to $\eta_0'$ then $X_2'\to Y'$ has empty $\eta'$-fiber.   

It remains to show that there exists $\eta'\in g\inv(\eta)$ reducing to
$\eta_0'$.  Let $\Spf(R)\inject\fY$ be the rig-point associated to $\eta$ as
in~\parref{par:reduction.map}, so $\Spf(R)_\rig = \{\eta\}$.
Consider the Cartesian squares
\[\xymatrix{
  {g\inv(\eta)} \ar[r] \ar[d] & {Y'} \ar[d]^g & & 
  {\fZ} \ar[r] \ar[d] & {\fY'} \ar[d]^\psi \\
  {\{\eta\}} \ar[r] & Y & & {\Spf(R)} \ar[r] & {\fY} 
}\]
where $\fZ$ is the fiber product $\fY'\times_\fY\Spf(R)$ in the category
of formal $\Spf(\sO_K)$-schemes; note that $\fZ$ is admissible because
$\psi$ is flat. 
Since the rigid generic fiber functor is compatible with fiber products,
the left square is canonically identified with the rigid generic fiber of
the right square.  The result now follows from the surjectivity of the
reduction map~\cite[Proposition~3.5]{bl:fmI}.\qed

\begin{lem} Let \label{lem:img.affinoid}
  $f:A\to B$ be a homomorphism of $K$-affinoid algebras and let
  $\phi:\Sp(B)\to\Sp(A)$ be the associated morphism of affinoid spaces.
  Let $\xi\in|\Sp(A)|$ be an element not contained in the image of $\phi$.
  Then there is an affinoid subdomain $U\subset\Sp(A)$ containing $\xi$ that
  is disjoint from the image of $\phi$.
\end{lem}

\pf Let $\fm\subset A$ be the maximal ideal corresponding to $\xi$, let 
$a_1,\ldots,a_r$ generate $\fm$, and let $a_i' = f(a_i)$.  
Since $f(\fm) B = B$, there exist $b_1,\ldots,b_r\in B$ such that 
$\sum_{i=1}^r a_i' b_i = 1$.  Let
$M > \max\{|b_1|_{\sup},\ldots,|b_r|_{\sup}\}$ with $M\in\Gamma$.  Then
for all $\eta\in|\Sp(B)|$ there is some $i$ such that
$|a_i'(\eta)| > 1/M$, so the Weierstrass subdomain
\[ U = \{ \xi'\in|\Sp(A)|~:~ |a_i(\xi')|\leq 1/M \text{ for all } i =
1,\ldots,r\} \]
satisfies our requirements.\qed

\pf[of~\parref{thm:controots2}]
By the theorem on semicontinuity of fiber
dimension for rigid spaces~\cite[Theorem~4.9]{ducros:relative_dimension}, the 
locus $Z$ of points $\eta\in|Y|$ not isolated in its fiber is a
Zariski-closed subset of the affinoid space $Y$, so $Z$ is the set
underlying an affinoid space.
By assumption $Z\cap Y_t=\emptyset$, so by~\parref{lem:img.affinoid},
after replacing $S$ with an affinoid 
subdomain containing $t$ we may assume that $Y\to S$ has finite fibers.
The flatness of $Y\to S$ now follows from the unmixedness theorem exactly
as in the proof of~\parref{thm:controots1}. 

By~\parref{prop:str.thm.finite.fiber} there is an
\'etale morphism $g: S'\to S$ and a point $t'\in |S'|$ in the fiber
over $t$ such that the fiber product $Y' = Y\times_{S}S'$
decomposes into a disjoint union 
$Y' = Y'_1\djunion Y'_2$, where $Y_1'\to S'$ is finite and 
$Y_2'$ has empty $t'$-fiber.  Replacing $S'$ with an affinoid neighborhood
of $t'$ disjoint from $g\inv(t)\setminus\{t'\}$ we may assume that
$S',Y',Y_1',Y_2'$ are all affinoid and 
that $g\inv(t)=\{t'\}$.  Then 
$Y_2'\to S$ has empty fiber over $t$, so again
by~\parref{lem:img.affinoid}, after replacing 
$S$ with an affinoid subdomain we may assume that
$Y_2' = \emptyset$, and therefore that $Y'\to S'$ is finite.
By \cite[Corollary~5.11]{bl:fmII} the image of 
$g$ is open, so we again shrink $S$ to assume that $g$ is surjective.
Then by descent theory for rigid spaces 
\cite[Theorems~4.2.7 and~4.2.2]{conrad:relampleness}, we have that
$Y\to S$ is finite.\qed

\section{Application: a local intersection multiplicity formula for rigid spaces}
\label{sec:newton.polygons}

\paragraph Osserman and Payne~\cite[\S5]{osserman_payne:lifting}
have proved a general theorem relating the multiplicities of
an intersection of subvarieties of a torus with the corresponding
multiplicities of the intersection of their tropicalizations.  In the case
of a dimension-zero complete intersection this theorem becomes a formula
for intersection numbers whose history begins with
Bernstein~\cite{bernstein:number_of_roots};
see~\parref{rem:bernstein.history}.  We use this multiplicity formula,
along with the continuity of roots theorem~\parref{thm:controots2} 
and a polynomial approximation argument, to derive an intersection multiplicity
formula~\parref{thm:multiplicity.ps} for
rigid spaces in the case of a complete intersection of dimension zero.
Theorem~\parref{thm:multiplicity.ps} is a natural generalization
of the theorem of the Newton polygon to a higher-dimensional setting;
see~\parref{eg:multiplicity.ps.dim1}.

Tropical intersection multiplicities are calculated in terms of the mixed
volume of a collection of polytopes (in the case of a dimension-zero
complete intersection): 

\begin{defn}
  The \emph{Minkowski Sum} of an $n$-tuple of polytopes
  $P_1,\ldots,P_n\subset N_\R$ is defined to be
  \[ P_1 + \cdots + P_n = \{ v_1 + \cdots + v_n~:~v_i\in P_i \}. \]
  For $\lambda\in\R_{\geq 0}$ we let 
  $\lambda P_i = \{\lambda v:v\in P_i\}$, and we define a function
  $V_{P_1,\ldots,P_n}:\R^n_{\geq 0}\to\R$ by
  \[ V_{P_1,\ldots,P_n}(\lambda_1,\ldots,\lambda_n) = 
  \vol(\lambda_1 P_1 + \cdots + \lambda_n P_n), \]
  where $\vol$ is a Euclidean volume form on $N_\R\cong\R^n$ normalized such that
  the volume of a fundamental domain for the lattice $N$ is one.  It 
  well-known that $V_{P_1,\ldots,P_n}$ is a homogeneous polynomial of
  degree $n$ in $\lambda_1,\ldots,\lambda_n$.  The \emph{mixed volume} 
  $\MV(P_1,\ldots,P_n)$ is defined to be the
  coefficient of the $\lambda_1\cdots\lambda_n$-term of
  $V_{P_1,\ldots,P_n}$. 
\end{defn}

\begin{eg} Fixing \label{eg:mv.det}
  a basis, we identify $N$ with $\Z^n$.
  Suppose that $P_i$ is the line segment connecting points
  $v_i,v_i'\in N=\Z^n$, and let 
  $w_i = v_i - v_i'$.  Then
  \[ V_{P_1,\ldots,P_n}(\lambda_1,\ldots,\lambda_n) 
  = |\det(\lambda_1 w_1,\ldots,\lambda_n w_n)|
  = \lambda_1\cdots\lambda_n |\det(w_1,\ldots,w_n)|, \]
  where $\det(w_1,\ldots,w_n)$ is the determinant of the matrix whose
  $i$th column is the column vector $w_i\in\Z^n$.  Therefore
  $\MV(P_1,\ldots,P_n) = |\det(w_1,\ldots,w_n)|$ in this case.  
\end{eg}

\begin{defn}
  Let $P\subset N_\R$ be an integral $\Gamma$-affine pointed polyhedron, let
  $f_1,\ldots,f_n\in K\angles{U_P}$, let $Y_i = V(f_i)$, and let 
  $Y = \bigcap_{i=1}^n Y_i$.  Let $v\in N_\Gamma\cap P$.  The
  \emph{intersection multiplicity of $Y_1,\ldots,Y_n$ over $v$}, denoted
  $i(v,\, Y_1\cdots Y_n)$, is defined to be the length of the space
  $Y\cap U_{\{v\}}$:
  \[ i(v,\,Y_1\cdots Y_n) \coloneq 
  \dim_K \Gamma(Y\cap U_{\{v\}}, \sO_{Y\cap U_{\{v\}}}). \]
\end{defn}

Note that $i(v,\,Y_1\cdots Y_n) < \infty$ if and only if 
$\dim(Y\cap U_{\{v\}})=0$, in which case
\[ i(v\,,Y_1\cdots Y_n) = \sum_{\trop(\xi)=v} \dim_K(\sO_{Y,\xi}). \]
Note also that 
$i(v,\,Y_1\cdots Y_n)$ only depends on the images of $f_1,\ldots,f_n$
in $K\angles{U_{\{v\}}}$.
The relation between $i(v,\,Y_1\cdots Y_n)$ and the Newton complexes of
$f_1,\ldots,f_n$ is as follows:

\begin{thm}[Katz; Osserman-Payne] Let \label{thm:multiplicity.lp}
  $f_1,\ldots,f_n\in K[M]$ and let
  $v\in\bigcap_{i=1}^n\Trop(f_i)$ be an isolated point.  For
  $i=1,\ldots,n$ let  $Y_i = V(f_i)$ and let
  $\gamma_i = \pi(\conv(\vertices_v(f_i)))\in\New(f_i)$
  be the polytope corresponding to $v\in\Trop(f_i)$ as
  in~\parref{par:newton.complex}.  Then
  \[ i(v,\,Y_1\cdots Y_n) = \MV(\gamma_1,\ldots,\gamma_n). \]
\end{thm}

Note that $v\in N_\Gamma$ since $\{v\}$ is an integral $\Gamma$-affine polytope.
The mixed volume $\MV(\gamma_1,\ldots,\gamma_n)$ is the stable tropical
multiplicity of the point $v\in\bigcap_{i=1}^n\Trop(f_i)$, as we will
discuss in~\S\ref{sec:non.transverse}.

\begin{remsub} Bernstein's \label{rem:bernstein.history} 
  theorem~\cite{bernstein:number_of_roots} can
  be seen as the generic coefficient (i.e.\ trivial valuation) case
  of~\parref{thm:multiplicity.lp}; see also~\cite[Chapter~3]{sturmfels:systems}
  for a proof in the bivariate case.  Theorem~\parref{thm:multiplicity.lp}
  is due to E.~Katz in the case of a nontrivial
  discrete valuation~\cite[Theorem~8.8]{katz:toolkit}.  Osserman and
  Payne~\cite{osserman_payne:lifting} develop an intersection theory over
  non-noetherian valuation rings in order to remove the
  noetherian hypothesis (in addition to proving a very general
  compatibility theorem).
\end{remsub}

\begin{eg} The \label{eg:multeg}
  following type of example arises in the analysis of the zeros of the
  logarithm of a $p$-divisible group over $\sO_K$ as in~\cite{jdr:thesis}.
  Let $p$ be a prime number and suppose that $\val(p)=1$.
  Let $N = M = \Z^2$ and let $x=x^{(-1,0)},\,y=x^{(0,-1)}\in K[M]$ as
  in~\parref{eg:tropical.line.1}.  Let
  \[ f_1 = px + x^p + y^p \sptxt{and} f_2 = y + x^p + y^p \quad\in\quad K[M]. \]
  The tropicalizations and Newton complexes of $f_1$ and $f_2$ are drawn in
  Figure~\ref{fig:multeg}.  They intersect in the two points
  $v_1 = (\frac 1{p-1},\frac p{p-1})$ and
  $v_2 = (\frac 1{p-1},0)$.  For $i,j=1,2$ let $\gamma_{i,j}$ be the cell
  in $\New(f_i)$ corresponding to $v_j\in\Trop(f_i)$ as
  in~\parref{par:newton.complex}.  Then
  \begin{align*} \gamma_{1,1} &= \conv\{(-1,0),(-p,0)\} &
    \gamma_{1,2} &= \conv\{(-1,0),(-p,0),(0,-p)\} \\
    \gamma_{2,1} &= \conv\{(0,-1),(-p,0)\} &
    \gamma_{2,2} &= \conv\{(0,-1),(0,-p)\}
  \end{align*}
  as indicated in the figure.  Since $\gamma_{1,1}$ and $\gamma_{2,1}$ are
  line segments, their mixed volume can be calculated as
  in~\parref{eg:mv.det}:
  \[ \MV(\gamma_{1,1},\gamma_{2,1}) = \left|\det
    \begin{pmatrix}
      p-1 & p \\ 0 & -1
    \end{pmatrix}\right| = p-1. \]
  By~\parref{thm:multiplicity.lp} there are exactly $p-1$ common zeros
  $\xi = (\xi_1,\xi_2)$ of $f_1,f_2$, counted with multiplicity, such that
  $\val(\xi_1) = \frac 1{p-1}$ and $\val(\xi_2) = \frac p{p-1}$.
  The calculation of $\MV(\gamma_{1,2},\gamma_{2,2})$ requires some
  grade-school geometry since $\gamma_{1,2}$ is all of $|\New(f_1)|$:
  we have
  \[ \vol(\lambda_1\gamma_{1,2} + \lambda_2\gamma_{2,2})
  = \lambda_1\lambda_2\big((p-1)^2 + p-1\big)
  + \lambda_1^2\frac{p(p-1)}2 \]
  for $\lambda_1,\lambda_2\geq 0$, so
  $\MV(\gamma_{1,2},\gamma_{2,2}) = p^2 - p$.
  Hence there are exactly $p^2-p$ common zeros $\xi$ of $f_1,f_2$ such
  that $\val(\xi_1) = \frac 1{p-1}$ and $\val(\xi_2) = 0$.

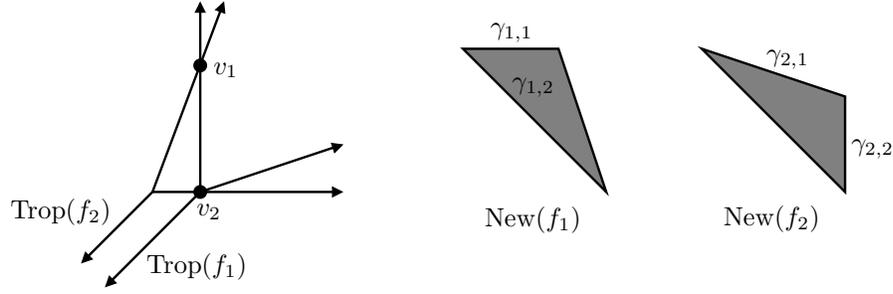
  \genericfig[ht]{multeg}{The tropicalizations and Newton complexes of
    $f_1=px+x^p+y^p$ and $f_2=y+x^p+y^p$.  See~\parref{eg:multeg}.}

\end{eg}

The goal of this section is to derive the following generalization
of~\parref{thm:multiplicity.lp}:

\begin{thm} Let \label{thm:multiplicity.ps}
  $P\subset N_\R$ be an integral $\Gamma$-affine pointed polyhedron,
  let $f_1,\ldots,f_n\in K\angles{U_P}$, and let 
  $v\in\bigcap_{i=1}^n\Trop(f_i)$ be an isolated point contained in the
  interior of $P$.  For $i=1,\ldots,n$ let $Y_i = V(f_i)$ and let
  $\gamma_i\in\New(f_i)$ be the polytope corresponding to $v\in\Trop(f_i)$ as
  in~\parref{thm:multiplicity.lp}.  Then
  \[ i(v,\,Y_1\cdots Y_n) = \MV(\gamma_1,\ldots,\gamma_n). \]
\end{thm}

\begin{eg}[The theorem of the Newton polygon]
  Let \label{eg:multiplicity.ps.dim1}
  $N = M = \Z$ and let $x = x^{(-1)}\in K[M]$ as
  in~\parref{eg:dim1.trop}.  Let $r\in\Gamma$ and $\rho=\exp(-r)$ and let 
  $P = [r,\infty)\subset N_\R$, so 
  $U_P = \B^1_K(\rho)$ and 
  $K\angles{U_P} = \{\sum_{n\geq 0} a_n x^n:|a_n|\rho^n\to 0 \}$
  as in~\parref{eg:trop.B1}.  Let $f\in K\angles{U_P}$ be nonzero, and
  assume for simplicity that $f(x)\neq 0$.  As explained
  in~\parref{eg:newtonpoly1}, a number $v > r$ is the valuation of a
  zero of $f$ if and only if $\conv(\vertices_v(f))$ is a line segment in
  the lower convex hull $\NP'(f)$ of $H(f,\{0\})$, in which case the slope
  of the segment is   $v$.  The polytope 
  $\gamma = \pi(\conv(\vertices_v(f)))\in\New(f)$ is the projection of
  $\conv(\vertices_v(f))$ onto the $x$-axis; it is a line
  segment whose length $L$ is exactly the horizontal length of 
  $\conv(\vertices_v(f))$.  See Figure~\ref{fig:newtonpoly}.
  Therefore~\parref{thm:multiplicity.ps} implies
  that there are exactly $L$ zeros $\xi$ of $f$, counted with
  multiplicity, such that $\val(\xi) = v$.  
\end{eg}

We will use the following consequence of~\parref{thm:controots2}:

\begin{cor}[to~\ref{thm:controots2}] Let
  \label{cor:controots2.special.case} 
  $P\subset N_\R$ be an integral $\Gamma$-affine pointed polyhedron and let
  $f_1,\ldots,f_n\in K\angles{U_P,t}\coloneq K\angles t\angles{U_P}$.  Let
  $Y_i\subset U_P\times\B^1_K$ be the subspace defined by $f_i$, let
  $\pi: Y_i\to\B^1_K$ be the projection onto the second factor, and for 
  $t_0\in|\B^1_K|$ let 
  $Y_{i,t_0} = \pi\inv(t_0)\subset\kappa(t_0)\hat\tensor_K U_P$.  Then for 
  any $v\in N_\Gamma\cap P$ such that 
  $i(v,\,Y_{1,0}\cdots Y_{n,0})<\infty$ there exists 
  $\epsilon\in|K^\times|$ such that
  \[ i(v,\,Y_{1,t_0}\cdots Y_{n,t_0}) = i(v,\,Y_{1,0}\cdots Y_{n,0}) 
  \sptxt{whenever} |t_0|\leq\epsilon. \]
\end{cor}

We will also need a device for approximating a power series by a sequence
of polynomials fitting into a one-parameter family:

\begin{lem} Fix \label{lem:polyapprox}
  a nonzero element $\varpi\in\fm_K$.  Let $P\subset N_\R$ be an
  integral $\Gamma$-affine pointed polyhedron with cone of unbounded directions
  $\sigma$, and let $f\in K\angles{U_P}$ be nonzero.  There is a power series
  $g\in K\angles{U_P,t}$ such that $g_0 = f$ and
  $g_{\varpi^m}\in K[S_\sigma]$ for all $m\geq 1$, where for
  $t_0\in|\B^1_K|$ we let $g_{t_0}$ denote the specialization of $g$ at
  $t=t_0$.   In particular, $g_{\varpi^m}\to f$ in $K\angles{U_P}$ as
  $m\to\infty$. 
\end{lem}

\pf For $m\geq 1$ we define
\[ q_m(t) = (t - \varpi)(t - \varpi^2)\cdots(t - \varpi^m)(t - (-1)^m \varpi^{-m(m+1)/2})
\in K[t], \]
so $q_m(\varpi^i) = 0$ for $i = 1,\ldots,m$ and $q_m(0) = 1$.  
Write $f = \sum_{u\in S_\sigma} a_u x^u$.  Choose a denumeration
$\delta: S_\sigma\isom\Z_{\geq0}$, and find a sequence of numbers
$m_N$, tending to $\infty$ as $N\to\infty$, such that 
\[ |q_{m_{\delta(u)}}|\cdot |a_u x^u|_P\To 0 \sptxt{as}
\delta(u)\To\infty, \]
where $|q_{m_{\delta(u)}}|$ denotes the supremum norm
of $q_{m_{\delta(u)}}$ in $K\angles t$.
Set \[ g = \sum_{u\in S_\sigma} q_{m_{\delta(u)}}(t) \,
a_u x^u\in K\angles{U_P,t}. \]
By construction, $g_{\varpi^m}\in K[S_\sigma]$ for all $m\geq 1$ 
and $g_0 = f$.\qed

\pf[of~\parref{thm:multiplicity.ps}]
It follows from~\parref{prop:automatically.finite} as applied to a
small polytope containing $v$ in its interior that 
$i(v,\,Y_1\cdots Y_n) < \infty$.  For $i=1,\ldots,n$ let
$g_i\in K\angles{U_P,t}$ be as in~\parref{lem:polyapprox}, 
so $p_{i,m} \coloneq g_{i,\varpi^m} \in K[M]$ for all
$m\geq 1$, and $p_{i,m}\to f_i$ as $m\to\infty$.  Let
$Y_{i,m} = V(p_{i,m})$.
By~\parref{lem:inu.finiteness} we have 
$\vertices_P(f_i) = \vertices_P(p_{i,m})$ for $m\gg 0$, so
$\New(f_i) = \New(p_{i,m})$ and 
$\Trop(f_i)\cap P = \Trop(p_{i,m})\cap P$ for all $i$ and all 
$m\gg 0$ (see~\parref{rem:new.trop.agree}); hence if 
$\gamma_{i,m} = \pi(\conv(\vertices_v(p_{i,m})))$ then
$\gamma_{i,m} = \gamma_i$ for $m\gg 0$.
By~\parref{cor:controots2.special.case} we likewise have
$i(v,\,Y_1\cdots Y_m) = i(v,\,Y_{1,m}\cdots Y_{n,m})$ for $m\gg 0$.  
Thus for $m\gg 0$, 
\[ i(v,\,Y_1\cdots Y_n) = i(v,\,Y_{1,m}\cdots Y_{n,m})
= \MV(\gamma_{1,m},\ldots,\gamma_{n,m}) = \MV(\gamma_1,\ldots,\gamma_n) \]
by~\parref{thm:multiplicity.lp}.\qed

\begin{rem}
  It would be interesting to investigate a more general relationship
  between the local intersection theory of tropical varieties with a
  non-Archimedean toric intersection theory along the lines of
  Osserman and Payne's work.
\end{rem}

\section{Application: tropically non-proper complete intersections}
\label{sec:non.transverse}

\paragraph Let $f_1,\ldots,f_n\in K[M]$ be nonzero, let 
$Y = \bigcap_{i=1}^n V(f_i)$, and let $C$ be a connected 
component of $\bigcap_{i=1}^n\Trop(f_i)\subset N_\R$.  If $C = \{v\}$
consists of a single point then~\parref{thm:multiplicity.lp} 
calculates the sum $\sum_{\trop(\xi)=v}\dim_K(\sO_{Y,\xi})$ in terms of a
mixed 
volume.  The main goal of this section is to generalize this result to the
case when $C$ is arbitrary.  More precisely, after 
taking the closure $\bar C$ of $C$ in an appropriate compactification
$N_\R(\Delta)$ of $N_\R$ and taking the closure $\bar Y$ of $Y$ in the
corresponding toric variety $X(\Delta)$, the size of the algebraic intersection
$\sum_{\trop(\xi)\in\bar C}\dim_K(\sO_{\bar Y,\xi})$ lying above $\bar C$ can be
calculated in terms of stable tropical intersection multiplicities.
See~\parref{thm:non.proper.multiplicity}. 
The compactification step is necessary: see~\parref{eg:non.proper.int}.
Along the way we will obtain a new proof that the stable tropical
intersection multiplicity is well-defined in the case of a dimension-zero
complete intersection.  

The idea is to translate each $V(f_i)$ by a generic point of the torus in order
to reduce our problem to~\parref{thm:multiplicity.lp}; the key ingredient
is the continuity of roots result~\parref{thm:controots1} which allows us
to relate the intersection multiplicities before and after the
translation.  It is important to notice that one is led to work with
families of translations parametrized by an \emph{affinoid} subspace of a
torus and not by a scheme; cf.~\parref{eg:trop.deform}.  This
rigid-analytic deformation technique is what makes the algebraic
result~\parref{thm:non.proper.multiplicity}  possible.

\begin{remsub}
  We have chosen work with Laurent polynomials in this section mainly for
  simplicity of formulation; most of the ideas also apply to power series.
\end{remsub}

\paragraph[Stable tropical intersection multiplicities]
There is a rich intersection theory of tropical varieties,
developed in many papers including~\cite{allermann_rau:first_steps,katz:tropical_int_thy,mikhalkin:tropical_geometry,sturmfels_et_al:first_steps}.
Basic to all of these theories is the notion of the stable tropical
intersection, which is entirely combinatorial.  As we are restricting
ourselves to the case of dimension-zero complete 
intersections, we will take a pedestrian approach and give a direct
definition of the stable
tropical intersection multiplicity of $n$ hypersurfaces in an
$n$-dimensional torus along a connected component. 

\begin{defn}
  Let $P = \bigcap_{i=1}^r \{v\in N_\R:\angles{u_i,v}\leq a_i\}$ be an integral
  $\Gamma$-affine polyhedron in $N_\R$, where $u_i\in M$ and
  $a_i\in\Gamma$.  A \emph{thickening} 
  of $P$ is a polyhedron of the form
  \[ P' = \bigcap_{i=1}^r \{v\in N_\R~:~\angles{u_i,v}\leq a_i+\epsilon\} \]
  for $\epsilon > 0$ contained in $\Gamma$.
  More generally, if $\Pi$ is a polyhedral complex
  then a \emph{thickening} $\cP$ of $\Pi$ is a set of the form
  $\cP = \{ P':P\in\Pi \}$,
  where $P'$ denotes a thickening of $P$.  We set
  \[ |\cP| = \bigcup_{P'\in\cP} P' \sptxt{and}
  \mathring\cP = \bigcup_{P'\in\cP}(P')^\circ\subset|\cP|^\circ.\]
  If $\cP' = \{P'':P\in\Pi\}$ is a second thickening of $\Pi$, we say
  that $\cP'$ \emph{dominates} $\cP''$ if $P''\subset(P')^\circ$ for all
  $P\in\Pi$. 
\end{defn}

\begin{rem} \label{rem:thickenings}
  \begin{enum}
  \item If $P'$ is a thickening of $P$ then $P$ is contained in the
    interior $(P')^\circ$ of $P'$, and hence if $\cP$ is a thickening of
    $\Pi$ then $|\Pi|\subset\mathring\cP\subset|\cP|^\circ$.
  \item If $P'$ is a thickening of $P$ then $\cU(P) = \cU(P')$.
  \item If $\Pi$ is a polyhedral complex and $C\subset|\Pi|$ is a
    connected component then $C$ is the support of the subcomplex $\Pi_C$
    of $\Pi$ whose cells are contained in $C$.  There is a thickening
    $\cP$ of $\Pi_C$ such that $|\cP|\cap|\Pi| = C$.
  \end{enum}
\end{rem}

Recall~\parref{par:trop.poly.cplx.1} that if $f_1,\ldots,f_n\in K[M]$ are
nonzero then each $\Trop(f_i)$ 
is (the support of) a canonical polyhedral complex, and therefore
$\bigcap_{i=1}^n \Trop(f_i)$ is also canonically a polyhedral complex. 
The following lemma is standard, but we include a proof for
completeness: 

\begin{lem}[Moving lemma] Let \label{lem:moving}
  $f_1,\ldots,f_n\in K[M]$ be nonzero, let $C$ be a connected
  component of $\bigcap_{i=1}^n \Trop(f_i)$, and let 
  $\cP$ be a thickening of (the complex underlying) $C$ such
  that $|\cP|\cap\bigcap_{i=1}^n\Trop(f_i) = C$.  Then there exist
  $v_1,\ldots,v_n\in N$ and $\epsilon\in\R_{>0}\cap\Gamma$ such that for
  all $t\in(0,\epsilon]$, the intersection
  \[ |\cP|\cap\bigcap_{i=1}^n \big(\Trop(f_i) + t v_i\big) \]
  is a finite set of points contained in $\mathring \cP$.
\end{lem}

\pf Each $\Trop(f_i)$ is a subset of a hyperplane arrangement in
$N_\R\cong\R^n$, so we can find $v_i$ and $\epsilon$ such that 
$|\cP|\cap\bigcap_{i=1}^n (\Trop(f_i) + t v_i)$ is a finite set of points
for $t\leq\epsilon$,  since the intersection of $n$ affine hyperplanes in
$\R^n$ generically contains 
zero or one points.  Furthermore, the union of the boundaries of the
polyhedra in $\cP$ is also contained in a hyperplane arrangement, so we can
choose $\epsilon$ such that 
$|\cP|\cap \bigcap_{i=1}^n(\Trop(f_i)+t v_i)\subset\mathring \cP$ for
$t\leq\epsilon$ as well since $n+1$ affine hyperplanes in $\R^n$
generically have no points of intersection.\qed

\paragraph Let \label{par:torus.translation}
$T = \Spec(K[M])$, let $v\in N_\Gamma$, and
choose $\xi\in T(K')$ with $\trop(\xi) = v$, where $K'$ is a suitable finite
extension of $K$.  Then $\xi$ induces the translation automorphism 
$\eta\mapsto\xi\cdot\eta$ of $T_{K'} = K'\tensor_K T$, which corresponds
to the automorphism $x^u\mapsto x^u(\xi) x^u$ of $K'[M]$.  We denote the
image of $f\in K'[M]$ under this automorphism by $\xi\cdot f$.
Since $\trop(\xi\cdot\eta) = \trop(\eta) + v$ we have
$\Trop(\xi\cdot f) = \Trop(f) + v$, and since $\Trop(f)$ and $\New(f)$ only
depend on the valuations of the coefficients of $f$, the complexes
$\Trop(\xi\cdot f)$ and $\New(\xi\cdot f)$ are independent of the choice
of $\xi\in\trop\inv(v)$.

\begin{defn}
  Let $f_1,\ldots,f_n\in K[M]$ be nonzero and let 
  $v\in\bigcap_{i=1}^n\Trop(f_i)$ be an isolated point.  The 
  \emph{stable tropical intersection multiplicity of
    $\Trop(f_1),\ldots,\Trop(f_n)$ at $v$} is defined to be
  \[ i\big(v,\Trop(f_1)\cdots\Trop(f_n)\big) = \MV(\gamma_1,\ldots,\gamma_n), \]
  where $\gamma_i\in\New(f_i)$ is the polytope corresponding to
  $v\in\Trop(f_i)$ as in~\parref{thm:multiplicity.lp}.  Now let
  $C\subset\bigcap_{i=1}^n\Trop(f_i)$ be a connected component and let
  $\cP$, $v_1,\ldots,v_n\in N$, and $\epsilon\in\R_{>0}\cap\Gamma$ be as
  in~\parref{lem:moving}.  The
  \emph{stable tropical intersection multiplicity of
    $\Trop(f_1),\ldots,\Trop(f_n)$ along $C$}
  is defined to be
  \[\begin{split}& i\big(C,\,\Trop(f_1)\cdots\Trop(f_n)\big) \\
  &\quad= \sum\bigg\{i\big(v,\,(\Trop(f_1)+\epsilon v_1)\cdots(\Trop(f_n)+\epsilon v_n)\big)
  ~:~ v\in |\cP|\cap\bigcap_{i=1}^n(\Trop(f_i)+\epsilon v_i) \bigg\},
  \end{split}\]
  which makes sense by~\parref{lem:moving}
  and~\parref{par:torus.translation}. 
\end{defn}

See~\parref{eg:non.proper.int} for an example.

\begin{remsub}
  The above definition of $i(C,\,\Trop(f_1)\cdots\Trop(f_n))$ agrees with the
  sum of the multiplicities of the points of the stable intersection 
  $\Trop(f_1)\cdots\Trop(f_n)$ contained in $C$;
  see~\cite[Theorem~4.6]{sturmfels_tevelev:elimination}. 
  Ordinarily one proves that this number is well-defined using the 
  balancing condition on a tropical variety, but it will also follow
  from~\parref{thm:non.proper.multiplicity}. 
\end{remsub}

\begin{notn}
  For a nonzero Laurent polynomial $f = \sum a_u x^u\in K[M]$ we denote
  the normal fan to  
  $|\New(f)|=\conv\{u:a_u\neq 0\}$ by $\Delta(f)$.
\end{notn}

Note that $\Delta(f)$ is a complete fan.

\begin{egsub}
  Let $M = N = \Z^2$ and let $x=x^{(-1,0)},\,y=x^{(0,-1)}\in K[M]$ as
  in~\parref{eg:tropical.line.1}.  Let $\lambda\in K^\times$ have 
  valuation $2$ and let $f = 1 + x + y + \lambda xy\in K[M]$, so
  $\Trop(f)$ and $\New(f)$ are drawn in Figure~\ref{fig:complexes}.
  The unbounded cells of $\Trop(f)$ are labeled $P_1,P_2,P_4,P_5$ in the
  figure; their cones of unbounded directions are the positive and negative
  coordinate axes.  The
  Newton polytope $|\New(f)|$ is the unit square, so $\Delta(f)$ is the
  fan of Figure~\ref{fig:normalfan}.  Note that the positive and negative
  coordinate axes are cones of $\Delta(f)$.
\end{egsub}

The positive-dimensional cones in $\Delta(f)$ represent the directions in
which $\Trop(f)$ is unbounded:

\begin{lem} \label{lem:Delta.unbounded}
  \begin{enum}
  \item Let $f = \sum a_u x^u\in K[M]$ be a nonzero Laurent polynomial
    and let $P$ be a cell of $\Trop(f)$.  Then $\cU(P)\in\Delta(f)$.
  \item Let $f_1,\ldots,f_n\in K[M]$ be nonzero and let
    $P$ be a cell of $\bigcap_{i=1}^n \Trop(f_i)$.  Then
    $\cU(P)$ is a cone of $\bigcap_{i=1}^n\Delta(f_i)$.
  \end{enum}
\end{lem}

\pf The second part follows from the first by~\parref{lem:face.unbounded}
and~\parref{lem:unbounded.intersection}, so we proceed with~(i).
By definition~\parref{par:trop.poly.cplx.1} there is a point
$v\in\Trop(f)$ such that 
\[ P = \gamma_v = \{ v'\in N_\R~:~\vertices_{v'}(f)\supset\vertices_v(f) \}. \]
Let $\check\gamma_v = \pi(\conv(\vertices_v(f)))\in\New(f)$ be the dual
cell as in~\parref{par:newton.complex}. 
We claim that 
\begin{equation} \label{eq:Delta.unbdd.1}
 \cU(P) = \{ v\in N_\R~:~\face_v(|\New(f)|)\supset\check\gamma_v \}. 
\end{equation}
First notice that the right side of~\eqref{eq:Delta.unbdd.1} is the
cone of $\Delta(f)$ corresponding to the minimal face of $|\New(f)|$
containing $\check\gamma_v$ (so $P$ is unbounded if and only if
$\check\gamma_v$ is contained in the boundary of $|\New(f)|$).
Let $v\in\cU(P)$ and let $(u,\val(a_u))\in\vertices_v(f)$; we want to show
that $\angles{u,v} = \max\{\angles{u',v}~:~a_{u'}\neq 0\}$.
Fix $v_1\in P$.  For any $\lambda\in\R_{\geq0}$ we have
$v_1 + \lambda v\in P$ by~\parref{lem:bounded.faces}, i.e.
$\vertices_{v_1+\lambda v}(f)\supset\vertices_v(f)$, so
\begin{equation} \label{eq:Delta.unbdd.2}
 \val(a_u) - \angles{u,v_1} - \lambda\angles{u,v} 
= \min\{\val(a_{u'}) - \angles{u',v_1} - \lambda\angles{u',v}~:~a_{u'}\neq 0\}. 
\end{equation}
If there were some $u'$ with $a_u\neq 0$ and $\angles{u',v}>\angles{u,v}$
then we could make~\eqref{eq:Delta.unbdd.2} false by taking
$\lambda\gg 0$.  This proves one inclusion of~\eqref{eq:Delta.unbdd.1}.
On the other hand, if $v\in N_\R$ satisfies
$\face_v(|\New(f)|)\supset\check\gamma_v$ then a similar argument
shows that $v_1+v\in P$ for any $v_1\in P$, so the other inclusion also
follows from~\parref{lem:bounded.faces}.\qed

Hence if $f_1,\ldots,f_n\in K[M]$ are nonzero Laurent polynomials then
$N_\R(\bigcap_{i=1}^n\Delta(f_i))$ is a natural compactification of $N_\R$
in which to take the closure of $\bigcap_{i=1}^n\Trop(f_i)$.

\begin{rem}
  Let $f_1,\ldots,f_n\in K[M]$ be nonzero, and suppose that the fan
  $\Delta = \bigcap_{i=1}^n\Delta(f_i)$ is not pointed.  Then there is a proper
  subspace $M'_\R\subset M_\R$ and elements $u_i\in M_\R$ such that
  $\New(f_i)\subset u_i + M'_\R$.  In this case, a Minkowski sum of cells
  of the $\New(f_i)$ is also contained in a translate of $M'_\R$, so all
  mixed volumes appearing in the definition of
  $i(C,\,\Trop(f_1)\cdots\Trop(f_n))$ are zero for any connected component
  $C\subset\bigcap_{i=1}^n\Trop(f_i)$.  This is the ``overdetermined'' or
  ``degenerate'' case, and for this reason we will generally assume that
  $\Delta$ is pointed.
\end{rem}

The rest of this section is devoted to proving that the purely
combinatorially defined quantity
$i(C,\,\Trop(f_1)\cdots\Trop(f_n))$ calculates algebraic intersection
multiplicities in the following sense:

\begin{thm} Let  \label{thm:non.proper.multiplicity}
  $f_1,\ldots,f_n\in K[M]$ be nonzero Laurent polynomials,
  and assume that the fan $\Delta = \bigcap_{i=1}^n\Delta(f_i)$ is pointed.  Let
  $C\subset\bigcap_{i=1}^n\Trop(f_i)$ be a connected component and let
  $\bar C$ be the closure of $C$ in $N_\R(\Delta)$.  Let $Y_i$ be the
  closure of $V(f_i)$ in $X(\Delta)$ and  let $Y = \bigcap_{i=1}^n Y_i$.
  Then
  \begin{equation}  \label{eq:non.transverse.mult}
    i\big(C,\Trop(f_1)\cdots\Trop(f_n)\big)
    = \sum_{\trop(\xi)\in\bar C} \dim_K(\sO_{Y,\xi}) 
  \end{equation}
  if the right side is finite.
\end{thm}

\begin{remsub}
  If $C$ is a polyhedron then the right side
  of~\eqref{eq:non.transverse.mult} is automatically finite
  by~\parref{prop:automatically.finite} and~(\ref{prop:trop.closure},iv)
  below. 
\end{remsub}

\begin{cor} Let \label{cor:well.defined}
  $f_1,\ldots,f_n\in K[M]$ be nonzero and let
  $C\subset\bigcap_{i=1}^n \Trop(f_i)$ be a connected component.  Then
  $i(C,\,\Trop(f_1)\cdots\Trop(f_n))$ is independent of
  all choices. 
\end{cor}

\begin{remsub} 
  The above corollary is a purely tropical result: it only
  depends on $f_1,\ldots,f_n$ through the valuations of their
  coefficients, and hence can be stated in terms of 
  tropical polynomials.  Thus it can be seen as an application of rigid
  geometry to ``pure'' tropical geometry.
\end{remsub}

\begin{eg} Let \label{eg:non.proper.int}
  $M = N = \Z^2$ and let $x=x^{(-1,0)},\,y=x^{(0,-1)}\in K[M]$ as
  in~\parref{eg:tropical.line.1}.  Choose 
  $\alpha,\beta\in K^\times$ of valuation zero, and let
  \[ f_1 = x + y + 1 \sptxt{and} f_2 = \alpha x + \beta y + 1
  \quad\in\quad K[M]. \]
  A picture of $\Trop(f_1) = \Trop(f_2)$ and $\New(f_1)=\New(f_2)$ can be
  found in Figure~\ref{fig:trop_new} (with $\lambda=1$).  Hence
  $C = \Trop(f_1)\cap\Trop(f_2)$ is a connected component, and one easily
  calculates that $i(C,\,\Trop(f_1)\cdot\Trop(f_2))=1$.  The fan
  $\Delta = \Delta(f_1) = \Delta(f_2)$ and the completion $N_\R(\Delta)$
  are described in~\parref{eg:proj.plane} and drawn in
  Figure~\ref{fig:tropdelta}; the associated toric variety is 
  $X(\Delta) = \bP^2_K$.  The closure $\bar C$ of 
  $C = \Trop(f_1) = \Trop(f_2)$ is 
  \[ \bar C = C\djunion\{[0:0:\infty]\}\djunion\{[0:\infty:0\}
  \djunion\{[\infty:0:0]\} \]
  with the notation in~\parref{eg:trop.projspace}.

  Algebraically, let $Y_i$ be the closure of $V(f_i)$ in 
  $\bP^2_K = X(\Delta)$ and 
  let $Y = Y_1\cap Y_2$ as in~\parref{thm:non.proper.multiplicity}. 
  Then $Y$ consists of the single point 
  \[ (\xi,\eta) = \bigg(\frac{\beta-1}{\alpha-\beta},\,
  \frac{\alpha-1}{\beta-\alpha}\bigg) \]
  as long as $(\alpha,\beta)\neq(1,1)$.  We can choose $\alpha$ and
  $\beta$ so that $\Trop(Y)$ is located anywhere on $\bar C$: 
  \begin{bullets}
  \item If $\val(\beta-\alpha)\gg 0$ but 
    $\val(\beta-1)=\val(\alpha-1)=0$ then $\Trop(Y)$ is a point on the ray
    $R_1$ of Figure~\ref{fig:trop_new}.
  \item If $\alpha=\beta$ then $\Trop(Y) = \{[0:0:\infty]\}$.
  \item If $\val(\alpha-1)\gg 0$ but $\val(\beta-1) = 0$ then 
    $\Trop(Y)$ is a point on the ray $R_2$ of Figure~\ref{fig:trop_new}.
  \item If $\alpha=1$ then $\Trop(Y) = \{[0:\infty:0]\}$.
  \item If $\val(\beta-1)\gg 0$ but $\val(\alpha-1)=0$ then
    $\Trop(Y)$ is a point on the ray $R_3$ of Figure~\ref{fig:trop_new}.
  \item If $\beta=1$ then $\Trop(Y)=\{[\infty:0:0]\}$.
  \end{bullets}
  Hence we need to consider all points $\xi\in|Y|$ with
  $\trop(\xi)\in\bar C$ in~\eqref{eq:non.transverse.mult}.

\end{eg}

We will prove~\parref{thm:non.proper.multiplicity}
and~\parref{cor:well.defined} below.  First we investigate the
relationship between the closure of a subscheme of a torus inside a toric
variety and the closure of its tropicalization.  For a different treatment
see~\cite[\S3]{osserman_payne:lifting}. 

\begin{prop} \label{prop:trop.closure}
  \begin{enum}
  \item Let $A$ be an integral domain, let $f\in A[M]$ be nonzero, let
    $\sigma'\in\Delta(f)$, and let $\sigma$ be an integral pointed cone
    contained in $\sigma'$.  Then there is a vertex $u\in M$ of
    $|\New(f)|$, depending only on $|\New(f)|$ and $\sigma$,
    such that $A[S_\sigma] x^{-u} f = (A[M]f)\cap A[S_\sigma]$.
  \item With the notation in~(i), suppose that $A = K$.  Let $T = \Spec(K[M])$, let 
    $Y = V(f)\subset T$, and let $\bar Y$ be the closure of $Y$ in
    $X(\sigma)$.  Then $\bar Y$ is cut out by $x^{-u}f$ and
    $\Trop(\bar Y,N_\R(\sigma))$ is the closure of $\Trop(f,N_\R)$ in
    $N_\R(\sigma)$. 
  \item Let $f_1,\ldots,f_n\in K[M]$ be nonzero and suppose that
    $\Delta = \bigcap_{i=1}^n \Delta(f_i)$ is pointed.
    Then the closure of $\bigcap_{i=1}^n\Trop(f_i)$ in $N_\R(\Delta)$ is
    equal to $\bigcap_{i=1}^n\bar{\Trop(f_i)}$, where
    $\bar{\Trop(f_i)}$ is the closure of $\Trop(f_i)$ in $N_\R(\Delta)$.
  \item With the notation in~(iii), let $Y_i$ be the closure of $V(f_i)$
    in $X(\Delta)$ and let $Y = \bigcap_{i=1}^n Y_i$.  Then
    $\Trop(Y)$ is contained in the closure of 
    $\bigcap_{i=1}^n\Trop(f_i)$ in $N_\R(\Delta)$.
  \end{enum}
\end{prop}

\pf We may assume that $\sigma'$ is a maximal cone of $\Delta(f)$, so
there is a vertex $u$ of $|\New(f)|$ such that
\[ \sigma' = \{ u'\in N_\R~:~u\in\face_{u'}(|\New(f)|) \}. \]
We claim that this $u$ works.  One checks that if $f' = x^{-u}f$ then 
$f'\in A[S_{\sigma'}]\subset A[S_{\sigma}]$.  Let $g\in A[M]$, and suppose
that $fg\in A[S_\sigma]$.  We have  $|\New(f'g)| = |\New(f')| + |\New(g)|$
by~\cite[Proposition~6.1.2(b)]{gkz:discriminants}, so since
$|\New(fg)|\subset\sigma^\vee$ and $0\in|\New(f')|$ we have 
$|\New(g)|\subset\sigma^\vee$, i.e.\ $g\in A[S_\sigma]$.  This
proves~(i).

The closure $\bar Y$ of $Y$ is the hypersurface in $X(\sigma)$ cut out
by $f'$ by the above.  Since $\Trop(f',N_\R) = \Trop(f,N_\R)$ we may
replace $f$ by $f'$.  Since $\Trop(f, N_\R(\sigma))$ contains
$\Trop(f,N_\R)$ it also contains the closure
$\bar{\Trop(f,N_\R)}$.  Let
$\tau\prec\sigma$ be nonzero, let $N_\R'=N_\R/\spn(\tau)$, let 
$\pi_\tau:N_\R\to N_\R'$ be the projection, and let 
$v_0\in\Trop(f,N_\R(\sigma))\cap N_\R'$.  
Fix $(u_0,\val(a_{u_0}))\in\vertices_{v_0}(f)\subset H(f,\tau)$ in the notation
of~\parref{par:height.graph}, and let 
$\alpha = \val(a_{u_0})-\angles{u_0,v_0}$.  Suppose for the moment that there
exists $v_1\in N_\R$ with $\pi_\tau(v_1)=v_0$ such that 
$\val(a_u)-\angles{u,v_1} > \alpha$ for all 
$u\in S_\sigma\setminus\tau^\perp$ (note
$\val(a_u)-\angles{u,v_1}=\val(a_u)-\angles{u,v_0}$ for 
$u\in\tau^\perp$).  Since $\angles{u,v} \leq 0$ for all
$u\in S_\sigma$ and $v\in\tau$, we therefore have 
$\val(a_u)-\angles{u,v_1+v} > \alpha$ when 
$u\in S_\sigma\setminus\tau^\perp$, and hence by~\eqref{eq:cuts.out.Trop}, 
$\vertices_{v_1+v}(f) = \vertices_{v_0}(f)$
for all $v\in\tau$.  But $\inn_{v_0}(f)$ is not a
monomial, so $\inn_{v_1+v}(f)$ 
is not a monomial, and hence $v_1+v\in\Trop(f,N_\R)$ for all 
$v\in\tau$ by~\parref{lem:trop.hypersurface}.  But
$v_1 + v\to v_0$ as $v\to\infty$, so
$\Trop(f,N_\R(\sigma))\subset\bar{\Trop(f,N_\R)}$. 

It remains to prove the existence of such an element $v_1$.  
Choose any $v_1\in\pi_\tau\inv(v_0)$.
For $u\in S_\sigma\setminus\tau^\perp$ we can find
$v\in\tau$ such that $\angles{u,v} < 0$; replacing
$v_1$ with $v_1 + \lambda v$ for $\lambda\gg 0$ allows us to assume that 
$\val(a_u) - \angles{u,v} > \alpha$.  Repeating this procedure for the
finitely many $u\in S_\sigma\setminus\tau^\perp$ for which $a_u\neq 0$
provides the required element $v_1$.  This completes the proof of~(ii).

Since $N_\R(\Delta)$ is covered by the open subspaces $N_\R(\sigma)$ for
$\sigma\in\Delta$, we will prove~(iii) with $N_\R(\sigma)$
replacing $N_\R(\Delta)$.  Let $\tau\prec\sigma$ be nonzero and define
$N_\R'$ and $\pi_\tau$ as above.  The inclusion
$\bar{\bigcap_{i=1}^n\Trop(f_i)}\subset\bigcap_{i=1}^n\bar{\Trop(f_i)}$ is
automatic, so let $v_0\in\bar{\bigcap_{i=1}^n\Trop(f_i)}\cap N_\R'$.  In
the proof of~(ii) we showed that there exists
$v_i\in\pi_\tau\inv(v_0)$ such that $v_i+\tau\subset\Trop(f_i)$ for each
$i=1,\ldots,n$.  Since $\tau$ spans $\ker(\pi_\tau)$, there is some element
$v\in\bigcap_{i=1}^n (v_i+\tau)$; then
$v+\tau\subset\bigcap_{i=1}^n\Trop(f_i)$, and $v+w\to v_0$ as
$w\to\infty$, which shows that $v_0\in\bar{\bigcap_{i=1}^n\Trop(f_i)}$.
This proves~(iii).

The final assertion follows immediately from~(ii) and~(iii),
since $\Trop(Y)\subset\bigcap_{i=1}^n\Trop(Y_i)$.\qed

\begin{lem} In \label{lem:dominator}
  the situation of~\parref{lem:moving}, suppose that
  $\Delta = \bigcap_{i=1}^n \Delta(f_i)$ is pointed.  Then
  there is a thickening
  $\cP'$ of the complex underlying $C$, dominated by $\cP$, such that
  $|\cP|\cap\bigcap_{i=1}^n(\Trop(f_i)+t v_i)\subset|\cP'|$ for all
  $t\in[0,\epsilon]$.
\end{lem}

\pf For $X\subset N_\R$ let $\bar X$ denote its closure in $N_\R(\Delta)$.  Let
$\Pi$ denote the polyhedral complex underlying $C$.  For any $P\in\Pi$ we
have $\cU(P)\in\Delta$ by~(\ref{lem:Delta.unbounded},ii), and hence
$\bar P\subset N_\R(\Delta)$.  If $P'$ is a thickening of $P$ then
$\bar P{}'\subset N_\R(\Delta)$ as well, and the interior of
$\bar P{}'$ is an increasing union of the closures of smaller thickenings
$P''$ of $P$.  Hence we can write
$\bigcup_{P'\in \cP}(\bar P{}')^\circ$ as an increasing union 
$\bigcup_{i=1}^\infty \bar{|\cP'_i|}$,
where each $\cP'_i$ is dominated by $\cP$.  Consider the set
\[ D = \bar{|\cP|}\cap\bigcup_{t\in[0,\epsilon]}
\bigcap_{j=1}^n\bar{(\Trop(f_i)+tv_i)}. \]
By~(\ref{prop:trop.closure},iii) when $t\in(0,\epsilon]$ we have
\[ \bigcap_{j=1}^n\bar{(\Trop(f_i)+tv_i)} 
= \bar{\bigcap_{j=1}^n(\Trop(f_i)+tv_i)}
= \bigcap_{j=1}^n(\Trop(f_i)+tv_i)\subset\mathring \cP \]
since the right side is a finite set of points contained in 
$\mathring \cP$, and clearly $\bigcap_{j=1}^n\bar{\Trop(f_i)}$ is contained
in $\bigcup_{P'\in \cP}(\bar P{}')^\circ$.  Hence $D$ is covered by
$\bigcup_{i=1}^\infty\bar{|\cP_i'|}$, so it suffices to show that
$D$ is \emph{compact}.

For $i=1,\ldots,n$ let
\[ D_i' = \bigcup_{t\in[0,\epsilon]} \{t\}\times
\bar{(\Trop(f_i)+tv_i)}
\subset [0,\epsilon]\times N_\R(\Delta), \]
so $D$ is the image of 
$([0,\epsilon]\times\bar{|\cP|})\cap\bigcap_{i=1}^n D_i'$
under the projection  
$[0,\epsilon]\times N_\R(\Delta)\to N_\R(\Delta)$.  Since 
$[0,\epsilon]\times N_\R(\Delta)$ is compact, it is enough to show that  
each $D_i'$ is closed.  But this is clear because 
$\bar{(\Trop(f_i)+tv_i)} = \bar{\Trop(f_i)} + tv_i$.\qed

Finally we note that in the case we will be interested in, schematic
closure respects fibers:

\begin{lem} Let \label{lem:closure.fibers}
  $f = \sum_{u\in M} a_u x^u\in K[M]$ be nonzero and let
  $\Delta$ be an integral pointed fan refining $\Delta(f)$.  Let $v\in N$,
  and define
  \[ g = \sum_{u\in M} a_u x^u t^{\angles{u,v}}\in K[M][t^{\pm 1}]. \]
  Let $Y\subset X(\Delta)\times\G_m$ be the closure of $V(g)$, let
  $\pi: Y\to\G_m$ be projection onto the second factor, and for
  $t_0\in|\G_m|$ let $Y_{t_0}=\pi\inv(t_0)$ and let
  $g_{t_0}$ be the specialization of $g$ at $t_0$.  Then $Y_{t_0}$ is the
  closure of $V(g_{t_0})$.
\end{lem}

\pf Fix $\sigma\in\Delta$.  By~(\ref{prop:trop.closure},i) there exists
$u_1\in M$ such that $Y\cap(X(\sigma)\times\G_m)$ is cut out by
$x^{-u_1} g\in K[S_\sigma][t^{\pm 1}]$ and such that the closure of
$V(g_{t_0})$ 
is cut out by $x^{-u_1} g_{t_0}$ (since $|\New(g_{t_0})|=|\New(g)|$).
But $Y_{t_0}\cap(X(\sigma)_{\kappa(t_0)}\times\G_m)$ is also cut out by
$x^{-u_1}g_{t_0}$.\qed

\pf[of~\parref{thm:non.proper.multiplicity} and~\parref{cor:well.defined}]
Let $\cP$ be a thickening of $C$ such that 
$|\cP|\cap\bigcap_{i=1}^n\Trop(f_i) = C$, let
$v_1,\ldots,v_n\in N$ and $\epsilon\in\R_{>0}\cap\Gamma$ be as
in~\parref{lem:moving}, and let $\cP'$ be as in~\parref{lem:dominator}.
We may assume without loss of generality that all polyhedra in question
are integral $\Gamma$-affine (and pointed).  Writing
$f_i = \sum_{u\in M} a_{i,u} x^u$, define
$g_i\in K[M][t^{\pm 1}]$ by
$g_i = \sum_{u\in M} a_{i,u} x^u t^{\angles{u,v_i}}$.
Let $Y_i\subset X(\Delta)\times\G_m$ be the closure of $V(g_i)$ and let
$Y = \bigcap_{i=1}^n Y_i$.  For $t_0\in|\G_m|$ let
$g_{i,t_0}$ be the specialization at $t_0$, let
$Y_{i,t_0}$ be the fiber of $Y_i$ over $t_0$, and let
$Y_{t_0} = \bigcap_{i=1}^n Y_{i,t_0}$.  By~\parref{lem:closure.fibers}
$Y_{i,t_0}$ is the closure of $V(g_{i,t_0})$.
If $\delta = -\val(t_0)$ then $\Trop(g_{i,t_0}) = \Trop(f_i)+\delta v_i$
by~\parref{par:torus.translation}. 

Let $S = \Sp(A)\subset\B^1_K$ be the annulus~\parref{eg:laurent.domains}
with inner radius $\exp(-\epsilon)$ and outer radius $1$, let
$Y_\cP = Y^\an\cap(U_\cP\times S)$ in the notation
of~\parref{par:controots1.setup}, and for $t_0\in|S|$ let
\[ Y_{\cP,t_0} = Y_{t_0}^\an\cap(U_\cP\times\{t_0\})
= Y^\an\cap(U_\cP\times\{t_0\}). \]
By~(\ref{prop:trop.closure},iv) we have that 
$\Trop(Y_{\cP,t_0})$ is contained in the closure of 
$|\cP|\cap\bigcap_{i=1}^n(\Trop(f_i)+\delta v_i)$ in $N_\R(\Delta)$.  When
$t_0=1$ this implies that $\Trop(Y_{\cP,1})\subset\bar C$, and 
$\Trop(Y_{\cP,t_0})\subset\bar|\cP'|$ when $\delta\in(0,\epsilon]$.
Therefore the hypotheses  
of~\parref{thm:controots1} are satisfied, so any two finite fibers
$Y_{\cP,t_0}$ have the same length.  By hypothesis
$Y_{\cP,1}$ is finite, and by~\parref{thm:multiplicity.ps}, the length of
$Y_{\cP,t_0}$ is equal to $i(C,\,\Trop(f_1)\cdots\Trop(f_n))$ when
$\delta\in(0,\epsilon]$.

The corollary is proved as follows.  
If $Y_{\cP,1}$ is finite then we are done, so suppose $Y_{\cP,1}$ is not
finite.  When
$\delta\in(0,\epsilon]$ the fiber $Y_{\cP,t_0}$ is still finite, so by
semicontinuity of fiber dimension there
exists some $t_1$ with $\val(t_1)=0$ such that $Y_{\cP,t_1}$ is
finite; we then apply~\parref{thm:non.proper.multiplicity} to
$g_{1,t_1},\ldots,g_{n,t_1}$ (note that $\Trop(g_{i,t_1})=\Trop(f_i)$ and
$\New(g_{i,t_1})=\New(f_i)$).\qed

\bibliographystyle{thesis}
\bibliography{thesis}

\end{document}